\documentclass[10pt,a4paper]{article}

\usepackage{geometry}
\usepackage[utf8]{inputenc}
\usepackage{multirow}
\usepackage{amsthm, amssymb, natbib, graphicx}
\usepackage{listings}
\usepackage[ruled,vlined]{algorithm2e}
\usepackage{caption}
\usepackage{amsmath,tikz}
\usepackage{subcaption}
\usepackage[normalem]{ulem}
\usetikzlibrary{positioning}
\usetikzlibrary{arrows.meta,shapes.geometric}
\usepackage{array, makecell}
\usepackage{sidecap, caption}
\usepackage{xcolor, color} 
\usepackage{comment} 
\usepackage{algorithm2e}
\usepackage[most]{tcolorbox}
\setcitestyle{numbers, comma,square}
\usepackage{bbm}
\usepackage{caption}
\usepackage{tabularray}
\usepackage{diagbox}
\UseTblrLibrary{booktabs}
\usepackage{csquotes}
\usepackage[
    pdfpagelabels,
    hypertexnames=false,  % Make link-names unique -- even in corner cases of subequations
    bookmarksopen=true,
    bookmarksopenlevel=1,  % PDF bookmarks are opened to the first level by default
    colorlinks=true,
    pdfstartview=FitV,
    linkcolor=blue,
    citecolor=blue,
    urlcolor=blue,
]{hyperref} % Needs to be loaded early, so it is included here
\hypersetup{
    colorlinks=true,
    linkcolor=blue,
    filecolor=magenta,      
    urlcolor=cyan,
    pdftitle={Overleaf Example},
    pdfpagemode=FullScreen,
    }
\usepackage{cleveref}

\newtheorem{theorem}{Theorem}[section]
\newtheorem{corollary}{Corollary}[section]
\newtheorem{definition}{Definition}[section]
\newtheorem{remark}{Remark}[section]
\newtheorem{lemma}{Lemma}[section]

\newtheorem{assumption}{Assumption}

\DeclareMathOperator*{\argmin}{argmin}

\DeclareMathOperator{\diag}{diag}

\DeclareMathOperator{\sign}{sign}

\newcommand{\R}{\mathbb{R}}
\def\s{{\tilde{s}}}
\def\H{{\tilde{H}}}
\def\g{{\tilde{g}}}

\def\W{{\tilde{W}}}
\def\m{{\tilde{m}}}

\usepackage[normalem]{ulem}
\newcommand{\rrrr}[1]{}

\newcommand{\bbbb}[1]{}
\newcommand{\ignore}[1]{}
\def\T{{\cal T}}
\def\B_d{{\cal B}}
\def\G_d{{\cal G}}
\def\K{{\cal K}}

\title{Global Optimality Characterizations and Algorithms for Minimizing Quartically-Regularized Third-Order Taylor Polynomials}
\author{\thanks{wenqi.zhu@maths.ox.ac.uk, University of Oxford, UK} Wenqi Zhu,  \thanks{cartis@maths.ox.ac.uk, University of Oxford, UK} Coralia Cartis}

\author{Coralia Cartis\thanks{The order of the authors is alphabetical; the second author (Wenqi Zhu) is the primary contributor.} \textsuperscript{\normalfont ,}\thanks{Mathematical Institute, Woodstock Road, University of Oxford, Oxford, UK, OX2 6GG.  {\tt cartis@maths.ox.ac.uk} } \quad and \quad
Wenqi Zhu\footnotemark[1] \textsuperscript{\normalfont,}\thanks{Mathematical Institute, Woodstock Road, University of Oxford, Oxford, UK, OX2 6GG.
{\tt wenqi.zhu@maths.ox.ac.uk}}}

\date{Submitted on 2025/4/28; Revision on 2026/7/22}

\begin{document}
\maketitle

\begin{abstract}
 High-order methods for convex and nonconvex optimization, particularly $p$th-order Adaptive Regularization Methods (AR$p$), have attracted significant research interest by naturally incorporating high-order Taylor models into adaptive regularization frameworks, resulting in algorithms with better {worst-case global evaluation complexity} and faster local convergence rates than first- and second-order methods. 
However, for $p \ge 3$, characterizing global optimality and efficiently minimizing high-order Taylor models with regularization terms remains an open question. 
This paper establishes global optimality conditions for general, nonconvex cubic polynomials with quartic regularization. 
These criteria generalise existing results, recovering the optimality results for regularized quadratic polynomials,  and can be further simplified in the low-rank and diagonal tensor cases. Under suitable assumptions on the Taylor polynomial, we derive a lower bound for the regularization parameter such that the necessary and sufficient criteria coincide, establishing a connection between this bound and the subproblem's convexification and sum-of-squares (SoS) convexification techniques.
{Leveraging the optimality characterization,} we develop a Diagonal Tensor Method (DTM)  for minimizing quartically-regularized cubic Taylor polynomials by iteratively minimizing a sequence of local models that incorporate both diagonal cubic terms and quartic regularization  (DTM model). We {prove} that the DTM algorithm is convergent, with a global evaluation complexity of $\mathcal{O}(\epsilon^{-3/2})$. Furthermore, when a special structure is present (such as low rank or diagonal), DTM can exactly solve the given problem (in one iteration). 
In our numerical experiments, we propose practical DTM variants that exploit local problem information for model construction, which we then show to be competitive with cubic regularization and other subproblem solvers, with superior performance on problems with special structure.
\end{abstract}

\section{Introduction and Problem Set-up}
\label{sec intro}
{Recent works show  that $p$th-order Adaptive Regularization Methods (AR$p$, with $p \geq 3$), 
which incorporate higher-order derivatives into adaptive regularization frameworks, 
achieve faster global and local convergence rates than first- and second-order methods, both in theory and in numerical experiments~\cite{birgin2017worst, cartis2020sharp, cartis2020concise, doikov2022local}. 
The minimization of regularized polynomials lies at the core of regularization-based high-order optimization methods.}
Given the unconstrained nonconvex optimization problem
$$
  \min_{x\in \mathbb{R}^n} f(x),
$$
where $f: \mathbb{R}^n \rightarrow \mathbb{R}$ is $p$-times continuously differentiable and bounded below, with $p \ge 1$.  
The AR$p$ framework is constructed as such, use a local $(p+1)$-degree polynomial model
\begin{eqnarray}
 \min_{s \in \mathbb{R}^n} m_p(x_k, s) 
 = T_p(x_k, s) + \frac{\sigma_k}{p+1}\|s\|^{p+1},
    \label{subprob p}
\end{eqnarray}
to approximate $f(x_k+s)$, where $T_p(x_k, s)$ is the $p$th-order Taylor expansion of $f(x_k+s)$ at $x_k$. 
Here $\sigma_k > 0$ is a regularization parameter that ensures the model is bounded below; it is adjusted adaptively to ensure progress toward optimality across iterations.  
At each iteration, $x_k$ is updated via  
$
s_k \approx \arg\min_{s \in \mathbb{R}^n} m_p(x_k, s),
$$
x_{k+1} := x_k + s_k,
$
whenever sufficient objective decrease is achieved.  
This process continues until an approximate local minimizer of $f$ is obtained, satisfying $\|\nabla f(x_k)\| \le \epsilon_g$.
% {In AR$p$ methods, a polynomial local model $ m_p(x, s) $ is typically constructed to approximate $ f(x+s) $ at the iterate $ x = x_k $. Then,}
% . {The model $m_p$ is based on the $p$th-order Taylor expansion of $f(x_k+s)$ at $x_k$, $T_p(x_k, s) $ which itself is a $p$th order polynomial. 
% To ensure the local model is bounded below, an adaptive regularization term is added to $ T_p $, leading to  $    m_p(x_k, s) = T_p(x_k, s)+ \frac{\sigma_k}{p+1}\|s\|^{p+1} $ where $\sigma_k > 0$. 
% Assuming Lipschitz continuity on $\nabla^p f$,}
It is shown in \cite{carmon2020lower, cartis2022evaluation} that the AR$p$ algorithm needs at most 
$\mathcal{O}(\epsilon_g^{-\frac{p+1}{p}})$ evaluations of $f$ and its derivatives to obtain a point with 
$\|\nabla f(x_k)\| \le \epsilon_g$, and this bound is optimal for this function class. This result demonstrates that as we increase the order $ p $, the evaluation complexity bound improves, which is a key motivation for implementing this method.

% $$
%    T_p(x_k, s) := f(x_k) + \sum_{j=1}^p \frac{1}{j!} \nabla_x^j f(x_k)[s]^j,
% $$
% where $\nabla_x^j f(x_k) \in \R^{n^j}$ is a $j$th-order tensor and $\nabla^j f(x_k)[s]^j$ is the $j$th-order derivative of $f$ at $x_k$ along $s \in \R^n$. 

When $p=1$, minimizing \eqref{subprob p} leads to variants of the steepest-descent algorithm. The case $p=2$ corresponds to  Newton-type algorithms, yielding adaptive cubic regularization (AR2/ARC) variants \cite{cartis2022evaluation}. In this paper, {we focus on the case $p = 3$ and \eqref{subprob p} becomes the minimization of nonconvex cubic polynomials with quartic regularization}
\begin{equation}
\tag{AR$3$ Model}
 m_3(x_k, s) = f(x_k) + \nabla_x f(x_k)^Ts + \frac{1}{2} \nabla_x^2 f(x_k) [s]^2 + \frac{1}{6}  \nabla_x^3 f(x_k) [s]^3+ \frac{ \sigma_k}{4} \|s\|_W^4
  \label{ar3 model}
\end{equation} 
where $\nabla^j f(x_k)[s]^j$ is the $j$th-order derivative of $f$ at $x_k$ along $s \in \R^n$ for $j =2,3$. $\nabla_x^2 f(x_k)$ is a symmetric matrix, $\nabla_x^3 f(x_k) \in \R^{n\times n \times n}$ is a {symmetric} third-order tensor\footnote{A {symmetric} tensor means that the entries are invariant under any permutation of its indices,} and 
$W$ is a symmetric, positive-definite matrix with norm\footnote{Unless otherwise stated, $\|\cdot\|$ denotes the Euclidean norm in this paper.} $\|v\|_W:= \sqrt{v^TWv}$. Note that typically the regularization term uses the Euclidean norm. Here, we use the $W$ norm for generality, which is particularly useful when applying a change of basis.

One of the most critical challenges for higher-order methods is the efficient minimization of the $p$th-order subproblem for $p\ge 3$. Developing effective algorithms for minimizing the AR$3$ subproblem remains an active area of research. Schnabel et al.~\cite{chow1989derivative, schnabel1991tensor} were among the first to explore third-order tensor models for unconstrained optimization without regularization terms.
Following that, for convex $m_3$, Nesterov proposed a series of methods tailored for its minimization \cite{Nesterov2021implementable, Nesterov2020inexact, Nesterov2021inexact, Nesterov2021superfast, Nesterov2022quartic, Nesterov2006cubic}. These methods leverage convex optimization techniques, including Bregman gradient approaches, high-order proximal operators, and iterative convex quadratic models with quartic regularization.
In the nonconvex setting, recent works by Cartis and Zhu introduced several efficient algorithms and proved their convergence, such as the Quadratic Quartic Regularization (QQR), Cubic Quartic Regularization (CQR), and Sums-of-Squares Taylor (SoS) methods \cite{zhu2022quartic, zhu2023cubic, cartis2023second, zhu2024global}. Subsequently, numerical implementation and adaptive regularization techniques are provided in \cite{cartis2024Efficient}.
There is also a body of literature on polynomial optimization that focuses on the (global) minimization of certain quartic polynomials, using tools like SoS techniques, SDP relaxation, and the Lasserre hierarchy \cite{ahmadi2023higher, lasserre2001global, lasserre2015introduction, laurent2009sums}.  {However, traditionally, these works were not associated with} Taylor models and adaptive regularization frameworks. Lately, \cite{ahmadi2023higher} proposed a high-order Newton's method incorporating an SoS-convex approximation of the AR$3$ subproblem, and a global convergence rate for SoS-based approximations within nonconvex adaptive regularization was established in \cite{zhu2024global}.

Despite advancements in these methods, the key question of global optimality characterization for this subproblem remains unresolved. This paper addresses this key research gap by characterizing the global optimality of the AR$3$ subproblem. 
% {(for fixed $k$), or equivalently, of multivariate quartically-regularized symmetric cubic polynomials (see \eqref{m3}). 
% Additionally, we link these necessary and sufficient conditions to SoS-convex polynomials, identifying the special cases where the necessary and sufficient criteria coincide. 
% We leverage these optimality conditions to design an efficient algorithm tailored specifically for AR$3$ subproblem minimization. }
The paper makes the following key contributions.

\begin{enumerate} 
\item To the best of our knowledge, we are the first to establish {necessary and sufficient} global optimality condition for a generally nonconvex cubic polynomial with quartic regularization.
 Unlike quadratic polynomials with regularization, in \eqref{ar3 model}, {our} necessary and sufficient global optimality conditions do not generally coincide.
Interestingly, we find that under suitable assumptions for the Taylor polynomial, when the regularization parameter \( \sigma_k \) is large enough, the necessary and sufficient criteria coincide, and global optimality can be fully characterized. We provide a lower bound for \( \sigma_k \) and discover that this bound is closely related to the convexification and SoS-convexification of the resulting polynomial.
Specifically, the size of the regularization parameter required to close the gap between the necessary and sufficient conditions differs by a constant from that required for SoS-convexification. This provides insight into how the convex and SoS-convex quartically-regularized polynomials are related. %Note that SoS convexity can be verified in polynomial time, whereas determining the convexity of quartic-regularized polynomials is strongly NP-hard \cite{ahmadi2023higher, ahmadi2013complete}. 

\item  {We show our} global optimality conditions can be seen as a generalization of optimality conditions for various forms of quadratic polynomials with simpler tensor terms or alternative regularization norms. For instance, our results recover the optimality conditions established by Cartis et al \cite[Thm 8.2.8]{cartis2022evaluation} for quadratic polynomials with quartic regularization. We also provide special cases where our global optimality conditions can be simplified, such as when the quartic regularized polynomial has a zero tensor term, a low-rank tensor term, or a diagonal tensor term.  %TUsing our global optimality conditions, we design an algorithm for minimizing the AR3 subproblem. he algorithm uses Cholesky decompositions and Newton's method for root finding.
{We developed an algorithm for minimizing the cubic quartically regularized model that applies the sufficient global optimality conditions to form a positive definite linear system and solves the resulting secular equation via Cholesky factorizations and Newton’s method.}

\item We propose a convergent algorithm, the Diagonal Tensor Method (DTM), for minimizing $m_3(s)$. The complexity for the algorithm is at least as good as ARC/AR$2$, namely, $\mathcal{O}(\epsilon^{-3/2})$, with improvements in specific cases. For example, if the third-order derivative of the objective function has {a low rank structure}, then through an appropriate change of basis, {the AR3 subproblem can be written in the form of the DTM polynomial, which has a diagonal tensor term. Therefore, minimizing the DTM subproblem recovers the minimizer of \(m_3(s)\).}
Compared to previous algorithms (such as QQR and CQR \cite{zhu2022quartic, zhu2023cubic, cartis2023second}), DTM offers greater flexibility in approximating the tensor term of $m_3$ and incorporating tensor information. Numerical experiments further indicate that DTM performs competitively with state-of-the-art algorithms in minimizing \eqref{ar3 model}.

\end{enumerate}

The paper is organized as follows. In \Cref{sec: ar3 global opt}, we establish the global optimality conditions for the AR$3$ model and identify circumstances under which the global necessary and sufficient conditions coincide, and we give special cases where the global optimality conditions can be simplified in \Cref{Sec: Optimality Conditions in Special Cases}. In \Cref{sec: Secular Equation dtm}, we present an algorithm for minimizing \eqref{ar3 model} based on the optimality conditions. %These include cases with a zero tensor term, a low-rank tensor term, or a diagonal tensor term. Additionally, we provide an algorithm for minimizing \eqref{ar3 model} based on the optimality conditions, utilizing Cholesky decompositions and Newton's method for root finding. 
In \Cref{Sec Iterative Algorithm for Minimizing AR3 Model}, we present an iterative algorithm, the Diagonal Tensor Method (DTM), which leverages these optimality conditions to efficiently minimize the AR$3$ model. We prove the convergence of the algorithm and analyze its complexity. In \Cref{sec numerics}, we introduce a practical variant of the algorithm and provide preliminary numerical results, showing improved performance for cases with low-rank tensor structures. Finally, \Cref{sec: conclusion} concludes the paper.

\section{Global Optimality Conditions for the AR3 Model}
\label{sec: ar3 global opt}

In this section, we examine  global optimality conditions for \eqref{ar3 model}. Necessary conditions are presented in \Cref{sec: Necessary}, while sufficient conditions are discussed in \Cref{sec sufficiency for CQR}. \Cref{sec: gap global Necessary and sufficient} explores the differences between the necessary and sufficient conditions for global optimality and under what circumstances this gap can be closed. In this section and generally in this paper, since we focus on AR$3$ subproblem for a single iteration with a fixed $x_k$, unless otherwise stated, for notational simplicity, we write $m_3(x_k, s)$ as
\begin{equation}
 m_3(s) = f_0 + g^Ts+ \frac{1}{2}  H [s]^2 + \frac{1}{6} T [s]^3+ \frac{\sigma}{4}  \|s\|_W^4,
 \label{m3}
\end{equation}
where $f_0 = f(x_k) =  m_3(x_k, 0) \in \R$, $g = \nabla_x f(x_k) \in \R^n$, $ H  = \nabla^2_x f(x_k) \in \R^{ n \times n}$ and $T  = \nabla^3_x f(x_k) \in \R^{ n^3}$. 
Let
\begin{eqnarray}
\B_d(s)  := H + \frac{1}{2}T [s] + \sigma  \|s\|_W^2W
, \qquad
\G_d(s) := H +  T[s] +  \sigma \|s\|_W^2 W. 
\label{B and G}
\end{eqnarray}
If $s_*$ is a second order local minimum of $m_3$, then the  second-order local optimality conditions of \eqref{ar3 model} are
\begin{eqnarray}
 \nabla  m_3(s_*) = 0   \quad &\Rightarrow&  \quad    \B_d(s_*) s_*  = -g,
  \label{local iff1 m3} 
\\
 \nabla^2   m_3(s_*) \succeq 0  \quad &\Rightarrow&  \quad 
\G_d(s_*)   + 2 \sigma (Ws_*)(Ws_*)^T \succeq 0.
  \label{local iff2 m3}
\end{eqnarray}
where $\nabla$ denote the derivative with respect to $s$ throughout this paper. We provide the definition and bounds for the tensor norms and related corollaries in \Cref{appendix def of norm}.

\subsection{Necessary Global Optimality Conditions for the AR3 Problem}
\label{sec: Necessary}

We give a technical lemma before proving the {the necessary} global optimality for the AR$3$ subproblem. Note that \Cref{technical lemma dtm} can be seen as a generalization of Lemma 4.7 in \cite{wright2006numerical}. The proof of \Cref{technical lemma dtm} uses a similar argument as that of Lemma 4.7 in \cite{wright2006numerical} and is provided in \Cref{appendix alternative lemma}.

\begin{lemma}
Let $\B_d(s)$ and  $\G_d(s)$ be defined as in \eqref{B and G}. For any vector $s, v\in \R^n$, we have 
\begin{eqnarray}
m_3(s+v)  - m_3(s) =  \big[g  + \B_d(s)s\big]^T  v +\frac{1}{2} \G_d(s)[v]^2 + \frac{1}{6} T  [v]^3 +  \frac{\sigma}{4} \bigg[ \|s+v\|_W^2 - \|s\|_W^2\bigg] ^2.
\label{universal difference}
\end{eqnarray}
\label{technical lemma dtm}
\end{lemma}

\begin{definition} (Weighted Tensor Norm)
Let $ T \in \R^{n^3}$ be a third-order tensor and $W$ be a symmetric positive-definite matrix. 
The weighted tensor norm is defined as
\begin{equation}
\Lambda_W := \max_{u,v \in \R^n} \frac{|T[u][v]^2|}{\|u\|_W \|v\|_W^2},
\end{equation}
which represents the `size' of the tensor. The asymmetric scaling in $u$ and $v$ is tailored to facilitate second-order bounds and tensor terms arising from Hessian operations\footnote{More details on $\Lambda_W$ are given in \Cref{lemma Tensor norm}.}.  $\Lambda_W$ admits an explicit upper bound (\Cref{lemma Tensor norm}; see \Cref{remark explicit LambdaW} for a closed-form bound) depending only on the tensor entries $T_{i,j,k}$ and the smallest eigenvalue of $W$, both of which are available from the problem data.
\label{def Lambda W}
\end{definition}

\begin{theorem}
\label{thm necessary tight}
\textbf{(Necessary Global Conditions)}
Let $\Lambda_W$ be defined as in \Cref{def Lambda W}. 
If $s_*$ is a global minimum of $m_3(s)$ over $\R^n$, then  $s_*$ satisfies $ \B_d(s_*) s_*:= -g$ and 
\begin{eqnarray}
H + \frac{2}{3} T [s_*] +\sigma \|s_*\|_W^2  W   + \frac{\Lambda_W}{3} \|s_*\|_W W\succeq 0 
  \label{necessary tight DTM}
\end{eqnarray}
where $\B_d$ and $\G_d$ are defined in \eqref{B and G}. 
\end{theorem}

\begin{proof}
When $s_*$ is a global minimum of $m_3(s)$ and $ \B_d(s_*) s_*:= -g$, we deduce from \eqref{universal difference} that
\begin{eqnarray}
0 \le 2 \big[m_3(s_*+v) - m_3(s_*)\big] =  \G_d(s_*) [v]^2 + \frac{1}{3}T[v]^3 +  \frac{\sigma}{2}\big(\|s_*\|_W^2-\|v+s_*\|_W^2\big)^2. 
  \label{mid result 2}
\end{eqnarray}
If $v^T W s_* = 0$, then \eqref{mid result 2} becomes
$$
0 \le \bigg[\G_d(s_*) \bigg[\frac{v}{\|v\|}\bigg]^2+ \frac{1}{3}T\bigg[\frac{v}{\|v\|}\bigg]^3+  \frac{\sigma}{2} \frac{\|v\|_W^4}{\|v\|^2} \bigg] \|v\|^2.
$$
%We choose the sign of $v$ such that $T[v]^3 <0$. while the second dominant term $T[v]^3$ is negative depending on the sign of $v$. 
As $\|v\| \rightarrow 0 $, we have that the first term $\G_d(s_*) \big[\frac{v}{\|v\|}\big]^2$ dominates, therefore $\G_d(s_*) [v]^2\ge 0$. Since $ -\frac{1}{3} T [s_*]   + \frac{\Lambda_W}{3}  [v]^2\ge 0$, we deduce that 
$
\big[H + \frac{2}{3} T [s_*] +\sigma \|s_*\|_W^2  W   + \frac{\Lambda_W}{3}\big] [v]^2\ge \G_d(s_*) [v]^2\ge 0.
$ If $\tilde{v}^T W s_* \neq 0$, since $\tilde{v}$ and $s_*$ are not orthogonal with respect to the $W$-norm, the line $s_* + \tilde{k} \tilde{v}$ intersects the ball centred at the origin with radius $\|s_*\|_W$  at two points, $s_*$ and $u$, with $  \|s_*\|_W = \|u\|_W.$ Denote $v := u-s$, and note that $v$ is parallel to $\tilde{v}$. Then, \eqref{mid result 2} can be written as
\begin{eqnarray*}
0\le 2 \big[m_3(s_*+v) - m_3(s_*)\big] &=& {v}^T \bigg[ H + \frac{2}{3}T [s_*] + \sigma  \|s_*\|_W^2 W \bigg] v + \frac{1}{3}T[u][v]^2 \notag
\\ &\le& {v}^T \bigg[ H + \frac{2}{3}T [s_*] + \sigma  W \|s_*\|_W^2 +\frac{\Lambda_W}{3} W \|u\|_W \bigg] v
\end{eqnarray*}
where the last inequality uses \Cref{lemma Tensor norm}, namely $T [u] [v]^2 \le \Lambda_W \|u\|_W \|v\|_W^2$. Finally, using $\|s_*\|_W = \|u\|_W$, we obtain \eqref{necessary tight DTM}. 
\end{proof}

\begin{remark} (Other Global Necessary Conditions) If $s_*$ is a global minimum of $m_3(s)$ over $\R^n$, then $s_*$ satisfies $ \B_d(s_*) s_*= -g $ where $\B_d(s_*)$ is defined in \eqref{B and G} and the following conditions also hold.
\begin{enumerate}
\item If $v^TWs = 0$, then \eqref{necessary tight DTM} becomes $\G_d(s_*) [v]^2\ge 0$. 
\item For all $v \in \R^n$, $\G_d(s_*) [v]^2 + \frac{1}{3} T [v]^3 \ge  0$
for $ \|s_*\|_W = \|v + s_*\|_W$. This is derived from \eqref{mid result 2}. 
\item For all $v \in \R^n$,  $ \G_d(s_*)  + \frac{2\Lambda_W}{3} \|s_*\|_W W\succeq 0$. This follows from \eqref{necessary tight DTM} since  $-\frac{1}{3} T [s_*][v^2] + \frac{\Lambda_W}{3} \|s_*\|_W \ge 0$. 
\end{enumerate}
\label{remark other necessary conditions}
\end{remark}

The proof of \Cref{thm necessary tight} (and \Cref{technical lemma dtm}) relies on reformulations using quadratic forms, drawing inspiration from Nocedal and Wright \cite[Thm 4.1]{nocedal1999numerical}. The same necessary conditions can be derived through the approach presented in Cartis et al \cite[Sec 8.2.1]{cartis2022evaluation}. 

\subsection{Sufficient Global Optimality Conditions}
\label{sec sufficiency for CQR}

\begin{theorem} \textbf{(Global Sufficient Condition)} {Let 
$
\Lambda_W 
$ be defined in \Cref{def Lambda W}.} Assume $s_* \in \R^n$ satisfies $ \B_d(s_*) s_*:= -g$ and  
\begin{eqnarray}
H +  \frac{2}{3}T [s_*]  +\sigma \|s_*\|_W^2 W  -\frac{\Lambda_W }{3} \|s_*\|_WW- \frac{\Lambda_W ^2}{18\sigma} W \succeq 0
\label{sufficient tight DTM}
\end{eqnarray}
where $\B_d$ and $\G_d$ are defined in \eqref{B and G}.  
Then, $s_*$ is a global minimum of $m_3(s)$ over $\R^n$. 
\label{thm sufficiency tight}
\end{theorem}

\begin{proof}
    From \eqref{universal difference}, when   $ \B_d(s_*) s_*:= -g.$ We have 
\begin{eqnarray*}
2\big[m_3(s_*+v)  - m_3(s_*)\big] = {v}^T \bigg[ \underbrace{H + \frac{2}{3}T [s_*] + \sigma  \|s_*\|_W^2 W }_{:=\mathcal{G}_2(s_*)}\bigg] v + \frac{1}{3}T[u][v]^2 +  {\frac{\sigma}{2}\big(\|s_*\|_W^2-\|u\|_W^2\big)^2}
\end{eqnarray*}
where $u= v+s_*$. Using $\frac{1}{3}T[u][v]^2 \ge -\frac{\Lambda_W }{3}\|u\|_W \|v\|_W ^2$, we have
\begin{eqnarray}
&&2\big[m_3(s_*+v)  - m_3(s_*)\big]\ge v^T \mathcal{G}_2(s_*)  v   -\frac{\Lambda_W }{3}  \big\|u \big\|_W   \big\|s+u \big\|_W ^2    +  {\frac{\sigma}{2}\big(\|s_*\|_W^2-\|u\|_W^2\big)^2}
\label{ineq tight 1}
\\ && \qquad =  v^T \mathcal{G}_2(s_*)  v   +  \|s_*   + u\|_W^2 \bigg[-\frac{\Lambda_W }{3}  \|u\|_W +  \frac{\sigma}{2}(\|s_*\|_W  -\|u\|_W )^2 \underbrace{\frac{ (\|s_*\|_W   + \|u\|_W)^2}{ \|s_*   + u\|_W^2}}_{\ge 1} \bigg].  \notag
\end{eqnarray}
Rearranging gives
\begin{eqnarray}
 2\big[m_3(s_*+v)  - m_3(s_*)\big] &\ge& v^T \mathcal{G}_2(s_*)  v   +  \|v\|_W^2 \bigg[-\frac{\Lambda_W }{3}  \|u\|_W +  \frac{\sigma}{2}(\|s_*\|_W  -\|u\|_W )^2 \bigg] \label{ineq tight 1.1}
\\ & =& v^T \mathcal{G}_2(s_*)  v   +  \|v\|_W^2 \bigg[  \frac{\sigma}{2}\|u\|_W ^2  - \bigg(\sigma\|s_*\|_W + \frac{\Lambda_W }{3}  \bigg)  \|u\|_W +   \frac{\sigma}{2} \|s_*\|_W^2  \bigg]   \notag
\\ & =& v^T \mathcal{G}_2(s_*)  v   +  \|v\|_W^2 \bigg[ \frac{\sigma}{2}\bigg(  \|u\|_W  - \big(\|s_*\|_W + \frac{\Lambda_W }{3 \sigma}\big) \bigg) ^2 - \frac{{\Lambda_W }^2}{18 \sigma} - \frac{\Lambda_W }{3} \|s_*\| \bigg].
\label{ineq tight 2}
\end{eqnarray}
Using $\mathcal{G}_2-\frac{{\Lambda_W }^2}{18 \sigma}W - \frac{\Lambda_W }{3 \sigma} \|s_*\| W\succeq 0$, we obtain that $m_3(s_*+v)  - m_3(s_*) \ge 0$ for all $u \in \R^n$ and therefore, $s_*$ is the global minimum. 
\end{proof}

\begin{remark} The bound in \Cref{thm sufficiency tight} is sharp and attained for some $m_3$. Assume for a certain subproblem, that the descent direction of $T$ is dominated by $s_*$, such that $T[\frac{s_*}{\|s_*\|}]^3 = -\Lambda_W$. 
Choose $u$ parallel to $s_*$ with norm $\|u\|_W  = \|s_*\|_W + \frac{t^*}{3 \sigma}$. We have $\frac{ \|s_*\|_W   + \|u\|_W}{ \|s_*   + u\|_W}=1$ and therefore,  \eqref{ineq tight 1} and \eqref{ineq tight 1.1} become equalities.
Also, since $\|u\|_W  = \|s_*\|_W + \frac{t^*}{3 \sigma}$, \eqref{ineq tight 2} becomes  
$$  
2\big[m_3(s_*+v)  - m_3(s_*)\big] = v^T \bigg[\mathcal{G}_2-\frac{{\Lambda_W }^2}{18 \sigma}W - \frac{\Lambda_W }{3} \|s_*\| W\bigg]  v. 
$$
%The first inequality uses $\|s_* +u\|_W\le \|s_*\|_W  +\|u\|_W$ with equality attained when $u$ (or equivalently $v = u+s_*$) is parallel to $s_*$.
\end{remark}

\textbf{Global Optimality and Local Optimality:} The global necessary conditions \eqref{necessary tight DTM} do not inherently imply the second-order local optimality condition \eqref{local iff2 m3}. This is because the necessary and sufficient global optimality conditions for \(m_3\) are not equivalent.
However, the global sufficient conditions \eqref{sufficient tight DTM} do lead to the second-order local optimality condition in \eqref{local iff2 m3}.  To see this, for all $v \in \R^n$,  we have
\begin{eqnarray*}
v^T\bigg[ H + T [s_*]  +\sigma \|s_*\|_W^2 +2 \sigma Ws(Ws)^T \bigg] v \underset{\eqref{sufficient tight DTM}}{\ge} \frac{1}{3} \underbrace{ [-T[s_*][v]^2 + \Lambda_W \|s_*\|_W \|v\|^2_W ]}_{\ge 0 \text{ by \Cref{lemma Tensor norm}}}+ \frac{t^*}{18\sigma} \|v\|^2_W +2 \sigma (v^TWs_*)^2 \ge 0. 
\end{eqnarray*}
Therefore, if $s_*$ satisfies {the premise of}  \Cref{thm sufficiency tight}, then $s_*$ also satisfies the second-order local optimality conditions in \eqref{local iff1 m3}--\eqref{local iff2 m3}.
%This is because the necessary and sufficient global optimality conditions are not the same for the cubic polynomial with quartic regularization. A similar behavior is also found in the CQR polynomials which have the form $ M_{\text{CQR}}(s) = f_i + g_i^T s+\frac{1}{2} H [s]^2  +  \frac{\beta}{6} \|s\|_W^3+ \frac{\sigma}{4} \|s\|_W ^4$. More details can be found \cite{zhu2023cubic}. Further details are discussed in \Cref{remark sufficient global and local condition}.

\subsection{The Gap Between the Global Necessary and Sufficient Conditions}
\label{sec: gap global Necessary and sufficient}

In this section, we explore the conditions under which the global necessary and sufficient conditions coincide for the AR$3$ subproblem. Clearly, the first-order condition \( \B_d(s_*) s_* = -g \) is identical for both the necessary and sufficient conditions. The second-order global necessary condition is 
$$ 
H + \frac{2}{3} T [s_*] + \sigma \|s_*\|_W^2 W + \frac{\Lambda_W}{3} \|s_*\|_W W \succeq 0.
$$
The second-order global sufficient condition is
$$ 
H + \frac{2}{3} T [s_*] + \sigma \|s_*\|_W^2 W - \frac{\Lambda_W}{3} \|s_*\|_W W - \frac{{\Lambda_W}^2}{18\sigma} W \succeq 0.
$$
Generally, when \(\sigma\) is large or when \(T\) is small, the gap between the global necessary and sufficient conditions narrows. 
{For instance, if \(T = 0\), then \(m_3(s)\) becomes a purely quadratic function with quartic regularization (QQR). We see that the conditions in \Cref{thm necessary tight} and \Cref{thm sufficiency tight} both become \( \B_d(s_*) s_* = -g \) and \( H + \sigma \|s_*\|_W^2 W \succeq 0 \), indicating that the global necessary and sufficient conditions coincide.} This optimality condition for the quadratic quartically regularized (QQR) polynomial is in line with that stated in \cite[Theorem 8.2.8]{cartis2022evaluation}. In the remainder of this section, we discuss the optimality conditions and convexification of $m_3$ where \(\sigma\) is large.

\subsubsection{\texorpdfstring{Locally Nonconvex $m_3$: Large $\sigma$  and Global Optimality}{Locally Nonconvex}}

We show that, under suitable assumptions, for a locally nonconvex $m_3$ (i.e., when $H$ is not positive definite), if $\sigma$ is sufficiently large, then the global necessary condition \eqref{necessary tight DTM} and the sufficient condition \eqref{sufficient tight DTM} coincide.

%Under suitable assumptions, if the first-order optimality condition of the objective function is not satisfied, namely $\|g\| \geq \epsilon_g$ for some positive tolerance $\epsilon_g$, then for a locally nonconvex $m_3$ (i.e., $H$ is not positive definite) and for $\sigma$ as large as $\mathcal{O}(\epsilon_g^2)$, the global necessary condition \eqref{necessary tight DTM} and the sufficient condition \eqref{sufficient tight DTM} coincide.
%Namely, a stationary point satisfying the global necessary condition \eqref{necessary tight DTM} is sufficient to characterize the global minimizer.

\begin{lemma}
\label{thm sigma large enough}
Assume that $s_* \in \R^n$ satisfies $ \B_d(s_*) s_*:= -g$ and $\sigma$ satisfies
\begin{eqnarray}
\sigma \ge 3 \max \bigg\{-\lambda_{\min}[H]\|s_*\|_W^{-2} \lambda_{\min}(W)^{-1}, \qquad  \frac{7}{3}\Lambda_W  \|s_*\|_W^{-1}  \bigg\}, 
   \label{sigma bound 0}
\end{eqnarray}
then \eqref{necessary tight DTM} and \eqref{sufficient tight DTM} are equivalent. Moreover, if $s_*$ and $\sigma$ satisfy $ \B_d(s_*) s_*:= -g$, \eqref{necessary tight DTM} and \eqref{sigma bound 0}, then $s_*$  is a global {minimizer}. 
\end{lemma}
\begin{proof} Clearly, \eqref{necessary tight DTM} implies \eqref{sufficient tight DTM}. To prove we converse, we rewrite \eqref{sufficient tight DTM} as 
\begin{eqnarray}
&&  \frac{2}{3} \bigg[H + \frac{2}{3} T [s_*] + \sigma \|s_*\|_W^2 W + \frac{\Lambda_W}{3} \|s_*\|_W W\bigg] + \notag
 \\ && \qquad \frac{1}{3} \bigg(\underbrace{ H + \frac{2}{3} T [s_*]  +\sigma \|s_*\|_W^2 W  -\frac{5}{3} \Lambda_W \|s_*\|_W W  - \frac{{\Lambda_W}^2}{6\sigma} W}_{:=\mathcal{J}_2} \bigg) \succeq 0.
  \label{alternative expression for suff}
\end{eqnarray}
Using \eqref{sigma bound 0}, we deduce that
\begin{eqnarray}
   \sigma  > 2^{-1/2} \Lambda_W  \|s_*\|_W^{-1} 
   \qquad &\Rightarrow& \qquad  \frac{\sigma}{3} \|s_*\|_W^2 >\frac{{\Lambda_W}^2} {6\sigma}, 
   \label{sigma bound 1}
\\ \sigma  \ge 7 \Lambda_W  \|s_*\|_W^{-1} 
\qquad &\Rightarrow& \qquad  
\frac{\sigma}{3}  \|s_*\|_W^2  \ge \frac{7}{3}\Lambda_W \|s_*\|_W ,
  \label{sigma bound 2}
\\\sigma  \ge - 3\lambda_{\min}[H]\|s_*\|_W^{-2}\lambda_{\min}(W)^{-1}
\qquad &\Rightarrow& \qquad   
\frac{\sigma}{3} \|s_*\|_W^2  \lambda_{\min}(W) \ge  -\lambda_{\min}[H].
  \label{sigma bound 3}
\end{eqnarray}
\eqref{sigma bound 2} gives $\frac{\sigma}{3} \|s_*\|_W^2 W \ge -\frac{2}{3}T[s_*][v]^2 + \frac{5}{3}\Lambda_W \|s_*\|_W$ for all unit vector $v \in \R^n$. \eqref{sigma bound 3} gives $\frac{\sigma}{3} \|s_*\|_W^2 W \succeq -H $. Consequently, we deduce that $\mathcal{J}_2$ defined in \eqref{alternative expression for suff} is positive semi-definite.
Thus, \eqref{sufficient tight DTM} implies  \eqref{necessary tight DTM}. 

To prove the last statement, if $\sigma$ and $s_*$ satisfy \eqref{sigma bound 0}, then \eqref{necessary tight DTM} and \eqref{sufficient tight DTM} are equivalent. With the condition in $ \B_d(s_*) s_*:= -g$, using \Cref{thm sufficiency tight}, we can deduce that $s_*$ is a global minimum. 
\end{proof}

Note that the size of $\sigma$ is closely related to the convexification of the AR$3$ subproblem with large $\sigma$. If $H$ is indefinite (i.e., $m_3$ is locally nonconvex at $s=0$), regardless of how large $\sigma > 0$ is chosen, it is not possible to fully convexify $m_3$ for all $s \in \R^n$. However, we can select $\sigma$ large enough such that $m_3$ is convex except in a small region around the origin. The proof is provided in \Cref{lemma Convexification nonconvex}

\begin{lemma} 
If $\lambda_{\min}[H] \leq 0$, let $s_0 > 0$ be any scalar, and suppose that $\sigma$ satisfies
\begin{eqnarray}
\sigma > 2 \max \bigg\{-\lambda_{\min}[H] \max \{s_0^{-2}, 1\} , \qquad  \Lambda_W  \max \{s_0^{-1}, 1\} \bigg\}>0,
\label{convex sigma bound DTM}
\end{eqnarray}
then $m_3(s)$ is convex for all $\|s\|_W \geq s_0$.  
\label{lemma Convexification nonconvex}
\end{lemma}

\begin{proof}
If $\lambda_{\min}[H] \le 0$, for any unit vector $v \in \R$,
\begin{eqnarray}
\nabla^2 m_3 (s) [v]^2 = \bigg[H + T[s] + \sigma \bigg(\|s\|_W W+ 2(Ws)(Ws)^T\bigg)\bigg] [v]^2 \ge  \lambda_{\min}[H]- \Lambda_W \|s\|_W +  \sigma  \|s\|_W^2.
\label{sencond order m3 temp}
\end{eqnarray}
Since $\sigma$ satisfies \eqref{convex sigma bound DTM}, for any $\|s\|_W\ge s_0$,
\begin{eqnarray}
\frac{1}{2}\sigma \|s\|_W \ge  \max \{s_0^{-1}, 1\} \Lambda_W  s_0 \ge \Lambda_W, \qquad \qquad \frac{1}{2}\sigma \|s\|_W^2 \ge \frac{1}{2}\sigma s_0^2  \ge \frac{1}{2}\min \{s_0^2, 1\} \ge -\lambda_{\min}[H].
\label{sencond order m3 temp 1}
\end{eqnarray}
Substituting \eqref{sencond order m3 temp 1} into \eqref{sencond order m3 temp}, we deduce that $\nabla^2 m_3 (s) [v]^2 \ge 0$ for all   $v \in \R$ and $\|s\|_W\ge s_0$. 
\end{proof}

The condition stated in \eqref{convex sigma bound DTM} indicates that for any arbitrarily small $s_0 > 0$, there exists a sufficiently large $\sigma$ such that $m_3$ is convex everywhere except in a small region around the origin where $\|s\|_W < s_0$.
However, the bound in \Cref{thm sigma large enough} depends on $s_*$, which is unknown \textit{a priori}. To obtain a bound that does not rely on $s_*$, similar to \cite[Lemma 3.1]{cartis2020concise}, we require a uniform lower bound for $\|s_*\|_W$. %This can be established  \cite[Lemma 3.1]{cartis2020concise}, with a minor modification to account for the $W$-norm.

\begin{theorem}
    Let $s_*$ be a first-order critical point of $m_3$ such that $\nabla m_3(x_k, s_*) = 0$. For a fixed tolerance $ \epsilon_T \in (0, 1)$, either $s_*$ is also an approximate first-order critical point of the Taylor series $T_3$, satisfying $\|\nabla T_3(x_k, s_*)\| \leq \epsilon_T$, or the following hold:
    \begin{enumerate}
        \item The following lower bound for the stepsize holds: $\|s_*\|_W \geq (\epsilon_T/\sigma)^{1/3}$.
        \item If $\sigma$ satisfies 
        \[
        \sigma \geq \max\bigg\{ (3 \lambda)^{3} \epsilon_T^{-2}, (7 \Lambda_W)^{3/2} \epsilon_T^{-1/2} \bigg\},
        \]
        where $\lambda := \max \{-\lambda_{\min}[H]\lambda_{\min}(W)^{-1}, 0\}$. Then, the necessary condition \eqref{necessary tight DTM}  and the sufficient condition \eqref{sufficient tight DTM} are equivalent.
    \end{enumerate}
    \label{thm epi T}
\end{theorem}

\begin{proof}
To prove the first result, we use $s_*$ is a first-order critical point, therefore  $\nabla m_3(x_k, s_*) = 0$ and
\begin{eqnarray*}
\epsilon_T < \big\|T_3(x_k, s_*)\big\| = \big\|T_3(x_k, s_*)- \nabla m_3(x_k, s_*)  \big\| + \underbrace{\big\| \nabla m_3(x_k, s_*)  \big\|}_{=0} = \sigma \|s_*\|_W^3 
\end{eqnarray*}
which gives the first result. 
To prove the second result, we use that $\sigma$ satisfies the condition in Result 2). Then, we have
\begin{eqnarray*}
\sigma^{1/3}\ge  3 \lambda  \epsilon_T^{-2/3}, \qquad 
 \sigma^{2/3}    \ge 7 \Lambda_W   \epsilon_T ^{-1/3}. 
\end{eqnarray*}
where $\lambda := \max \{-\lambda_{\min}[H]\lambda_{\min}(W)^{-1}, 0\}$. 
From this, we deduce that
\begin{eqnarray*}
\sigma \ge   3 \max \bigg\{ \lambda   \sigma ^{2/3} \epsilon_T^{-2/3},   \frac{7}{3}\Lambda_W   \sigma ^{1/3} \epsilon_T^{-1/3} \bigg\} \ge   3 \max \bigg\{\lambda \|s_*\|_W^{-2},   \frac{7}{3}\Lambda_W  \|s_*\|_W^{-1}  \bigg\}
\end{eqnarray*}
where the last inequality uses $\|s_*\| \ge (\epsilon_T/\sigma)^{1/3}$. 
Using \Cref{thm sigma large enough}, we obtain that  \eqref{necessary tight DTM} and \eqref{sufficient tight DTM} are equivalent.
\end{proof}

 Combining the result of partially convexifying a locally nonconvex $m_3$ (\Cref{lemma Convexification nonconvex}) with \Cref{thm epi T}, we deduce that,  if $\sigma \geq 4(-3 \lambda)^3 \epsilon_T^{-2} = \mathcal{O}(\epsilon_T^{-2})$ and $\|\nabla T_3(x_k, s_*)\| \geq \epsilon_T$,  then, we have the following three results. Firstly, $m_3$ is convex except in a small region around the origin $\|s\| < \mathcal{O}(\epsilon_T)$.  Secondly,  $\|s_*\| \geq \big[\frac{\epsilon_T}{\sigma}\big]^{1/3} = \mathcal{O}(\epsilon_T^{1/3})$, which indicates that the global minimizer lies within the convex region of $m_3$ (Result 1 of \Cref{thm epi T}). Thirdly, the necessary and sufficient conditions coincide, namely,  the global minimum of $m_3$ is characterized by $\B_d(s_*) s_* = -g$ and \eqref{necessary tight DTM}. %Additionally, by \Cref{lemma lower bound for s_*}, we have $\|s_*\| \geq \big[\frac{\epsilon_T}{\sigma}\big]^{1/3} = \mathcal{O}(\epsilon_T^{1/3})$, which ensures that the global minimum lies within the convex region of $m_3$.

As discussed in \Cref{sec intro}, the model \( m_3 \) is typically used in the AR\(3\) algorithm and is based on the third-order Taylor expansion of \( f(x_k + s) \) at \( x_k \). 
A point \( s_* \) satisfying \( \|\nabla T_3(x_k, s_*)\| \leq \epsilon_T \) is considered an approximate first-order critical point of the third-order Taylor expansion \( T_3 \). Since the Taylor series approximates \( f(x) \) in the vicinity of \( x_k \), the point \( x_k + s_* \) can be regarded as a good approximation to a first-order critical point of \( f(x) \) near \( x_k \). However, if \( s_* \) lies too far from \( x_k \), the Taylor approximation may no longer be accurate, and \( s_* \) may fail to approximate a critical point of \( f(x) \). 
To address this, the AR\(3\) framework checks at each iteration whether a sufficient decrease in the objective function has been achieved. If \( s_* \) does not yield sufficient decrease, it is rejected, and the regularization parameter \( \sigma_k \) is increased, typically leading to more conservative steps and smaller step sizes in subsequent iterations.
%In the AR\(3\) framework, the regularization parameter \( \sigma_k \) is adjusted adaptively so that a first-order local minimizer of \( f \) is found, satisfying \( \|\nabla f(x_k)\| \leq \epsilon_g \). 
A similar analysis of the size of the regularization weight, ensuring that the sufficient and necessary global optimality conditions coincide, can also be applied to the locally nonconvex subproblems in the AR\(3\) algorithm. More details can be found in \Cref{appendix Difference in Objective Function and Taylor Expansion)}.

%Let $0 < \epsilon_g \ll 1$ be the tolerance of the AR$3$ algorithm, then AR$3$ framework compute a first-order local minimizer that satisfies $\|\nabla f(x_k)\| \leq \epsilon_g$. 

\subsubsection{\texorpdfstring{Locally Convex $m_3$: Convexification  and Link to Sum-of-squares (SoS) Method}{Locally Convex}}

In this section, we demonstrate that if $m_3$ is strictly locally convex ($H \succeq \delta I_n >0$) and $\sigma$ is large enough, then $m_3$ remains convex for all $s \in \R^n$, and the global necessary and sufficient conditions coincide. We also discuss how convexification links to SoS-convex. Firstly, we provide definitions of sums of squares (SoS), and SoS-convex.

\begin{definition} \textbf{(SoS polynomial)}
    A ${\hat{p}}$-degree polynomial $q(s): \R^n \rightarrow \R$ where $s \in \R^n$ is a SoS if there exist polynomials $\Tilde{q}_1, \dotsc, \Tilde{q}_r: \R^n \rightarrow \R$, for some  $r\in \mathbb{N}$, such that $q(s) =\sum_{j=1}^r \Tilde{q}_j(s)^2$ for all $s \in \R^n$ \cite[Def. 1]{ahmadi2023higher}. 
    \label{def SoS polynomial}
\end{definition}

\noindent
\Cref{def SoS polynomial} implies that
the degree ${\hat{p}}$ of an SoS polynomial $q$ has to be even, and that the maximum degree of each $\Tilde{q}_j$ is $\frac{\hat{p}}{2}$. Let $\R[s]^{n\times n}$ be the real vector space of  $n \times n$ real polynomial matrices, where each entry of such a matrix is a polynomial with real coefficients.
%, denoted by $h_{i, j}(s): \R^n \rightarrow \R$ for all $1 \le i, j \le n$. 

\begin{definition} \textbf{(SoS-convex)}  
A polynomial $p(s): \R^n \rightarrow \R$ is SoS-convex if its Hessian $H(s):=\nabla^2 p(s)$ is an SoS matrix. An SoS matrix is a symmetric polynomial matrix
which can be expressed as $H(s) = \tilde{G}^T(s) \tilde{G}(s)$  where  \( \tilde{G}(s) \in \mathbb{R}[s]^{n \times n} \). More details can be found in \cite[Def. 2.4]{ahmadi2013complete}. 
\label{def SoS convex}
\end{definition}

SoS-convexity is a tractable sufficient condition for the convexity of polynomials through SoS decomposition \cite{ahmadi2009sum}.  Here, we prove that the bound we propose for  $\sigma$ shares the same order of magnitude as the SoS convexification bound proposed in Ahmadi et al. \cite{ahmadi2023higher, zhu2024global}. 

\begin{lemma}
\label{lemma Convexification}
(Convexification for Locally Convex $m_3$)
If $H \succeq \delta I_n \succ 0$ and $\sigma$ satisfies 
\begin{equation}
\sigma > \frac{1}{4} \Lambda_W^2\delta^{-1},
   \label{sigma bound convex}
\end{equation}
then $m_3(s)$ is convex, and any $s_*$ satisfying $\B_d(s_*)s_* = -g$ is the global minimum. 
\end{lemma}

\begin{proof}
If $H \succeq \delta I_n \succ 0$ , for any unit vector $v \in \R$,
\begin{eqnarray}
\nabla^2 m_3 (s) [v]^2 = \bigg[H + T[s] + \sigma \bigg(\|s\|_W W+ 2(Ws)(Ws)^T\bigg)\bigg] [v]^2 \ge \delta - \Lambda_W \|s\|_W +  \sigma \|s\|_W^2.
\end{eqnarray} 
The discriminant of the quadratic equation $\delta - \Lambda_W \|s\|_W +  \sigma \|s\|_W^2$ is $\Lambda_W ^2 - 4 \delta  \sigma < 0$. The condition in \eqref{sigma bound convex} ensures that the discriminant is negative. Therefore, we can ensure $\nabla^2 m_3 (s) [v]^2 > 0$ for all $s, v \in \R^n$, thus $m_3$ is convex.  
If $m_3$ is convex, any stationary point satisfying $\B_d(s_*)s_* = -g$ is the global minimum.
\end{proof}

\paragraph{Link to SoS Convexification:}  
Recently, Ahmadi et al.~\cite{ahmadi2023higher} proved that a sufficiently regularized locally convex polynomial becomes SoS-convex.  
Specifically, if $m_3$ is locally convex (i.e., $H \succeq \delta I_n \succ 0$), then for
\begin{equation}
\label{sigma bound sos dtm} 
{\sigma \ge 8 \Lambda_W^2 \delta^{-1},}
\end{equation}
$m_3(s)$ is SoS-convex for all $s \in \mathbb{R}^n$; see \cite[Lemma~3]{ahmadi2023higher} and \cite[Thm.~3.3]{zhu2024global}.  
{For such locally convex regularized polynomials, both determining a valid regularization parameter $\sigma$ and globally minimizing the resulting SoS-convex model reduce to solving semidefinite programs in polynomial time.}
The bound in \eqref{sigma bound sos dtm} is of the same order as \eqref{sigma bound convex}, with \eqref{sigma bound sos dtm} being tighter by a constant factor.  
{This is consistent with the fact that SoS-convexity implies convexity~\cite{ahmadi2009sum}.  
In this regime, the global minimizer is unique and the necessary and sufficient optimality conditions coincide.  This comparison illustrates the gap between convexity and SoS-convexity, and clarifies how increasing the regularization parameter $\sigma$ eventually forces the problem to satisfy necessary and sufficient 
global optimality conditions and to become tractable. A full discussion of this connection is beyond the scope of the current paper, and we refer the reader to~\cite{zhu2026sufficiently}.}
By examining the size of the regularization parameter, \Cref{lemma Convexification} together with \eqref{sigma bound convex} provides insight into how the class of SoS-convex quartically regularized polynomials compares with the class of convex quartically regularized polynomials.

\section{Optimality Conditions in Special Cases} 
\label{Sec: Optimality Conditions in Special Cases}

While the optimality conditions for \eqref{ar3 model} (i.e., \Cref{universal difference}, \Cref{thm necessary tight}, and \Cref{thm sufficiency tight}) represent the most general form of optimality conditions for polynomials in this class, this section discusses several special cases where the global optimality conditions \eqref{necessary tight DTM} and \eqref{sufficient tight DTM} can be significantly simplified. For certain special cases, such as when \( m_3 \) has \( T = 0 \) or the tensor term is replaced by a constant cubic term, the optimality conditions naturally connect with results in the literature. This indicates that the global optimality conditions in \eqref{necessary tight DTM}, \eqref{sufficient tight DTM}, \Cref{thm necessary tight}, and \Cref{thm sufficiency tight} represent a generalization of optimality conditions for this class of polynomials.
In \Cref{sec: low-rank}, we demonstrate that when \eqref{ar3 model} involves a low-rank tensor term or a diagonal tensor term, the constant \( \Lambda_W \) in the optimality conditions can be explicitly expressed. Additionally, we prove that \eqref{ar3 model} with a low-rank tensor term (with rank less than $n$) is equivalent to one with a diagonal tensor term through an appropriate change of basis.
%Additionally, we prove that an \eqref{ar3 model} with a rank-\( P \) tensor term (with $P \le n$) can be transformed into a Diagonal Tensor Model through an appropriate change of basis. 

%They can be simplified to the optimality conditions of various polynomials with either a simpler tensor term $T$ or with alternative regularization norms.

\textit{Quadratic polynomials with regularization} are crucial components of the (adaptive) cubic regularization framework. By setting \(T = 0\), \(m_3\) simplifies to a quadratic polynomial with quartic regularization,
\begin{equation}
\tag{QQR}
    m_{\textbf{QQR}}(s) = f_0 + g^T s + \frac{1}{2} H[s]^2 + \frac{\sigma}{4} \|s\|_W^4.
\end{equation}
Polynomials of this form have been extensively studied in \cite{zhu2022quartic, cartis2023second, cartis2022evaluation}. For these polynomials, since \(T = 0\), we have \(\B_d(s) = \G_d(s)\), and \Cref{universal difference} reduces to
\[
    m_{\textbf{QQR}}(s + v) -     m_{\textbf{QQR}}(s) = \big[g + \B_d(s) s\big]^T v + \frac{1}{2} \B_d(s)[v]^2 + \frac{\sigma}{4} \bigg[\|s + v\|^2 - \|s\|^2\bigg]^2.
\]
As a result, the optimality conditions in \Cref{thm necessary tight} and \Cref{thm sufficiency tight} simplify to \(g + \B_d(s_*) s_* = 0\) and \(\B_d(s_*) \succeq 0\), which align with Theorem 8.2.8 in \cite{cartis2022evaluation}.
    
Also, by replacing the term \(T[s]^3\) with \(\beta \|s\|^3\) where $\beta \in \R$, we obtain
\begin{equation}
\tag{CQR}
\label{CQR poly}
m_{\textbf{CQR}}(s) = f_0 + g^T s + \frac{1}{2} H[s]^2+ \beta \|s\|^3 + \frac{\sigma}{4} \|s\|_W^4.
\end{equation}
These polynomials are useful for approximating AR$p$ subproblems and are studied in detail in \cite{zhu2023cubic}. For these polynomials, since \(T[s]^3\) is replaced by \(\beta \|s\|^3\), \Cref{universal difference} becomes
\begin{eqnarray*}
m_{\textbf{CQR}}(s + v) - m_{\textbf{CQR}}(s) = \big[g + \B_d(s) s\big]^T v + \frac{1}{2} \B_d(s)[v]^2 + \frac{\sigma}{4} \bigg[\|s + v\|^2 - \|s\|^2\bigg]^2+ \\ + \frac{\beta}{12} \bigg(\|s + v\| - \|s\|\bigg)^2 \bigg(2\|s + v\|^2 + \|s\|\bigg).
\end{eqnarray*}
The global sufficiency condition in \Cref{thm sufficiency tight} becomes \(g + \B_d(s_*) s_* = 0\), \(\B_d(s_*) \succeq 0\), and \(\beta \geq -\frac{3}{2} \sigma \|s_*\|_W\), which coincides with the conditions given in \cite[Thm 2.2]{zhu2023cubic}.

The optimality conditions can also be extended to polynomials with a diagonal third-order tensor term and a\textit{ separable quartic regularization (SQR)} component. Let $s = [s_1, \dotsc, s_n]^T \in \R^n$, we express these polynomials as
\begin{equation}
\tag{SQR}
m_{\textbf{SQR}}(s) = f_0 + g^T s + \frac{1}{2} H[s]^2+  \frac{1}{6} \sum_{j=1}^n t_j s_j^3+ \frac{1}{4} \sum_{j=1}^n \sigma_j s_j^4.
\label{separable Polynomial}  
\end{equation}
where $\{t_j\}_{1 \le j \le n}$ are scalars and  $\{\sigma_j\}_{1 \le j \le n}$ are positive scalars. Let $\T = \diag\{t_j\}_{1 \le j \le n} \in \R^{n \times n \times n}$,
\begin{eqnarray}
\hat{\B_d}(s) := H + \frac{1}{2} \T [s] + \Sigma [s]^2, \qquad \text{and}  \qquad  \hat{\G_d}(s) :=H + \frac{1}{3}\T [s] + \Sigma [s]^2-\diag\big\{\frac{t_j^2}{18 \sigma_j} \big\}
\label{hat B and GS}
\end{eqnarray}
where  $\Sigma [s]^2 := \diag\{\sigma_j s_j^2\}_{1 \le j \le n} \in \R^{n \times n}$. Then \Cref{universal difference} becomes
\begin{eqnarray}
m_{\textbf{SQR}}(s+v)  - m_{\textbf{SQR}}(s) =  \big[g  + \hat{\B_d}(s)s\big]^T  v +\frac{1}{2}\hat{\G_d}(s)[v]^2 +  \frac{1}{4}\bigg\|\boldsymbol{\sigma}^{\frac{1}{2}} \cdot v \cdot \big( v +2s+\frac{t \cdot \boldsymbol{\sigma} ^{-1}}{3}\big)\bigg\|^2. 
\label{universal difference seperable}
\end{eqnarray}
where $\boldsymbol{\sigma} = [\sigma_1, \dotsc, \sigma_n]^T \in \R^n$ and $t = [t_1 , \dotsc, t_n]^T \in \R^n$. 
The global sufficiency condition in \Cref{thm sufficiency tight} becomes  
\begin{eqnarray}
\hat{\B_d}(s_*)s_* =-g, \qquad \hat{\G_d}(s_*) \succeq 0.
\label{sufficient seperable}
\end{eqnarray}
The global necessary condition in \Cref{thm necessary tight} becomes 
\begin{eqnarray}
\hat{\B_d}(s_*)s_* =-g, \qquad  \hat{\G_d}(s_*) \succeq - 2 \|H_0\|_1  I_n
\label{necessary seperable}
\end{eqnarray}
where  $\|H_0\|_1 =  \sum_{\iota, j=1, \iota \neq j}^n |H_{\iota, j}| $ is the sum of the absolute value of the off-diagonal entries of $H$. For \eqref{separable Polynomial} polynomials, if the second order term $H$ is a diagonal matrix,  $\|H_0\|_1 = 0$. The global necessary condition \eqref{necessary seperable} and sufficient condition  \eqref{sufficient seperable}  coincide. Moreover, for univariate case with $n=1$, $m_{\textbf{SQR}}(s) = m_3(s)$, and the global optimality condition for univariate $m_3(s)$ (or equivalently  $m_{\textbf{SQR}}(s)$ with $n=1$) becomes 
$\nabla m_3(s_*) = 0$ and $H + \frac{1}{3}T[s_*] + \sigma \|s_*\|^2 \ge \frac{T^2}{18 \sigma}.$ Proofs and more details can be found in \Cref{sec: CQR separable Global Opt}. Despite the necessary and sufficient conditions not being equivalent in the general case, conditions such as \eqref{necessary tight DTM}, \eqref{sufficient tight DTM}, \eqref{sufficient seperable}, and \eqref{necessary seperable} are still useful algorithmically. The algorithm implementation can be found in \Cref{sec: Secular Equation dtm}. 

%{(KZ: In fact, \( m_{\textbf{SQR}}(v) - m_{\textbf{SQR}}(s_*) \) is also a sum of squares (SoS). This can be partly observed from the proof in \Cref{sec: CQR separable Global Opt}.)}
%Let me know if we should include the optimality condition for \( m_{\textbf{SQR}}(v) \) and \Cref{sec: CQR separable Global Opt} here, or if we should move it to the next paper.}

\subsection{Transforming a Low-rank Quartically Regularized Polynomial into  a Diagonal Tensor Polynomial}
\label{sec: low-rank}

%This section proves that a \eqref{ar3 model} with a rank-$P$ tensor term can be transformed into a Diagonal Tensor Model under a suitable change of basis. This property makes the Diagonal Tensor Model useful for minimizing cubic quartically regularized polynomials with a rank $P$ tensor term. 
%The subsection is arranged as follows. In \Cref{thm change of basis dtm}, we give more details about the change of basis. In \Cref{subsec rank one rank p tensor}, we investigate the objective functions that possess a rank one (rank $P$) tensor structure. Additionally, we conducted a literature review on methods for obtaining rank-one approximations for third-order derivative tensor terms.

In this subsection, we prove that if the tensor term in \eqref{ar3 model}, denoted as $T^{\mathbbm{P}}$, has rank $P$ with $0 \le P \le n$, then by applying an appropriate change of basis, the model 
\begin{equation}
m_3^{\mathbbm{P}}(s) := f_0 + g^T s + \frac{1}{2} H [s]^2 +\frac{1}{6}T^{\mathbbm{P}}[s]^3 + \frac{\sigma}{4}\|s\|_W^4
\label{rank P m dtm}
\end{equation}
can be transformed into a quartically regularized polynomial with a diagonal tensor term  (i.e. Diagonal Tensor Model),
\begin{equation}
m_3^{\mathbbm{D}}(s) := f_0 + \g^T  \s  + \frac{1}{2} \H [\s]^2  +\frac{1}{6} \T[\s]^3  + \frac{\sigma}{4}\big\|\s \big\|_{\W}^4. 
\label{m3 diagonal}
\end{equation}
Note that $ \g$, $\H $, $\s$ are derived from  $ g$, $H $, $s$ after a change of basis and $\T$ is a diagonal tensor with $P$ non zero entries. Note that $T^{\mathbbm{P}}$ could have a dense structure while $\T$ has a sparse structure that only has $P$ nonzero terms on the diagonal. 
The special diagonal structure of the tensor in \eqref{m3 diagonal} allows the optimality conditions to be simplified and provides an explicit expression for $\Lambda_W$. According to \Cref{example special T}, for $\T = \diag\{t_j\}_{1 \le j \le n}$, we have 
\begin{equation}
\Lambda_W = \hat{t} = \max \{|t_1|, |t_2|, \dotsc, |t_n|\} [ \lambda_{\min}(W)]^{-3/2}.
\end{equation}
The global necessary conditions in \eqref{necessary tight DTM} and the sufficient conditions in \eqref{sufficient tight DTM} can be simplified by replacing $\Lambda_W$ with $\hat{t}$.

To show that \eqref{ar3 model} with a low-rank tensor term (with rank less than $n$) is equivalent to one with a diagonal tensor term through an appropriate change of basis, we begin by defining the rank $P$ tensor, following the conventions in \cite{schnabel1971tensor, schnabel1984tensor, schnabel1991tensor}.
\begin{definition}
A tensor $T^{\mathbbm{P}}$ is said to have rank $P$ if there exist $P$ linearly independent vectors $\boldsymbol{a}^{(1)}, \dotsc, \boldsymbol{a}^{(P)}$ such that
\begin{equation}
    T^{\mathbbm{P}} = \sum_{k=1}^P \boldsymbol{a}^{(k)} \otimes \boldsymbol{a}^{(k)} \otimes \boldsymbol{a}^{(k)}, \qquad \text{and} \qquad T^{\mathbbm{P}}[s]^3 = \sum_{k=1}^P \big(a^{(k)}_1 s_1 + \dots + a^{(k)}_n s_n\big)^3.
    \label{rank P dtm}
\end{equation}
Here, $\otimes$ denotes the outer product, $\boldsymbol{a}^{(k)} = [a^{(k)}_1, \dots, a^{(k)}_n]^T \in \mathbb{R}^n$, and $s = [s_1, s_2, \dots, s_n]^T \in \mathbb{R}^n$. The $(l, j, \kappa)$-th entry of the tensor is given by $(T^{\mathbbm{P}})_{lj\kappa} = \sum_{k=1}^P a^{(k)}_l a^{(k)}_j a^{(k)}_\kappa$. 
\label{tensor approx dtm}
\end{definition}

Under a suitable change of basis, a rank $P$ third-order term can be diagonalized to have $P$ non-zero diagonal entries, while a rank-one third-order term reduces to a single non-zero element.
\begin{eqnarray*}
   && \text{Rank $1$ Tensor:}\qquad   T^{\mathbbm{1}}[s]^3 =(a^{(1)}_1 s_1 +, \dotsc, + a^{(1)}_n s_n)^3 \qquad \underset{\text{change of basis}}{\leadsto} \qquad \tilde{s}_1^3
    \\
  &&  \text{Rank $P$ Tensor:}\qquad  T^{\mathbbm{P}}[s]^3 =\sum_{k=1}^P (a^{(k)}_1 s_1 +, \dotsc, + a^{(k)}_n s_n)^3 \qquad \underset{\text{change of basis}}{\leadsto} \qquad \tilde{s}_1^3 +\tilde{s}_2^3 +\dotsc+ \tilde{s}_P^3 .
\end{eqnarray*}
We provide the result in \Cref{thm change of basis dtm}.

\begin{theorem}
 \label{thm change of basis dtm}
Let $m_3^{\mathbbm{P}}(s)$ be defined as in \eqref{rank P m dtm}, and let $T^{\mathbbm{P}}$ be defined as in \Cref{tensor approx dtm} with $P \le n$. Define
\begin{align}
\hat{C} := \begin{bmatrix} \boldsymbol{a}^{(1)}, & \dotsc, & \boldsymbol{a}^{(P)} \end{bmatrix}^T \in \R^{P \times n}, \qquad C := \begin{bmatrix} \hat{C}, & \hat{C}_\perp \end{bmatrix}^T \in \R^{n \times n},
\label{C matrix dtm}
\end{align}
where {$\hat{C}_\perp = [\boldsymbol{a}_\perp^{(P+1)}, \dotsc , \boldsymbol{a}_\perp^{(n)}]^T$ with ${\boldsymbol{a}^T}^{(i)} \boldsymbol{a}_\perp^{(j)} = 0$ for all $i = 1, \dotsc, P$ and $j = P+1, \dotsc, n$.} Define the transformed variables 
\begin{align}
\tilde{s} = Cs, \quad \g = C^{-T} g, \quad \H = C^{-T} H C^{-1},  \quad \W = C^{-T}WC^{-1}, \quad \T = \diag\{t_j\}_{1 \le j \le n}
\label{change of variable dtm}
\end{align}
where $t_1 = \dotsc = t_P =1, t_{P+1} = \dotsc = t_n = 0.$ {Note that  \(\tilde{W}\) is derived from \(W\) after a change of basis. } Then,
\begin{align*}
\m(\s)  = f_0 + \g^T  \s  + \frac{1}{2}\s ^T \H \s  +\frac{1}{6} \T[\s]^3  + \frac{\sigma}{4}\big\|\s \big\|_{\W}^4 = m_3^{\mathbbm{P}}(s).
\end{align*}
\end{theorem}

\begin{proof}
Note that $Cs = \big[({\boldsymbol{a}^{(1)}})^Ts, \dotsc, ({\boldsymbol{a}^{(P)}})^Ts, \dotsc\big]^T \in \R^n$, and let $[Cs]_k$ denote the $k$th entry of $Cs$. Then
\[
\T[Cs]^3 = \sum_{k=1}^n t_k^3 ([Cs]_k)^3 = \sum_{k=1}^P ({\boldsymbol{a}^{(k)}}^Ts)^3 \underset{\eqref{rank P dtm}}{=} T^{\mathbbm{P}}[s]^3,
\]
where the first equality arises from the tensor operation on the diagonal tensor $\T$, and the second equality holds because $t_1 = \dotsc = t_P =1$ and $t_{P+1} = \dotsc = t_n =0$.

Since $\hat{C}$ has linearly independent columns. Consequently, $C := \begin{bmatrix} \hat{C} & \hat{C}_\perp \end{bmatrix}^T$ is full rank, and $C^{-1}$ exists. Using $T^{\mathbbm{P}}[s]^3 = \T[Cs]^3$, we rewrite \eqref{rank P m dtm} as follows:
\begin{eqnarray}
m_3^{\mathbbm{P}}(s)  &=& f_0 + g^T C^{-1} C s  + \frac{1}{2}(C^{-1} C s)^T H C^{-1} C s + \frac{1}{6}\T[Cs]^3 + \frac{\sigma}{4}\big\|C^{-1} C s\big\|_W^4 \notag \\
&\underset{\tilde{s} = Cs}{=}& f_0 + (C^{-T} g)^T  \s  + \frac{1}{2}\s ^T (C^{-T} H C^{-1}) \s  +\frac{1}{6} \T [\s]^3  + \frac{\sigma}{4} \big[\s^T (C^{-T}WC^{-1}) \s\big]^2.\label{low-rank temp}
\end{eqnarray} 
Note that $\H = C^{-T} H C^{-1}$ is a symmetric matrix, and $\W = C^{-T}WC^{-1} \succ 0$. Therefore, the norm $\|\s\|_{\W}^2 = \s^T (C^{-T}WC^{-1}) \s$ is well-defined. Substituting $\g = C^{-T} g$, $\H = C^{-T} H C^{-1}$, and $\W = C^{-T}WC^{-1}$ into \eqref{low-rank temp} yields the desired result.
\end{proof}

In practice, many objective functions have a rank-one third-order derivative or a low-rank tensor term. A basic example is a function $f : \R^n \rightarrow \R$ with $x = [x_1, \dotsc, x_n]^T \in \R^n$ given by
\[
f(x) = (a_1 x_1 + \dotsc + a_n x_n)^3 + \text{lower-degree separable cross terms}.
\]
In this case, the third-order derivative of $f$ is rank-one. The lower-order separable cross terms may include combinations such as $x_jx_k$, $x_j^2$, and $x_j$ for any $j, k = 1, \dotsc, n$.

Furthermore, tensor data from real-world scenarios (e.g., in image processing \cite{loh2011high}, computer vision \cite{panagakis2021tensor}, machine learning \cite{kolda2009tensor, wang2008tensor}, and bioinformatics \cite{broadbent2024deciphering}) often lie in low-dimensional subspaces and exhibit low-rank structures. The problem of low-rank tensor learning, where only a limited number of samples are available to infer the underlying tensor, can be formulated as
\[
\min_{\mathcal{X} \in \hat{D}} f_{\hat{n}, \mathcal{Y}}(\mathcal{X}) + \hat{\beta} \cdot \text{rank}(\mathcal{X}),
\]
where $f_{\hat{n}, \mathcal{Y}}(\cdot)$ is a loss function dependent on the number of samples $\hat{n}$ and the observed tensor $\mathcal{Y}$. Here, $\hat{D}$ represents a given constraint set, $\hat{\beta} > 0$ is a regularization parameter, and the parameter tensor $\mathcal{X} \in \R^{n^3}$ and observed tensor $\mathcal{Y}$ have low-rank structures. 
% Additionally, objective functions $f: \R^n \rightarrow \R$ may have separable third-order derivatives that satisfy
% $
% \frac{\partial^3{f}}{\partial{x_\iota}\partial{x_j}\partial{x_\kappa}} =  \frac{\partial{f}}{\partial{x_\iota}} \frac{\partial{f}}{\partial{x_j}} \frac{\partial{f}}{\partial{x_\kappa}}.
% $
% The third-order derivative, in this case, has a rank-$1$ structure
% $
% \mathbb{T}^{(1)} = \nabla f(x) \otimes \nabla f(x) \otimes \nabla f(x).
% $
{In \Cref{sec: sphere}, we give an illustrative example—namely, the sphere-packing problem—where low-rank and spare third-order tensors emerge empirically at the optimized configuration. Whether such behaviour persists more broadly in large-scale real-world problems remains an open question for future investigation.
}

Extensive literature provides methods for finding the best rank-$1$ approximation for {symmetric} third-order tensors, including higher-order singular value decomposition (HOSVD), a multilinear extension of the singular value decomposition \cite{de2000multilinear}, the higher-order power method (HOPM) \cite{de1996independent}, and symmetric-HOPM \cite{de2000best, kofidis2002best}. 
While these methods typically operate on existing tensors, the tensor term in \eqref{ar3 model} possesses a special structure, as it represents third-order derivative information and can thus be estimated from lower-order data. Recently proposed Generalized quasi-Newton methods \cite{welzel2023generalizing} approximate the third-order tensor using Hessian evaluations, resulting in a low-rank structure. Additionally, \cite{schnabel1991tensor, schnabel1984tensor} employ $p$ interpolation points (not necessarily consecutive) to derive tensor approximation expressions in the form
\begin{equation}
\tag{Schnabel Model}
m^{\mathcal{S}}(s) := f_0 + g^T s + \frac{1}{2} H [s]^2 + \frac{1}{2} \sum_{k=1}^P \big({\boldsymbol{a}^{(k)}}^T s\big)^2 \big({\boldsymbol{b}^{(k)}}^T s\big) + \frac{1}{6} \sum_{k=1}^P \gamma_k \big({\boldsymbol{a}^{(k)}}^T s\big)^4,
\label{sch dtm}
\end{equation}
where $\gamma_k > 0$ and $\boldsymbol{a}^{(k)}, \boldsymbol{b}^{(k)} \in \R^n$.
In the case where the high-order derivatives are low-rank, and we have ${\boldsymbol{a}^{(k)}} = {\boldsymbol{b}^{(k)}}$, using a similar change of basis as in \Cref{thm change of basis dtm}, \eqref{sch dtm} can be transformed into a separable quartically regularized polynomial \eqref{separable Polynomial}, 
$$
m^{\mathcal{S}}(\s) = f_0 + \g^T \s + \frac{1}{2} \H [\s]^2 + \frac{1}{2} \sum_{k=1}^P \s_k^3 + \frac{1}{6} \sum_{k=1}^P \gamma_k \s_k^4 = f_0 + \g^T \s + \frac{1}{2} \H [\s]^2 + \frac{1}{2} \T [\s]^3 + \frac{1}{6} \Sigma [\s]^4,
$$
where $\s = [\s_1, \dotsc, \s_n]^T \in \R$, $\T \in \R^{n^3}$ is a diagonal tensor with $P$ nonzero entries, and $\Sigma \in \R^{n^4}$ is a diagonal tensor with nonzero entries $\{\gamma_1, \dotsc, \gamma_P\}$. %{(KZ: Last time you asked me under what situation we have ${\boldsymbol{a}^{(k)}} = {\boldsymbol{b}^{(k)}}$ for \eqref{sch dtm}. The answer to this is not so straightforward.  \eqref{sch dtm} is also quite different from our $m_3$ model. \eqref{sch dtm} was originally used to approximate the Taylor expansion up to $4$th order derivative, while we use the quartic term as regularization. )}
Developing low-rank tensor approximations for our subproblem and efficiently incorporating these approximations into \Cref{TDTM} is a standalone area of study. Since this paper focuses on designing a convergent iterative algorithmic framework for \eqref{ar3 model}, we leave the tensor approximation for future work.

\section{An Algorithm for Minimizing the AR3 Model}
\label{sec: Secular Equation dtm}
%{(KZ: This section is amended for a general $T$. Remarks for a diagonal $T$ are added. Major changes are highlighted in blue. )}

%{look at when to close the gap, look at what to do with optimality condition. Newton! When is Newton easy to solve?  diagonal.}

%{Before DTM on full tensor. We highlight the case of diagonal and low-rank and say how to solve it by Newton. Newton can be applied to optimality conditions in full. local min vs global min.}
 
% {(Optimality conditions can apply directly, but not tensor-free while solving the minimization. Or a change of basis which is tensor involved, but you can minimize the DTM after change of basis as tensor free.)}

% To minimize  \eqref{ar3 model}, we use the global optimality characterizations outlined in {\Cref{thm necessary tight} and \Cref{thm sufficiency tight}}.   

In this section, we describe how to exploit the global optimality characterizations established in \Cref{thm necessary tight} and \Cref{thm sufficiency tight} to minimize \eqref{ar3 model} (see \Cref{Newton algo dtm}). 
While \Cref{Newton algo dtm} could apply to a general \eqref{ar3 model}, it is more scalable for certain structured settings, such as the diagonal or low-rank cases. 
Motivated by this observation, in \Cref{Sec Iterative Algorithm for Minimizing AR3 Model} we develop a more implementable and scalable framework (\Cref{TDTM}) that reduces the general \eqref{ar3 model} to a sequence of diagonal tensor subproblems, each solved in the inner loop using \Cref{Newton algo dtm}.

In this section, we first introduce the secular equation, defined in \eqref{secular1 dtm}--\eqref{secular2 dtm}, based on the sufficient condition. 
This section aims to formulate an algorithm for locating the root of the (nonlinear) secular equation, which is defined in \eqref{secular1 dtm}--\eqref{secular2 dtm}. We employ  Cholesky factorization
 and Newton's iteration to find a root of the secular equation. 
A convergence analysis for the root-finding method is deferred to future work. We seek  a vector $s(\lambda)$ which satisfies the system of equation,
\begin{align}
& \B_d(s) s := \bigg(H + \frac{1}{2}\Gamma + \lambda W \bigg)s = -g , 
\label{secular1 dtm}
\\&  H +  \frac{2}{3} \Gamma  + \lambda W  -\frac{L_W}{3}  W- \frac{\Lambda_W ^2}{18\sigma} W \succeq 0,  
\label{secular3 dtm}
\\\text{where} \qquad &\lambda=  \sigma  \|s\|_W^2,  \qquad L_W :=  \Lambda_W \|s\|_W, \qquad  \Gamma  = T s.
\label{secular2 dtm}
\end{align}
Note that in the case of $m_3$ with diagonal tensor \eqref{m3 diagonal}, we have $ L_W =  \hat{t} \|s\|_W, \Gamma  = \diag\big\{t_js_j\big\}_{1\le j \le n}.$ Let $D_\Gamma  \in \R^{n \times n} $ be a diagonal matrix with diagonal entries $\lambda_1 \le \lambda_2 \le \dotsc \le \lambda_n $  as the eigenvalues of $H + \frac{1}{2}\Gamma$. Let  $U \in \R^{n \times n}$ with each column $u_1, u_2, \dotsc, u_n \in \R^n$ correspond to the generalized eigenvectors. Then,
\begin{eqnarray}
  \bigg(H + \frac{1}{2}\Gamma \bigg) U = W   U D_\Gamma , \text{ where }   U^T W U =  U W U^T = I_n.
  \label{gener e value}
\end{eqnarray}
Define $\hat{g}:=U^T g$ and $\hat{s}=W U^Ts$.  From this, we may deduce that  $s = U \hat{s}$ and $g = W U\hat{g}$. It follows from \eqref{secular1 dtm} that 
\begin{eqnarray*}
   -WU \hat{g} = -g = \bigg(H + \frac{1}{2}\Gamma+ \lambda W\bigg) s =    \bigg(H + \frac{1}{2}\Gamma+ \lambda W\bigg) U  \hat{s} \underset{\eqref{gener e value}}{=} WU\big(D_\Gamma + \lambda I_n\big)  \hat{s} . 
\end{eqnarray*}
Since $U$ is non-singular, therefore $ \hat{s}:= (D_\Gamma + \lambda I_n)^{-1} \hat{g}. $ For $\lambda \neq \lambda_j(\Gamma)$ for any $1 \le j\le n$, 
\begin{eqnarray}
 \Psi_d(\lambda) := \|s(\lambda)\|_W^2 = \|\hat{s}\|^2= \bigg\|(D_\Gamma   + \lambda W )^{-1}\hat{g}\bigg\|^2 = \sum_{j=1}^n \frac{\hat{g}_j^2}{\big(\lambda +\lambda_j\big)^2}
\label{general eig value dtm}
\end{eqnarray}
where $\hat{g}_j$ is the $j$th component of $\hat{g}=U^T g$. 
Using \eqref{general eig value dtm},  we convert  $\lambda=  \sigma \|s\|_W^2$ to the following system of nonlinear equations
\begin{eqnarray}
 \phi_d(\lambda):= {\frac{1}{\sqrt{\Psi_d(\lambda)}}} -  \frac{1}{ \K(\lambda)} = 0 \quad
 \text{where} \quad  \K(\lambda) := \big[\lambda/\sigma\big]^{1/2}. 
  \label{phi_1 dtm}
\end{eqnarray}
Assume that the sufficient condition \eqref{secular3 dtm} is satisfied, then we have 
\begin{eqnarray}
\mathcal{B}(s)s  = H + \frac{1}{2} \Gamma  +   \lambda W \succ 
\bigg(H + \frac{1}{2} \Gamma  +   \lambda W\bigg) +  \underbrace{\bigg(\frac{1}{6} \Gamma - \frac{1}{3}  L_W\bigg) }_{\prec 0}  \underset{\eqref{secular3 dtm}}{\succ} 0.
\label{b psd}
\end{eqnarray}
To numerically determine the root of \eqref{phi_1 dtm}, we employ a univariate Newton iteration combined with root-finding techniques based on Cholesky factorization. The update steps for the Newton iteration are outlined in \Cref{Newton lemma dtm}, with the proof framework based on \cite[Lemma 2.2]{zhu2023cubic}, \cite[Lemma 6.1]{cartis2011adaptive}, and \cite[Lemma 7.3.1]{conn2000trust}.

\begin{lemma} 
Assume $H + \frac{1}{2} \Gamma  +   \lambda W \succ 0$ and $g \neq 0$,  the Newton iteration updates are 
\begin{eqnarray}
 \Delta \lambda^{(l)} =  -\frac{\phi_d(\lambda^{(l)})} {\nabla_{\lambda} \phi_d (\lambda^{(l)})} = \frac{\|s(\lambda^{(l)})\|_W^{-1} - \K(\lambda^{(l)})^{-1}}{\|\omega(\lambda^{(l)})\|_W^2\|s(\lambda^{(l)})\|_W^{-3} - \big[\nabla_{\lambda} \K(\lambda^{(l)})^{-1}\big]},    \qquad \lambda^{(k+1)} = \lambda^{(l)}+     \Delta \lambda^{(l)}
     \label{Newton correction dtm}
\end{eqnarray}
where $\omega (\lambda)$ is given by $(H + \frac{1}{2}\Gamma + \lambda W) = L(\lambda) L^T(\lambda)$ and $W^{-1}L(\lambda) \omega (\lambda) = s(\lambda)$. Also, 
\begin{eqnarray*}
\K(\lambda)^{-1} = \big[\sigma / \lambda\big]^{1/2} , \qquad  \big[\nabla_{\lambda} \K(\lambda)^{-1}]' = - \frac{1}{2}{\sigma}^{1/2} \lambda^{-3/2}. 
\end{eqnarray*}
\label{Newton lemma dtm}
\end{lemma}

\begin{proof} Since $H + \frac{1}{2} \Gamma  +   \lambda W \succ 0$, the Cholesky factorization gives $H_W(\lambda) :=(H + \frac{1}{2}\Gamma + \lambda W) = L(\lambda) L^T(\lambda)$. Differentiating $(H+ \frac{1}{2}\Gamma+ \lambda W) s(\lambda) = -g$ with respect to $\lambda$ gives $    (H+ \frac{1}{2}\Gamma+ \lambda W)  \nabla_\lambda s(\lambda) + W s(\lambda) =0. $ Therefore, 
\begin{eqnarray}
\nabla_\lambda  s(\lambda) = -H_W^{-1}(\lambda) Ws(\lambda) = L(\lambda)^{-T}   L(\lambda)^{-1} {W}s(\lambda) .
\label{ds appendix dtm}
\end{eqnarray}
On the other hand, $\Psi_d(\lambda) = \|s(\lambda)\|^2_W$ has derivative $\nabla_\lambda \Psi_d(\lambda)  =  2 \big\langle Ws(\lambda), \nabla_\lambda s(\lambda)\big\rangle$ where $\big\langle \cdot, \cdot\big\rangle$ denotes the Euclidean inner product. 
Substituting \eqref{ds appendix dtm} gives, 
\begin{eqnarray}
\nabla_\lambda \Psi_d(\lambda)  &=& 
  - 2  \big\langle Ws(\lambda),  
  L(\lambda)^{-T}  L(\lambda)^{-1} {W} s(\lambda) \big\rangle \notag
\\  &=& -2  \big\langle L(\lambda)^{-1} W s(\lambda),  L(\lambda)^{-1} {W} s(\lambda) \big\rangle 
= -2 \big\|L(\lambda)^{-1}   W s(\lambda)\big\|= -2 \big\|\omega(\lambda)\big\|^2
\label{grad phi}
\end{eqnarray}
where ${W^{-1}}L(\lambda) \omega (\lambda) = s(\lambda)$. Consequently, 
\begin{eqnarray}
\nabla_\lambda (\Psi_d(\lambda)^{-1}) = \nabla_\lambda (\Psi_d(\lambda)^{-1/2}) = - \frac{1}{2}\Psi_d(\lambda)^{-3/2} \nabla_\lambda \Psi_d(\lambda)  \underset{\eqref{grad phi}}{=}  \|s(\lambda)\|_W^{-3}\|\omega(\lambda)\|^2.
\label{derivative dtm}
\end{eqnarray}
Substituting \eqref{derivative dtm} into \eqref{phi_1 dtm} and then applying the Newton updates $\Delta \lambda^{(l)} = -\frac{\phi_d(\lambda^{(l)})} {\nabla_{\lambda} \phi_d (\lambda^{(l)})}$ yields the result.
\end{proof}

An algorithm to solve \eqref{phi_1 dtm} is provided in \Cref{Newton algo dtm}, where we alternately update $\Gamma$ and apply Newton steps to the minimizer of $m_3$. 
The approach is as follows: at each $k$th iteration, we first compute and fix $\Gamma := T [s^{(\kappa)}]$. Using this fixed $\Gamma$, we then perform a univariate Newton iteration in the second while loop to find $s^{(l)}$ such that $\phi_d(\lambda^{(l)}) = 0$. After this step, we update $s^{(\kappa+1)} := s^{(l)}$ and reset $\Gamma$ based on the new value of $s^{(\kappa+1)}$. This iterative process is repeated until convergence, ultimately yielding an approximate root of \eqref{secular1 dtm}.

{
\begin{remark} (Connection of \Cref{Newton algo dtm} to Global Optimality Conditions) The global optimality conditions are implicitly built into \Cref{Newton algo dtm}. Specifically, the first-order condition yields the linear system \eqref{secular1 dtm}
and the sufficient global condition ensures that the matrix system is positive definite (i.e.,  $\mathcal{B}(s) \succeq 0$  see \eqref{b psd}). These allow us to reformulate the minimisation of \(m_3\) as solving a sequence of positive definite linear systems.  When Algorithm~1 converges, then we have an (approximate) stationary point satisfying the global necessary condition
\[
H + \tfrac23 \Gamma + \lambda W + \tfrac{\Lambda_W}{3}\|s_\ast\|_W =
(H + \tfrac12 \Gamma + \lambda W) + \underbrace{(\tfrac16 \Gamma + \tfrac{\Lambda_W}{3}\|s_\ast\|_W)}_{\succeq 0}
\succeq 0.
\]
We show in \Cref{thm epi T} that, if for sufficiently large \(\sigma\), the necessary and sufficient conditions coincide, then the algorithm is guaranteed to obtain the global minimiser of \(m_3\) (whether convex or nonconvex).
\end{remark}
}

\begin{algorithm}
\caption{\small \textbf{Algorithm for Minimizing the AR3 Model} using global optimality characterizations. It locates a root of the associated nonlinear secular equation defined in \eqref{secular1 dtm}--\eqref{secular2 dtm}, using Cholesky factorization together with Newton iterations  to solve \eqref{phi_1 dtm}.}

\textbf{Input}: $g$, $H$, $T$, $\sigma > 0$. Maximum iterations $l_{\max}$ and $\kappa_{\max}$; accuracy levels $ \epsilon_{\kappa} >\epsilon_l > 0$. \\
\textbf{Initialization}: Set $s^{(0)} = \boldsymbol{0}$, $\Gamma =  T[s^{(0)}]$, $\lambda^{(0)} := \max\{-\lambda_{\min}[H], 0\}$. Initialize $l = 0$, $\kappa = 0$. 
\\\textbf{Main Step:}
\While{$\kappa \le \kappa_{\max}$ and $\big\|(H + \frac{1}{2}\Gamma + \lambda^{(l)} W)s^{(\kappa)} - g\big\| > \epsilon_{\kappa}$}
{    \begin{equation}
    \tag{$\lambda$-correction step}
    \lambda^{(0)} := \max \bigg\{\lambda^{(l)},  \lambda_{\min}\big[H + \frac{1}{2}\Gamma + \lambda^{(l)} W\big], 0\bigg\};
    \label{lambda correction step}     
    \end{equation}
    Reset $l = 0$.
    \\ 
    \While{$l \le l_{\max}$ and $\big|\phi_d(\lambda^{(l)})\big| \ge \epsilon_l$}
    {   \textbf{Newton Step:} Fix $\Gamma$ in this while loop, update $s^{(l)}$ and $\lambda^{(l)}$ using Newton’s updates.
        \\ Cholesky factorize $H + \frac{1}{2}\Gamma + \lambda^{(l)} W = LL^T$. 
        \\\If{Cholesky factorization fails} 
        {Return, go to \textbf{Safeguard Step}.}
        Solve $LL^Ts(\lambda^{(l)}) = - g$ and compute $\omega(\lambda^{(l)}) = L^{-1} W s(\lambda^{(l)})$. \\
        Update $\lambda^{(l)}$ using the Newton step in {\Cref{Newton correction dtm} and \eqref{Newton correction dtm}} and $l := l + 1$. 
    }
    Set $s^{(\kappa+1)} := s(\lambda^{(l)})$, $\Gamma := T[s^{(\kappa+1)}]$, $\kappa := \kappa + 1$. 
}
\textbf{Safeguard Step:}
\If{$ \kappa = \kappa_{\max}$ or $l = l_{\max}$ or Cholesky factorization fails}
{Return to solve using standard solvers, i.e., ARC solver.}
\label{Newton algo dtm}
\end{algorithm}

\begin{remark} ($\lambda$ Correction Step)
    Since we alternately update $\Gamma := T [s^{(\kappa)}]$ and use the Newton step to compute \( s^{(\kappa+1)} \), the value \(\lambda^{(l)}\) obtained from the Newton step ensures that \( H + \frac{1}{2} T [s^{(\kappa)}] + \lambda^{(l)} W \succeq 0 \). This is why \eqref{lambda correction step} is used to apply a slight adjustment, ensuring that \( H + \frac{1}{2}  T [s^{(\kappa)}] + \lambda^{(l)} W \succeq 0 \) holds.
\end{remark}

\begin{remark} (Special Case for Diagonal $T$)  
%Note that for the full tensor $T$, \Cref{Newton algo dtm} is not tensor free. However, there are cases where the Newton system is easier to solve. 
If \( T \) is a diagonal tensor, we have 
$\Gamma  = \diag\big\{t_js_j\big\}_{1\le j \le n} = \textbf{t} \cdot s^{(\kappa)}$
where \( \textbf{t} := [t_1, \dotsc, t_n]^T \) and \( L_W \) can be explicitly expressed as 
$
L_W = \hat{t} \|s\|_W.
$
The sparse structure of the diagonal tensor enables the computation of \( \Gamma \) without requiring tensor-vector multiplications; instead, only vector multiplications are needed. This significantly improves both storage efficiency and computational time. Numerical results demonstrating the computational time improvement achieved with diagonal tensors that avoid tensor-vector multiplications are illustrated in \Cref{fig computational time}.
\end{remark}

\begin{remark} 
   In \Cref{Newton algo dtm}, the second-order global sufficient optimality condition ensures that $\mathcal{B}$ is positive semidefinite. However, it is not guaranteed that $\mathcal{B}^{(l)} = H + \frac{1}{2}\Gamma + \lambda^{(l)} W$ remains positive semidefinite throughout the iterations, nor that the minimum obtained by \Cref{Newton algo dtm} satisfies the second-order global necessary optimality condition (\Cref{fig: increase t or sigma}). As a safeguard, whenever Cholesky factorization fails in \Cref{Newton algo dtm} or the maximum number of iterations is exceeded, \Cref{Newton algo dtm} reverts to the ARC method proposed by Cartis et al. \cite{cartis2011adaptive}.  The numerical results in \Cref{sec numerics} indicate that failures in the Cholesky factorization of $\mathcal{B}^{(l)}$ usually happens due to the ill-conditioning, and $\mathcal{B}^{(l)}$ typically remains positive definite throughout the algorithm.
\end{remark}

While \Cref{Newton algo dtm} efficiently finds stationary points of {the  quartically regularized cubic subproblem}, further research is necessary to establish its convergence.
% In \Cref{Newton algo dtm}, the second-order global necessary optimality condition ensures that $\mathcal{B}$ is positive semidefinite. 
% Due to a gap between the sufficient and necessary conditions for \eqref{M_d}, there may exist a global minimum that satisfies the sufficient condition \eqref{cqr global sufficiency condition} but does not fulfill \( H + \frac{1}{2}\Gamma + \lambda^{(l)} W \succeq 0 \). In such cases, Cholesky factorization may not be applicable.
After obtaining a stationary point from \Cref{Newton algo dtm}, we always test whether it satisfies the second-order local optimality condition \eqref{local iff2 m3} and the second-order sufficient global optimality condition \eqref{sufficient tight DTM}. Our numerical experiments in \Cref{sec numerics} indicate that the former is always satisfied, while the latter is usually satisfied for simple or small tensor terms or when $n$ is small.
Additionally, a hard case may occur when the eigenvector corresponding to the largest eigenvalue of \( H + \frac{1}{2}\Gamma \), denoted as \( u_1 \), is orthogonal to \( g \), i.e., \( u_1^T g = 0 \). In this scenario, a similar analysis to the trust-region subproblem \cite[Sec 8.3.1]{cartis2022evaluation} would need to be applied, resulting in non-unique global minimizers for the subproblem. We refer readers to \cite[Sec 8.3.1]{cartis2022evaluation} for a detailed discussion. 
To conclude, note that in this subsection, we employ Cholesky factorization and Newton's method to address the secular equation, though this is only one of several possible approaches. For instance, scalable iterative algorithms could be developed using Krylov methods, eigenvalue-based formulations, or subspace optimization. While we acknowledge the potential of these techniques for minimizing the diagonal tensor polynomial our primary aim here is to utilize the diagonal tensor polynomial in minimizing the \( m_3 \) subproblem. Consequently, a comprehensive analysis of these alternative approaches and the convergence of \Cref{Newton algo dtm} is left for future work.

\section{A Diagonal Tensor Algorithmic Framework for Minimizing the AR3 Model}
\label{Sec Iterative Algorithm for Minimizing AR3 Model}

In this section, we develop an implementable framework for solving the general \eqref{ar3 model} by iteratively minimizing a sequence of diagonal tensor subproblems. This framework is outlined in \Cref{Sec Iterative Algorithm for Minimizing AR3 Model} and \Cref{TDTM}.  At each iteration, we approximate \eqref{ar3 model} by a quartically regularized polynomial with a diagonal tensor term, referred to as the Diagonal Tensor Model (DTM). The outer iteration reduces the general problem to a sequence of diagonal tensor models, each of which is solved in the inner loop using \Cref{Newton algo dtm}. 
Compared to directly applying \Cref{Newton algo dtm} to the full model \eqref{ar3 model}, this approach is advantageous because the diagonal tensor structure avoids operations involving dense tensor entries and leads to simpler optimality conditions.

We introduce the notations used in the iterative algorithm. Since we are working with a single \eqref{ar3 model} with fixed $x_k$, in this section, we simplify notation by dropping the $x_k$ and writing the fourth-order Taylor expansion of $m_3 (s^{(i)}+s)$ at $s^{(i)}$ as 
\begin{equation}
M(s^{(i)}, s)  := f_i + g_i^T s + \frac{1}{2}s^T H_i s +\frac{1}{6}T_i[s]^3 + \frac{\sigma}{4}\|s\|_W^4,
\label{M(s, s^i)}
\end{equation} 
where $f_i =  m_3(s^{(i)}) \in \R$, $g_i = \nabla m_3(s^{(i)}) \in \R^n$, $H_i = \nabla^2 m_3(s^{(i)}) \in \R^{ n \times n}$ and $T_i = \nabla^3 m_3(s^{(i)}) \in \R^{ n^3}$ and $\nabla$ denotes the derivative with respect to $s$.  Since $m_3(s)$ is a fourth-degree multivariate polynomial, the fourth-order Taylor expansion is exact. Therefore, we have
\begin{eqnarray*}
M(s^{(i)}, s) =m_3(s+s^{(i)}),   \qquad  \text{ and } \qquad \min_{s\in \R^n}m_3(s) = \min_{s\in \R^n}M(s^{(i)}, s). 
\end{eqnarray*}
For all $s^{(i)}$,  the Hessian $H_i$ is a symmetric matrix, and $T_i \in \R^{n\times n \times n}$ is a {symmetric} third-order tensor. The first-order and second-order derivatives of $M(s^{(i)}, s)$ are given by $
    \nabla M(s^{(i)}, s) = g_i + H_i s +\frac{1}{2} T_i[s]^2 + \sigma\|s\|_W^2 Ws,$ and $ \nabla^2 M(s^{(i)}, s) = H_i +T_i[s] + \sigma \big[\|s\|^2 W + 2 (Ws)(Ws)^T\big]$, respectively.
%We refer to the model $M(s^{(i)}, s)$ as locally convex at $s^{(i)}$ if $H_i \succeq 0$, and it is locally nonconvex at $s^{(i)}$ if the smallest eigenvalue of $H_i$, denoted as $\lambda_{\min}[H_i]$, is negative. 
Our aim is to find a sequence of $\{s^{(i)}\}_{i \ge 0}$ such that $s^{(i)}$ converges to a  local minimizer of $m_3$ satisfying
$
\big\|g_i\big\| \leq \epsilon, 
$
where $\epsilon$ are the tolerance for minimizing the AR$3$ subproblem. 
A recent advancement introduced in \cite{cartis2023second}, the Quadratic Quartic Regularisation (QQR) method, is an iterative algorithm specifically designed to solve \eqref{ar3 model}. In QQR, the third-order tensor term $T_i[s]^3$ is approximated by a linear combination of quadratic and quartic terms. Following QQR, the Cubic Quartic Regularization (CQR) framework was introduced in \cite{zhu2023cubic}. In the CQR method, the third-order tensor term is approximated by a single-parameter cubic term $\beta\|s\|^3$. However, since $T_i[s]^3$ can span multiple tensor directions, it is natural to consider approximating the tensor term with a more comprehensive tensor expression.

Motivated by this, in this section, we introduce the diagonal tensor model that approximates $T_i$ using a diagonal tensor term $\T_i \in \R^{n^3}$ with diagonal entries $\{t^{(i)}_1, \dotsc, t^{(i)}_n\}$. 
The Diagonal Tensor Model (or Diagonal Tensor Polynomial, i.e., DTM) is defined as
\begin{equation}
M_d(s^{(i)}, s)   = f_0 + g_i^T s+\frac{1}{2} s^TH_i s  +  \frac{1}{6} \underbrace{\sum_{j=1}^n t^{(i)}_j s_j^3}_{:=\T_i[s]^3}+ \frac{1}{4} \underbrace{(\sigma+d_i)}_{:=\sigma_d^{(i)}}\|s\|_W ^4,
\label{M_d}
\tag{Diagonal Tensor Model}
\end{equation}
where $W \succ 0$ is a symmetric matrix, the matrix norm is defined as $\|v\|_W:= \sqrt{v^TWv}$ and $d_i >0$.
Note that each $t^{(i)}_j \in \R$ is a scalar and $t^{(i)}_j$ can be chosen as negative. For instance, we can use  $t^{(i)}_j$ as the diagonal entries of $  \nabla^2 m_3(s^{(i)})$. More details for choosing  $t^{(i)}_j$ can be found in \Cref{algo: pdtm} and \Cref{sec: low-rank}.  
The parameter $d_i$ is an adaptive regularization weight; $\sigma_d^{(i)} = \sigma+ d_i \ge 0$  ensures that $M_d$ remains bounded from below.

There are several advantages and motivations for using  \eqref{M_d}  as an approximation for  \eqref{ar3 model}.
Firstly,  a diagonal tensor term provides $n$ degrees of freedom in the tensor direction to approximate and capture information from the third-order derivative. 
Many  tensor data from real-world
scenario lies in low-dimensional subspaces and exhibits low-rank structures. In \Cref{sec: low-rank}, we will show that if \eqref{ar3 model} has a tensor term with rank $<n$, then \eqref{M_d} can exactly recover \eqref{ar3 model} under a change of basis (\Cref{sec: low-rank}). 

Since \eqref{M_d} is a quartic polynomial, the optimality conditions specified in \Cref{thm necessary tight} and \Cref{thm sufficiency tight} apply to it. Due to the special diagonal structure of the tensor in \eqref{M_d}, the optimality conditions are significantly simplified, allowing for an explicit expression of $\Lambda_W$. According to \Cref{example special T}, we have 
\begin{equation}
\Lambda_W = \hat{t}_i = \max \{|t_1^i|, |t_2^i|, \dotsc, |t_n^i|\} [ \lambda_{\min}(W)]^{-3/2}.
\label{LW for dt}
\end{equation}
The global necessary conditions for \eqref{M_d} state that if $s_d^{(i)} = \argmin_{s \in \mathbb{R}^n} M_d(s^{(i)}, s)$, let $\B_d_i(s) := H_i + \frac{1}{2}\T_i [s] + \sigma_d^{(i)} W \|s\|_W^2$.
%$\G_d_i(s) := H_i +  \T_i[s] +  \sigma_d^{(i)} W\|s\|_W^2 $
Then,
\begin{equation}
%\tag{Global Necessary Condition}
    \B_d_i(s_d^{(i)}) s_d^{(i)} = -g,   \qquad
 H +  \frac{2}{3} \T [s_d^{(i)}]  +\sigma \|s_d^{(i)}\|_W^2 W  +\frac{\hat{t}_i}{3} \|s_d^{(i)}\|_WW\succeq 0.
    \label{cqr global necessary condition}
\end{equation}
On the other hand,  the global sufficiency conditions for \eqref{M_d} gives if
\begin{equation}
%\tag{Global Sufficiency Condition}
\B_d_i(s_d^{(i)}) s_d^{(i)}  = -g, \qquad
 H +  \frac{2}{3} \T [s_d^{(i)}]  +\sigma \|s_d^{(i)}\|_W^2 W  -\frac{\hat{t}_i}{3} \|s_d^{(i)}\|_WW- \frac{\hat{t}_i^2}{18\sigma} W \succeq 0.
\label{cqr global sufficiency condition}
\end{equation}
Then, $s_d^{(i)} = \argmin_{s\in \R^n} M_d(s^{(i)}, s)$.

Lastly, notice that the sparse structure of the diagonal tensor term and the expression of $\Lambda_W$ in \eqref{LW for dt} allow us to develop a tensor-free framework for minimizing \eqref{ar3 model}. Operations related to the third-order term and the computation of $\hat{t}_i$ can be effectively converted into vector operations with size $\mathcal{O}(n)$.
Specifically, let $\textbf{t}^{(i)} = [t^{(i)}_1, \dotsc, t^{(i)}_n]^T$, then we have $\T_i[s] = \textbf{t}^{(i)} \cdot s$, $\T_i[s]^2 = \textbf{t}^{(i)} \cdot s \cdot s$, and $\T_i[s]^3 = \textbf{t}^{(i)} \cdot s \cdot s \cdot s$, where $\cdot$ represents the entrywise dot product.
Therefore, the evaluation of $M_d(s^{(i)}, s)$ and its derivatives, remains tensor-free, as we avoid operating with the $\mathcal{O}(n^3)$ (dense) tensor term $T$. This framework saves both storage and computational time, allowing for a more efficient algorithm. Numerical experiments illustrating this efficiency can be found in \Cref{sec numerics} and \Cref{fig computational time}. %{The key numerical bottleneck in applying tensor methods is the computation and storage of the third-order tensor term. Since this term involves $n^3$ variables and is typically dense, direct computation becomes infeasible for large-scale problems. To overcome these hurdles, we will adopt several innovative approaches.}
We introduce in \Cref{TDTM} the DTM algorithmic framework for minimizing \eqref{ar3 model} following a similar adaptive regularization framework as \cite{zhu2023cubic}. 

It is worth noting that the  Diagonal Tensor framework can be applied directly to the minimization of an objective function $f: \R^n \rightarrow \R$ under suitable assumptions (not just the AR$3$ subproblem). Specifically, if $f$ is three times continuously differentiable, bounded below, and the third derivative is globally Lipschitz continuous, then \Cref{TDTM} can also be used for the minimization of a potentially nonconvex objective function. However, our interest here is not to create a general method for optimizing $f$; instead, we focus on efficiently minimizing the local model $m_3(s)$. %\Cref{TDTM} outlines the framework of the CQR algorithm. For generality, we use ${\mathcal{m}}$ to denote the general objective function to be minimized.

\begin{algorithm}
\caption{\small Diagonal Tensor Algorithmic Framework with Adaptive Regularization}
\textit{Initialization}: Set $s^{(0)} = \boldsymbol{0} \in \R^n$ and $i=0.$ An initial regularization  
  parameter $d_0 = 0$, constants $ \eta_1 > \eta > 0$, $\gamma > 1 > \gamma_2 >0$.
\\
\textit{Input}: $f_0:=m_3(s^{(0)})$,  $g_0 := \nabla m_3(s^{(0)})$, $H_0 : = \nabla^2 m_3(s^{(0)})$, a symmetric $n \times n$ matrix $W \succ 0$, $\sigma>0$ a given regularization parameter;  an accuracy level $\epsilon>0$. 

\textbf{Step 1: Test for termination.} If $\|g_i\| = \|\nabla m_3(s^{(i)})\| \le \epsilon$, terminate with the $\epsilon$-approximate first-order minimizer $s_\epsilon = s^{(i)}$. 

\textbf{Step 2: Step computation.}
% {Compute\footnotesize{$^a$} \normalsize{}$\T_i = \diag\{t^{(i)}_j\}_{1 \le j \le n}$ to give a diagonal approximation of the (local) tensor information in $m_3$.} 
Use the coefficient $g_i, H_i, \T_i, \sigma, d_i$ to the \eqref{M_d}, $M_d$ and compute $s_d^{(i)} = \argmin_{s \in \R^n} M_d(s^{(i)}, s)$ using \Cref{Newton algo dtm} (see \Cref{remark solver s_d DTM}). When $T_i$ is diagonal or low-rank, the change of basis in the Low-Rank Rule of \Cref{algo: pdtm} reduces $M_d$ to a diagonal model, so that the matrices in the secular system \eqref{secular1 dtm}--\eqref{secular2 dtm} are diagonal; the step $s_d^{(i)}$ is then obtained directly via \Cref{algo: pdtm}. 

\textbf{Step 3: Acceptance of trial point.} Compute $m_3(s^{(i)}+s_d^{(i)})$ and define
\begin{equation}
\tag{Ratio Test}
\rho_i = \frac{m_3(s^{(i)}) -m_3(s^{(i)}+s_d^{(i)})}{m_3(s^{(i)}) - M_d(s^{(i)}, s_d^{(i)})}.
\label{ratio test}
\end{equation}
\textbf{Step 4: Regularization Parameter Update.} 
\If{ $\rho_i \ge \eta$.   \emph{Successful Iter.}}{
   Set $s^{(i+1)} = s^{(i)} + s_d^{(i)} $, update $g_{i+1}, H_{i+1}$ at $s^{(i+1)}$; get $\T_{i+1}$ using \Cref{algo: pdtm}; $d_{i+1} = d_i$. 
}
\ElseIf{$\rho_i \ge \eta_1$.     \emph{Very Successful Iter.}}{
 Set $s^{(i+1)} = s^{(i)} + s_d^{(i)} $, update $g_{i+1}, H_{i+1}$ at $s^{(i+1)}$; get $\T_{i+1}$ using \Cref{algo: pdtm}; decrease $d_{i+1} = \gamma_2 d_i$. 
}
\ElseIf{$\rho_i \le \eta$.   \emph{Unsuccessful Iter.} }{
 Set  $s^{(i+1)} = s^{(i)}$, $g_{i+1} = g_i$, $H_{i+1} = H_i$, $\T_{i+1} = \T_i$; increase $d_{i+1} = \gamma_2 d_i$. 
}

Repeat with $i:=i+1$. 

\label{TDTM}
\end{algorithm}

\begin{remark} (\textbf{Minimization of $M_d$}) 
While our proof in \Cref{sec: Convergence and Complexity of dtm} (in particular \Cref{lemma: lower bound of norms DTM}) assumes exact minimization of $\nabla M_d(s^{(i)}, s_d) = 0$, similar bounds can be derived when considering approximate minimization of $M_d$ with a step-termination condition, such as $\|\nabla M_d(s^{(i)}, s_d)\| \le \Theta \|s_d\|_W^3$ for some $\Theta \in (0, 1)$. Such step-termination conditions are also used in \cite{cartis2020concise}, with alternative variants available \cite{cartis2011adaptive}.  
\label{remark solver s_d DTM}
\end{remark}

\subsection{Convergence and Complexity}
\label{sec: Convergence and Complexity of dtm}
\Cref{algo: DTM variant 1} is a variant of the diagonal tensor framework (\Cref{TDTM}) for minimizing \(m_3\), for which we establish an optimal worst-case complexity bound. The convergence and complexity analysis of \Cref{algo: DTM variant 1} is an extension of the CQR algorithm in \cite{zhu2023cubic}. Since the complexity proofs in this section follow a similar framework to those presented in \cite[Sec 3.1]{zhu2023cubic} and \cite[Sec 2.4.1]{cartis2022evaluation}, we give the details of the proof in \Cref{appendix details of proof} and provide only a proof sketch and the conclusions in this section.

%{(KZ: I moved most proofs to appendix for this section)}

\begin{algorithm}
\caption{\small Diagonal Tensor Variant of \Cref{TDTM} for minimizing $m_3$}
Fix parameters $\hat{B}>0$ and  $\alpha \in (0, \frac{1}{2})$. Initialize and proceed through Steps 1 to 3 as \Cref{TDTM}. In Step 4, we use the following rule for updating the regularization parameter, $\T_{i+1}$ and $s^{(i+1)}$. 

\If{$\|\nabla m_3(s^{(i)}+s_d^{(i)})\| \le \epsilon$}
{{\textbf{Terminate} with an $\epsilon$-approximate first-order minimizer $s_\epsilon$ of  $m_3(s^{(i)}+s_d^{(i)})$.}}
\Else{ 
Let $\beta_i : = \T_i \left[\frac{s_d^{(i)}}{\|s_d^{(i)}\|_W} \right]^3$ where {$\T_i$ is updated using diagonal or low rank rule in \Cref{algo: pdtm}.}
\\
\uIf{$\rho_i \ge \eta$ and $\beta_i  \ge \alpha  $} 
{\textit{{Successful Step}}}
\uElseIf{$\rho_i \ge \eta$ and $\beta_i  \le -4\alpha  $ and $\sigma_d^{(i)}  \notin  \big[\frac{1}{6}, \frac{2}{3}\big] \big(-\beta_i + \alpha\big)\|s_d^{(i)}\|_W^{-1}$}{\textit{{Successful Step}}}
\uElseIf{$\rho_i \ge \eta$ and $\beta_i \in [-4\alpha, \alpha] $ and $
    \sigma_d^{(i)} \ge \frac{2}{3} \big(-\beta_i + \alpha\big)\|s_d^{(i)}\|_W^{-1} $}
{\textit{{Successful Step}}
\\ In all \textit{Successful Step}, update $s^{(i+1)} := s^{(i)} + s_d^{(i)} $, update $0 \le  d_{i+1} \le d_i$, update $\T_i$ with each diagonal entry $t^{(i)}_j[\lambda_{\min}(W)]^{-3/2} \le \hat{B} $, and  compute $g_{i+1}$. }
\Else{
\textit{Unsuccessful Step}: $s^{(i+1)} := s^{(i)}$,  $\beta_{i+1}:= \beta_{i}$ and increase $d_{i+1} := \gamma \max\{1, d_{i}\}. $
}
}
\label{algo: DTM variant 1} 
\end{algorithm}

In this variant, we introduce a positive parameter $\hat{B}$ that controls the magnitude of the third-order approximations $\T_i$,  preventing them from becoming excessively large.

\begin{assumption} 
For every iteration, $i$, 
the scalars $\{t^{(i)}_j\}_{1 \le j \le n}$ are uniformly bounded. Specifically, \\ $ \max_{i, 1 \le j \le n} |t^{(i)}_j|  [\lambda_{\min}(W)]^{-3/2}\le \hat{B} $.
\label{assumption1 DTM}
\end{assumption}
\begin{remark}
    Using \Cref{example special T} and \Cref{assumption1 DTM}, we obtain that  for every iteration, $i$, 
    $$
    \max_{u, v \in \R^n}\T_i[u][v]^2 \le  \max \{|t_1^{(i)}|, |t_2^{(i)}|, \dotsc, |t_n^{(i)}|\} [ \lambda_{\min}(W)]^{-3/2} \le \hat{B}\|v\|_W^2\|u\|_W. 
    $$
    In other words, \Cref{assumption1 DTM} implies that, for all $i$,  $|\T_i[s]^3|\le \hat{B} \|s\|_W^3$  and  $\|\T_i[s]^2\|\le \hat{B} \|s\|_W^2$  for any $s \in \R^n$. 
    \label{remark for assumption}
\end{remark}

In this variant, similar to the CQR Algorithm in \cite{zhu2023cubic}, we introduce an additional constant, $\alpha \in [0, \frac{1}{2}]$, to further control the step length $s_d ^{(i)}$. With this parameter, we account for two additional cases where the step is rejected and regularization is increased.  Define $\beta_i := \T_i \left[\frac{s_d^{(i)}}{\|s_d^{(i)}\|_W} \right]^3$. The first case occurs when $ \beta_i \le -4 \alpha$ and $\sigma_d$ falls within the narrow range $\big[\frac{1}{6}, \frac{2}{3}\big] \big(-\beta_i + \alpha\big)|s_d^{(i)}|^{-1}$. The second case arises when $\beta_i$ lies in the range $[-4 \alpha, \alpha]$ and $\sigma_d^{(i)} \ge \frac{2}{3} \big(-\beta_i + \alpha\big)|s_d^{(i)}|^{-1}$. For all other cases where $\beta_i$ and $\sigma_d^{(i)}$ do not match these scenarios, the standard procedure in Step 4 \Cref{TDTM} is sufficient to ensure the success of the step. More details can be found in \Cref{thm cubic Upper Bound DTM}. 

It is worth noting that whenever $\rho_i \ge \eta$,  even without the specific requirement for $\sigma_d^{(i)}$ in \Cref{algo: DTM variant 1}, one could ensure a decrease in the value of $m_3$ at each successful iteration (\Cref{remark cauchy analysis DTM}), such that
$
 m_3(s^{(i)}) -m_3(s^{(i+1)}) \ge \eta \big[m_3(s^{(i)}) - M_d(s^{(i)}, s_d^{(i)}) \big] >0
$
where $\eta \in (0, 1]$. Therefore, the Diagonal Tensor framework (\Cref{TDTM}) guarantees convergence to an approximate first-order local minimum of $m_3$. However, our goal in this section is to establish a global complexity bound; hence, we introduce the additional control parameter {(i.e., $\alpha, \hat{B}$ and updating rules for $\beta_i$)} in \Cref{algo: DTM variant 1}.

To establish the complexity bound for \Cref{algo: DTM variant 1}, we first derive an upper bound on the step size. Then, use the upper bound of step size to give an iteration-independent upper bound on the tensor term. The proof for the upper bound on step size is adapted from \cite[Thm 3.1]{zhu2023cubic}. The key modification is changing the regularization norm from $\|\cdot\|$ to $\|\cdot\|_W$ and substituting $ \big|\beta_i \|s_d^{(i)}\|^3 \big| \le \hat{B} \|s_d^{(i)}\|$ with $ |\T_i [s_d^{(i)}]^3| \le \hat{B} \|s_d^{(i)}\|_W^3$ , more details can be found in \Cref{thm upper bound on step size DTM} and \Cref{corollary upper bound for Ti DTM}.

Next, we establish a lower bound for the step size.  The lower bound on step size is essential for ensuring a reduction in the values of $m_3$ (\cite[Lemma 2.3]{birgin2017worst} and \cite[Lemma 3.3]{cartis2020concise}),  as it shows that the step cannot diminish to an arbitrarily small size relative to the criticality conditions ($\|\nabla m_3(s^{(i)} + s_d^{(i)})\| > \epsilon$ and $\|g_i\| > \epsilon$), for both successful and unsuccessful iterations. With  slight modification by replacing $ \big|\beta_i \|s_d\|^2 \big|\le \hat{B} \|s_d\|^2 $ with $ \|\T_i [s_d]^2\| \le \hat{B} \|s_d\|_W^2$, 
techniques in \cite[Lemma 3.1]{zhu2023cubic} give that the lower bound on the step is 
$$
\|s_d^{(i)}\|_W >   \min\bigg\{(\hat{B}  +L_H )^{-1/2}\epsilon^{1/2}, \frac{1}{2}{d_i}^{-1/3}  \epsilon^{1/3} \bigg\}
$$
where $\epsilon$ represents the prescribed first-order optimality tolerance

Moreover, the upper and lower bounds for $ d_i $ are
$$
        d_i > \frac{\hat{B} + L_H}{6}\|s_d^{(i)}\|_W^{-1}>0, \qquad  d_i \le d_{\max} :=  \gamma (\hat{B} + L_H)^{3/2}\epsilon^{-1/2}.
 $$
The full proofs are given in \Cref{lemma: lower bound of norms DTM}, \Cref{thm cubic Upper Bound DTM} and \Cref{lemma upper bound on di DTM} for completeness. 
The rationale behind upper bound  for $ d_i $ is that whenever the regularization parameter surpasses $d_{\max}$, 
 $M_d(s^{(i)}, s)$ overestimates the value of $m_3$ at the minimizer of $M_d$ (as stated in \Cref{thm cubic Upper Bound DTM}). This overestimation implies that any reduction achieved in minimizing $M_d(s^{(i)}, s)$ is as a guaranteed lower bound on the decrease in $m_3$, thereby ensuring that each iteration is successful and does not need to further increase the regularization parameter.

\begin{remark}[Distinction between \Cref{TDTM} and \Cref{algo: DTM variant 1} ] \Cref{TDTM} provides a simplified and more general framework, and it is also used in the final implementation. In contrast, \Cref{algo: DTM variant 1} is introduced to establish the global complexity bound. To this end, additional control parameters (namely \(\alpha\), \(\hat{B}\), and the updating rules for \(\beta_i\)) are incorporated in \Cref{algo: DTM variant 1}. In particular, \(\alpha \in [0, \tfrac{1}{2}]\) is introduced to further regulate the step length \(s_d^{(i)}\). 
 \Cref{algo: DTM variant 1} is globally convergent and also practical. Rather, in practical implementations, it is unnecessary to distinguish between the multiple subcases introduced in \Cref{algo: DTM variant 1}. These subcase distinctions are primarily required for the technical proof of the complexity analysis.
\end{remark}

 Next, we show that the diagonal tensor model constructed in successful steps of \Cref{algo: DTM variant 1} yields a decrease in the value of $m_3$ by at least $\mathcal{O}(\epsilon^{3/2})$. 

\begin{theorem}
\label{thm M3 value decrease DTM}  
\textbf{(Bounding the $m_3$ decrease  on successful iterations\footnote{The proof of \Cref{thm M3 value decrease DTM}  is adapted from \cite[Thm 3.3]{zhu2023cubic} and provided in \Cref{appendix proof of complexity thm}.})}
Suppose that  \Cref{assumption1 DTM} holds and \Cref{algo: DTM variant 1}  is used. For $\alpha \in (0, \frac{1}{2})$ and $\eta \in (0, 1]$, assume that $i+1$ is a successful iteration {which does not meet the termination condition, i.e., $\|g_{i+1}\|> \epsilon$ and $\|g_i\| > \epsilon$}. Then,
\begin{equation}
  m_3(s^{(i)})-  m_3(s^{(i+1)})
 \ge  \frac{\alpha \eta}{24} \|s_d^{(i)}\|_W^3 \ge \kappa_s \epsilon^{3/2}
  \label{general model decrease DTM}
\end{equation}
where {$\kappa_s = \frac{\alpha \eta}{192\gamma}(\hat{B} + L_H)^{-3/2}$. } 
\end{theorem}

Lastly, we establish the bound for unsuccessful iteration and total iterations, thereby giving the complexity bound for \Cref{algo: DTM variant 1}.

\begin{theorem} \textbf{(Complexity bound for \Cref{algo: DTM variant 1}\footnote{The proof of \Cref{thm: dtm complexity} is adapted from \cite[Thm 3.4]{zhu2023cubic} and provided in \Cref{appendix dtm complexity}.}}) {Let $m_{\text{low}}$ be a lower bound on $m_3(s)$ for $s \in \R^n$.} Suppose that  \Cref{assumption1 DTM}  holds.  Then, there exists a positive constant $\kappa_s$ (see \Cref{thm M3 value decrease DTM}),  $\kappa_c$ such that \Cref{algo: DTM variant 1} requires at most
$$
 \kappa_s \frac{m_3(s^{(0)})-m_{\text{low}}}{\epsilon^{3/2}}  + \kappa_c + \log{\epsilon^{-1/2}} = \mathcal{O}(\epsilon^{-3/2}). 
$$
function evaluations (i.e., $f_i$) and at most
$$
 \kappa_s \frac{m_3(s^{(0)})-m_{\text{low}}}{\epsilon^{3/2}}  + 1
= \mathcal{O}(\epsilon^{-3/2}). 
$$
derivatives evaluations (i.e., $g_i, H_i, T_i$) to compute an iterate $s_\epsilon$ such that $\|g_i\| \le \epsilon$ {or $\|g_{i+1}\| \le \epsilon$} .
\label{thm: dtm complexity}
\end{theorem}

\textbf{Improved Complexity for Special Cases: }
Although Algorithm \ref{algo: DTM variant 1} achieves a function value reduction of $\mathcal{O}(\epsilon^{3/2})$ for a general nonconvex $m_3$, in certain special cases, the DTM algorithm or its variants exhibit improved complexity bounds during specific iterations.  One such instance is when the \textit{tensor term is sufficiently small}. If $\|T_i\| \le \epsilon^{1/3}$, and we entries of $\T_i$ with absolute value less than  $\epsilon^{1/3}$, then the decrease in value of $m_3$ in such an iteration is $\mathcal{O}(\epsilon^{4/3})$. The analysis for this case follows similarly as in \cite[(eq 48)]{zhu2023cubic}. 

Secondly, if $m_3$ has a \textit{diagonal tensor,}  we can set $\T_0 = T$ and $d_0 =0$, \eqref{M_d} exactly recovers \eqref{ar3 model}. As a result,  \Cref{TDTM} converges in one iteration.     
Lastly, \Cref{thm change of basis dtm} indicates that if \eqref{ar3 model} has \textit{a tensor with rank $P \le n$} as illustrated in \eqref{rank P dtm}, then by setting $\g$, $\H$, the diagonal tensor $\T$ and $\W$ as indicated in \Cref{thm change of basis dtm}, we can minimize $ m_3^{\mathbbm{P}}(s)$ by minimizing the following diagonal tensor model
    \begin{eqnarray}
      \s_* =\min_{\s \in \R^n}  f_0 + \g^T  \s  + \frac{1}{2}\s ^T \H \s  +\frac{1}{6} \T[\s]^3  + \frac{\sigma}{4}\big\|\s \big\|_{\W}^4 
    \end{eqnarray}
and recover the minimizer by $s_* = \min_{s\in \R^n}  m_3^{\mathbbm{P}} =C^{-1}\s_*.$
A slight modification to \Cref{TDTM} can yield an algorithm for minimizing \eqref{ar3 model} with a tensor rank $P \le n$. By adjusting the inputs of \Cref{TDTM} as follows: set $d_0=0$, $f_0 := m_3(s^{(0)})$, $g_0 := C^{-T} \nabla m_3(s^{(0)})$, $H_0 := C^{-T} \nabla^2 m_3(s^{(0)}) C^{-1}$, and $\W := C^{-T}WC^{-1}$. Additionally, modify Step 2 as $s_d^{(i)} = C^{-1} \argmin_{s \in \R^n} M_d(s^{(i)}, s)$. With these adjustments, \Cref{TDTM} converges to the minimizer of $m_3$ in a single iteration, $i$. The modified variant of \Cref{TDTM} for minimizing $m_3^{\mathbbm{P}}$ with low-rank tensor is outlined in \Cref{algo: pdtm} with numerical examples in \Cref{sec: numerical testing dtm}.

\section{Numerical Implementation and Preliminary Results}
\label{sec numerics}

In this section, we offer two options for updating $\T_i$ in \Cref{algo: pdtm}, thereby completing the structure of \Cref{TDTM}. The first option for updating $\T_i$ is suitable for a general $m_3$ or $m_3$ with a diagonal tensor. The second option for updating $\T_i$ is suitable for $m_3^{\mathbbm{P}}$ with a low-rank tensor (rank $P \leq n$) and converges to the minimizer of $m_3^{\mathbbm{P}}$ in a single iteration.

\begin{algorithm}
\caption{\small An Implementable Diagonal Tensor Method (DTM) Variant for Minimizing $m_3(s)$}
Set a fixed bound $\hat{B} := 10^6$. Initialize and proceed through Algorithm \ref{TDTM} with specifications as follows. 
In Steps 2 and 4, update $\T_i$ and $W$ using one of the following rules:
\begin{enumerate}
    \item \textbf{Diagonal Rule}: Set $\T_i = \diag\{t^{(i)}_j\}_{1 \le j \le n}$, where $t^{(i)}_j = T_i[j,j,j]$ are the diagonal entries of $T_i$. 
    \\If $T_i[j,j,j] > \hat{B}$, set $t^{(i)}_j := \sign(T_i[j,j,j])\hat{B}$.
    % {\mbox{Set $M_d$ as in \eqref{M_d} and compute $s_d^{(i)} = \argmin_{s \in \R^n} M_d(s^{(i)}, s)$ using \Cref{Newton algo dtm}.}}
    \item \textbf{Low Rank Rule}: Compute a rank-$P$ approximation of $T_i$, denoting $\hat{C} := \begin{bmatrix} \boldsymbol{a}^{(1)}, & \dotsc, & \boldsymbol{a}^{(P)} \end{bmatrix}^T$. Set %$C$ and $\T$ as in \eqref{C matrix dtm}--\eqref{change of variable dtm}
    $$
     C := \begin{bmatrix} \hat{C}, & \hat{C}_\perp \end{bmatrix}^T, \quad \g_i = C^{-T} g_i, \quad \H_i = C^{-T} H_i C^{-1},  \quad \W_i = C^{-T}WC^{-1}, \quad \T_i = \diag\{t_j\}_{1 \le j \le n}
    $$
    % where $t_1 = \dotsc = t_P =1$ and $t_{P+1} = \dotsc = t_n = 0$. Define
    % $
    %  \tilde{M}_d(s^{(i)}, \s)  = f_0 + \g_i^T  \s  + \frac{1}{2}\s ^T \H_i \s  +\frac{1}{6} \T_i[\s]^3  + \frac{\sigma}{4}\big\|\s \big\|_{\W_i}^4.
    % $ 
    Recover $s_d^{(i)} = C^{-1} \s_d^{(i)}.$
\end{enumerate}
\label{algo: pdtm}
\end{algorithm}

\begin{remark} (Role of $\hat{B}$)
In Algorithm \ref{algo: pdtm}, the parameter $\hat{B}$ prevents the diagonal terms $|t_j^{(i)}|$ from becoming excessively large, thereby preventing potential numerical issues. 
\label{remark bound B dtm}
\end{remark}

\begin{remark} (Tensor-free nature) {Note that in \Cref{algo: pdtm}, $T_i$ is the full third-order derivative, represented as an $n \times n \times n$ tensor. 
$\mathcal{T}_i$ is the diagonal tensor used in the DTM to approximate $T_i$, and it has a sparse structure.} 
The Diagonal Rule is tensor-free as we only access the $n$ diagonal elements of $T_i$. All tensor-vector, tensor matrix multiplication is evaluated using $\mathcal{O}(n)$ operations, $\T_i[s] = \textbf{t}^{(i)} \cdot s$, $\T_i[s]^2 = \textbf{t}^{(i)} \cdot s \cdot s$, and $\T_i[s]^3 = \textbf{t}^{(i)} \cdot s \cdot s \cdot s$, where $\cdot$ represents the entrywise dot product.
The Low-Rank Rule obtains the  rank-$P$ approximation of $T_i$ using \texttt{Eigenproblems} Functionality in Tensor Toolbox \cite{bader2006algorithm, bader2008efficient}. Besides the step of extracting low-rank approximation, the rest of the Low-Rank Rule also does not involve tensor operations. 
\end{remark}

{We give the flow chart of the DTM framework in \Cref{fig:dtm-flowchart}. \Cref{TDTM} is the main DTM framework. \Cref{TDTM} calls \Cref{algo: pdtm} as a subroutine to update the diagonal tensor~$\mathcal{T}_i$, and calls \Cref{Newton algo dtm} as a subroutine to compute the step by minimizing the DTM model. \Cref{algo: DTM variant 1}, by contrast, is a more detailed DTM variant used exclusively for establishing worst-case global evaluation complexity bounds in \Cref{sec: Convergence and Complexity of dtm}.}
%The framework of \Cref{algo: pdtm} (i.e., \Cref{TDTM}) ensures a reduction in the value of $m_3$ at each successful iteration. 

\begin{figure}[ht]
\centering
\caption{\small \textbf{Flowchart of  AR3 minimization algorithmic framework.} The implementable variant (\Cref{TDTM}) and the theoretical variant (\Cref{algo: DTM variant 1}) are shown, together  with subroutines for step computation (\Cref{algo: pdtm}) and tensor updates (\Cref{Newton algo dtm}). The framework takes initial data and returns an approximate minimizer of $m_3$ satisfying first-order optimality conditions with tolerance level $\epsilon$.}

\begin{tikzpicture}[
  node distance=1.7cm and 2.8cm,
  >=Stealth,
  every node/.style={font=\small},
  block/.style={rectangle, draw, rounded corners, align=center, minimum width=3.6cm, minimum height=1cm},
  io/.style={trapezium, trapezium left angle=70, trapezium right angle=110,
             draw, align=center, minimum width=3.6cm, minimum height=1cm},
  line/.style={draw, -{Stealth}},
  dashedblock/.style={rectangle, draw, rounded corners, dashed, align=center, minimum width=3.6cm, minimum height=1cm}
]

% Nodes for main DTM framework (\Cref{TDTM})
\node[dashedblock, io] (start) {\textbf{Initial data}\\ $x_0, H_0, T_0, f, \epsilon \dotsc$};
\node[block, fill=yellow!20, below=1cm] (alg2) {
\begin{tabular}{c | c}
\textbf{Implementable Variant} & \textbf{Theoretical Variant} 
 \\
\Cref{TDTM} & \Cref{algo: DTM variant 1} \\
\footnotesize (practically minimize AR3) & \footnotesize (for worst-case complexity)
\end{tabular}
}; 
\node[block, fill=blue!20, below left=0.7cm and 0.7cm of alg2] (alg4) {\Cref{algo: pdtm}:\\ Update diagonal tensor $\mathcal{T}_i$};
\node[block, fill=blue!20,  below right=0.7cm and 0.7cm of alg2] (alg1) {\Cref{Newton algo dtm}:\\ Compute step $s^{(i)}_d$};
\node[dashedblock, io, below=2cm of alg2] (end) {\textbf{Return Output:}\\ an approx. min. of $m_3$\\ s.t  $\|\nabla m_3(s^{(i)})\| \le \epsilon$};

% Main arrows
% \path[line] (start) -- (alg2);
\path[dotted]  (alg4)  -- (alg2)
node[pos=0.25, above,  rotate=15] {\footnotesize 
subroutine};
\path[dotted] (alg1) --  (alg2) 
node[pos=0.25, above,  rotate=345]  {\footnotesize 
subroutine};

\end{tikzpicture}
\label{fig:dtm-flowchart}
\end{figure}

\subsection{\texorpdfstring{Minimizing $m_3$ with Diagonal, Low Rank and Full Tensor}{Minimizing Diagonal and Low Rank m3}}
\label{sec: numerical testing dtm}

The natural test cases for {the DTM algorithmic framework (\Cref{TDTM}) with diagonal updates or low-rank updates (\Cref{algo: pdtm} as subroutine)} are, respectively, the minimisation of \eqref{ar3 model} with diagonal tensor terms and the minimisation of \eqref{ar3 model} with low-rank tensor terms.
Therefore, in this subsection, we test the DTM algorithm to these test sets, along with a standard test set of $m_3$ with a full-rank tensor term.

\begin{enumerate}
    \item \textbf{Diagonal tensor test set}: The objective polynomial is  
$
m_3(s)  = g^T s + \frac{1}{2} H [s]^2  +  \frac{1}{6} \sum_{j=1}^n t_j s_j^3 + \frac{\sigma}{4} \|s\| ^4 
$
where $g$, $H$, and $t_j$ are generated as follows
\begin{eqnarray}
g = \texttt{a*randn(n, 1)}, \quad  H = \texttt{b*symm(randn(n, n))},  \quad  t = [t_1, \dotsc, t_n]^T = \texttt{c*randn(n, 1)}.
\label{set 1}
\end{eqnarray}
    \item
\textbf{Low rank test set}: The objective polynomial is 
$
m_3(s)  = g^T s + \frac{1}{2} H [s]^2 + \frac{1}{6} \sum_{k=1}^P ({\boldsymbol{a}^{(k)}}^T s)^3 + \frac{\sigma}{4}\|s\|^4. 
$
Specifically, $g$, $H$ are generated the same as the Diagonal tensor test set and $\boldsymbol{a}^{(k)} = \texttt{c*randn(n, 1)}$
for $j = 1, \dotsc, p$.
\item 
\textbf{AR$3$ subproblem test set}: The objective polynomial is the same as \eqref{ar3 model}. 
Specifically, $g$, $H$ are generated the same as the Diagonal tensor test set, and $ T = \texttt{c*symm(randn(n, n, n))}$ is generated using the Tensor Toolbox \cite{bader2006algorithm, bader2008efficient}.
\end{enumerate}
Here, \texttt{n} represents the dimension of the problem ($s \in \R^n$), and \texttt{symm(randn())} denotes a symmetric matrix with entries following a normal distribution with mean zero and variance one. The parameters \texttt{a}, \texttt{b}, \texttt{c}, and $\sigma$ are selected differently to test the algorithm's performance under various scenarios. We occasionally assign specific values to $H$ and $t$ to test the algorithm’s performance on ill-conditioned Hessians and tensors. These specifications are detailed in the numerical experiments. All tests were conducted on in MATLAB R2023b on an Intel(R) Core(TM) i7-4770 CPU @ 3.40GHz   3.40 GH processor with 16 GB of RAM. The code is written in double precision. Unless otherwise stated, the parameters are set at $\eta = 0.1, \eta_1 = 2, \gamma_2 = 0.5, \gamma = 2.0.$ In our numerical experiments, {unless otherwise specified, we set the tolerance and stopping criterion as $\|\nabla^2 m_3(s_\epsilon) \| \le \texttt{TOL}$, where $\texttt{TOL} = 10^{-4}$} ($\epsilon_l = 10^{-8}$, $\epsilon_{\kappa} = 10^{-7}$). After obtaining $s_\epsilon$, we tested both first-order and second-order criticality conditions ($\nabla^2 m_3(s_\epsilon) \succeq - \texttt{TOL}$). In all our test examples, the algorithm converged to an approximate second-order minimum of \eqref{ar3 model}.

\begin{figure}[ht]
    \centering
\caption{\small \textbf{Diagonal Tensor Test Set}: Averages over 10 random trials are plotted, with error bars showing one standard deviation. The same legend applies to all subplots.  Standard Set: \texttt{a = 20, b = 20, c = 20}, $\sigma = 100$ in  \eqref{set 1}.
\\ Ill-Conditioned Hessian Set: \texttt{a = 10, c = 5}, $H$ is diagonal matrix uniformly distributed in $[10^{-6}, 10^3]$, $\sigma = 500$ in  \eqref{set 1}.
\\ Ill-Conditioned Tensor Set: \texttt{a = 10, b = 20}, $t$ has entries uniformly distributed in $[10^{-6}, 10^2]$, $\sigma = 500$ in  \eqref{set 1}. }
    \label{fig: diagnonal}
    \includegraphics[width =5cm, height =6cm]{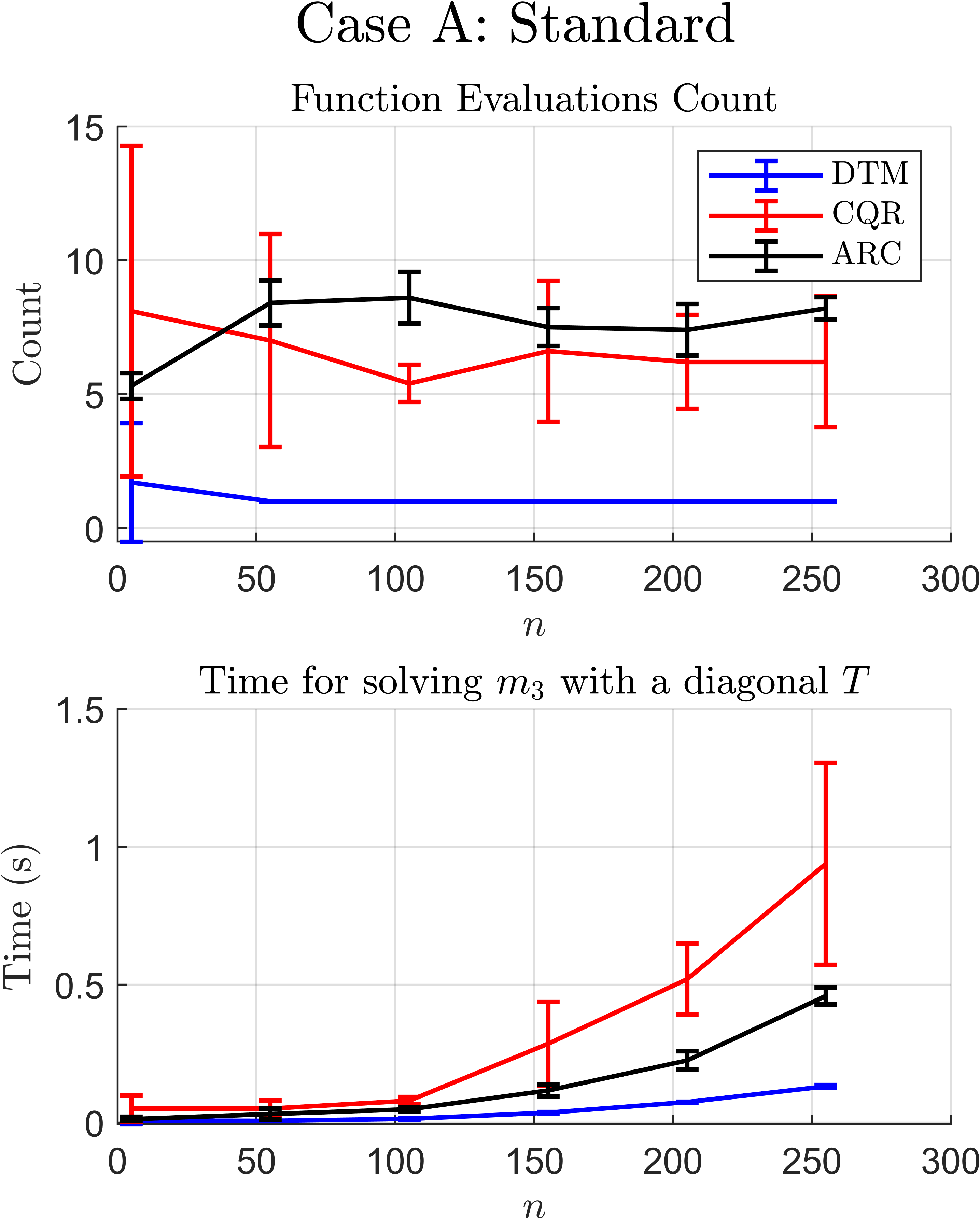}
    \includegraphics[width =5cm, height =6cm]{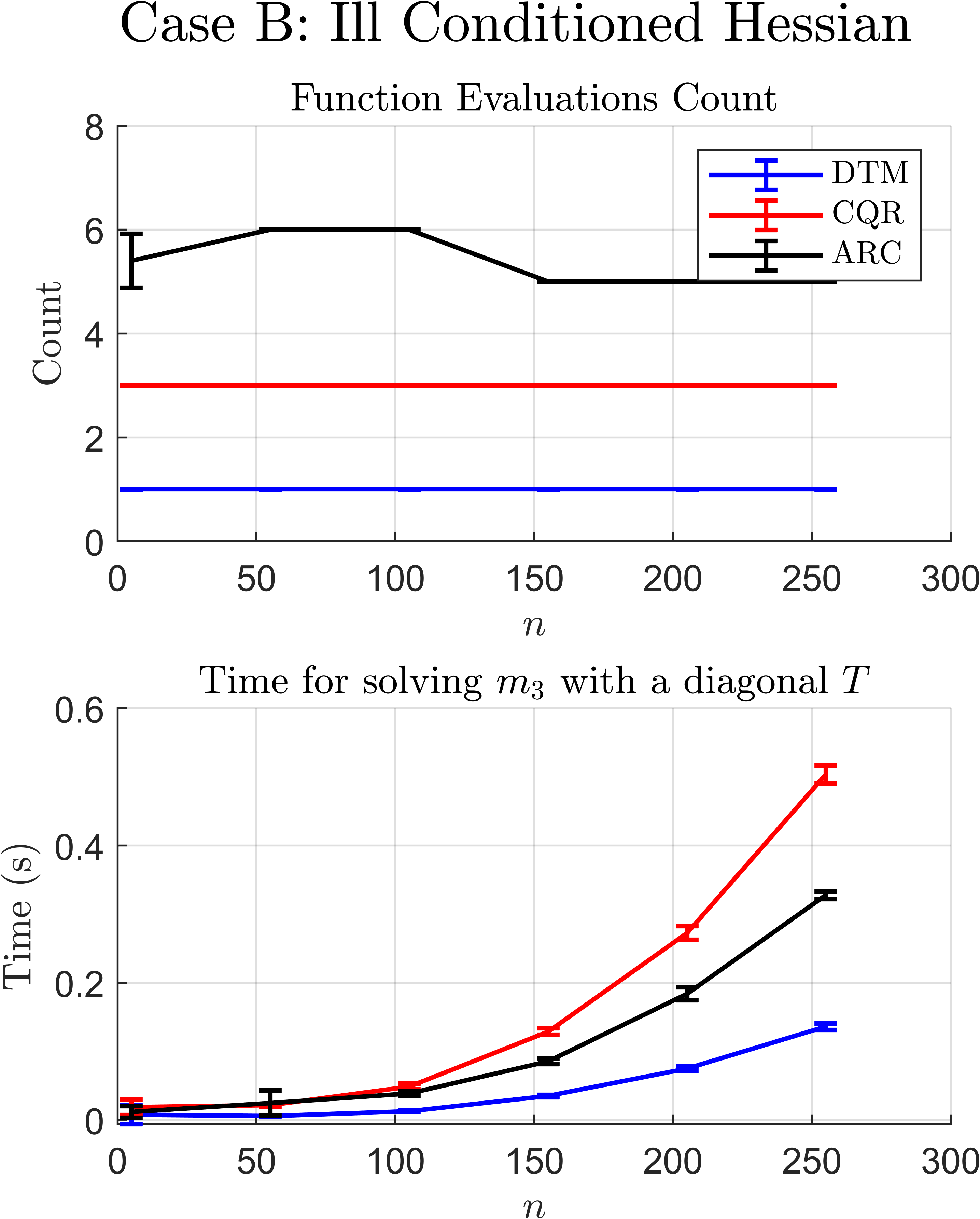}
    \includegraphics[width =5cm, height =6cm]{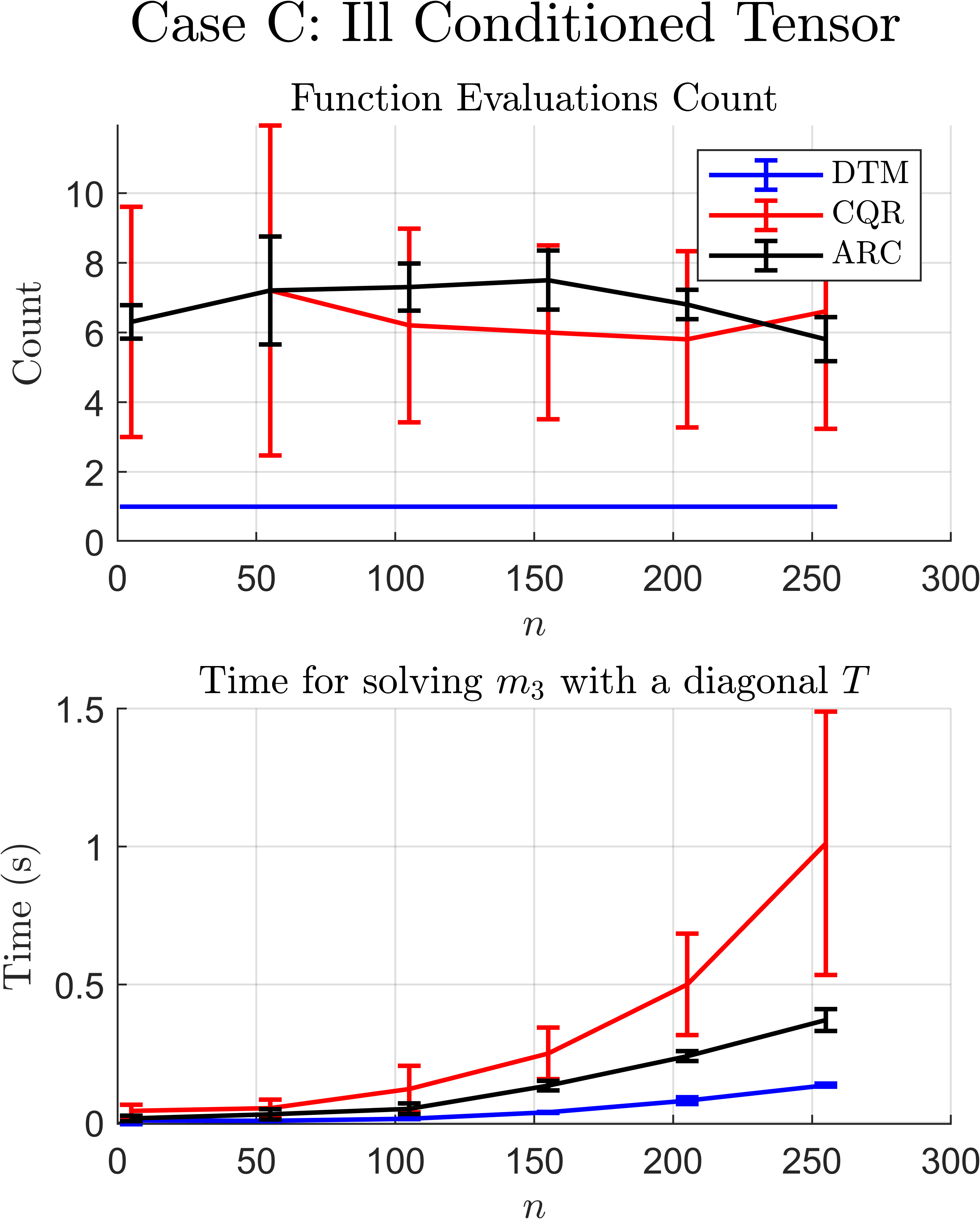}
\label{fig diagonal test set}
\end{figure}

\paragraph{Diagonal tensor test set:} 
In \Cref{fig diagonal test set}, we tested the {DTM algorithmic framework (\Cref{TDTM}) with diagonal updates (\Cref{algo: pdtm} as subroutine)} on the diagonal tensor test set for problem dimensions up to 250 and compared it to the state-of-the-art {Adaptive Regularization with Cubics (ARC) method \cite{cartis2011adaptive}\footnote{The code for the ARC algorithm can be found at URL: \url{https://www.galahad.rl.ac.uk/packages/}} and
 Cubic Quartic Regularization (CQR) method \cite{zhu2023cubic}}. 
 We tested three subsets of diagonal tensor in \eqref{set 1}: (A) randomly generated $g$, $H$, and $t$; (B) ill-conditioned Hessian; and (C) ill-conditioned tensors. For each subcase, 10 randomized trials were performed, and the averages across these trials {with error-bar} are displayed in \Cref{fig diagonal test set}. 
Since $M_d$ exactly approximates our test examples, DTM algorithm consistently converges in a single iteration (requiring one evaluation of the function and derivatives), while the ARC and CQR algorithm typically requires 5--10 function evaluations. 
{Also, DTM achieves the lowest computational time, particularly for large dimensions. We also observed that, for ARC (and CQR), more than $50\%$ of the total time is spent evaluating the function and its derivatives. Thus, since ARC (and CQR) requires $5$--$10$ function and derivative evaluations while DTM requires only one, this is the main reason its computational time becomes larger than that of DTM, especially for larger $n$.} In terms of the total number of Cholesky factorizations, DTM outperforms ARC in cases A and C, and ARC outperforms DTM in case B. In all cases, the second-order local optimality condition is satisfied, as well as the necessary global optimality condition outlined in \Cref{thm necessary tight}. 
{Among the three methods, CQR usually takes more evaluations than DTM but fewer than ARC. This behaviour is expected. In CQR, the tensor term is approximated using a single parameter~$\beta$. When the tensor has only one dominant diagonal entry---or when a single direction in the tensor is large---CQR performs well. However, when there are $n$ significant diagonal entries, a single~$\beta$ cannot adequately capture all tensor directions, and consequently, more function evaluations are needed. More details for the CQR method can be found in~\cite{zhu2023cubic}. The need to more effectively capture multiple tensor directions was one of the motivations behind developing DTM as a more generalized version of CQR. }{We also experimented with a tighter tolerance, such as $10^{-9}$, and the same behaviour was observed.}

\begin{figure}[ht]
    \centering
\caption{\small \textbf{Low-Rank Tensor Test Set}: Parameter configurations are as follows. \texttt{a = 10, b = 20, c = 20}, $\sigma = 400$, with rank $P = 1$ and $P = 4$, respectively. {The change of basis is implemented via the Low-Rank Rule for DTM in \Cref{algo: pdtm} and calculated into the computation time. No change of basis for other algorithms, as their algorithmic frameworks do not naturally adapt this feature.}}
    \label{fig: low-rank}
    \includegraphics[width =5cm]{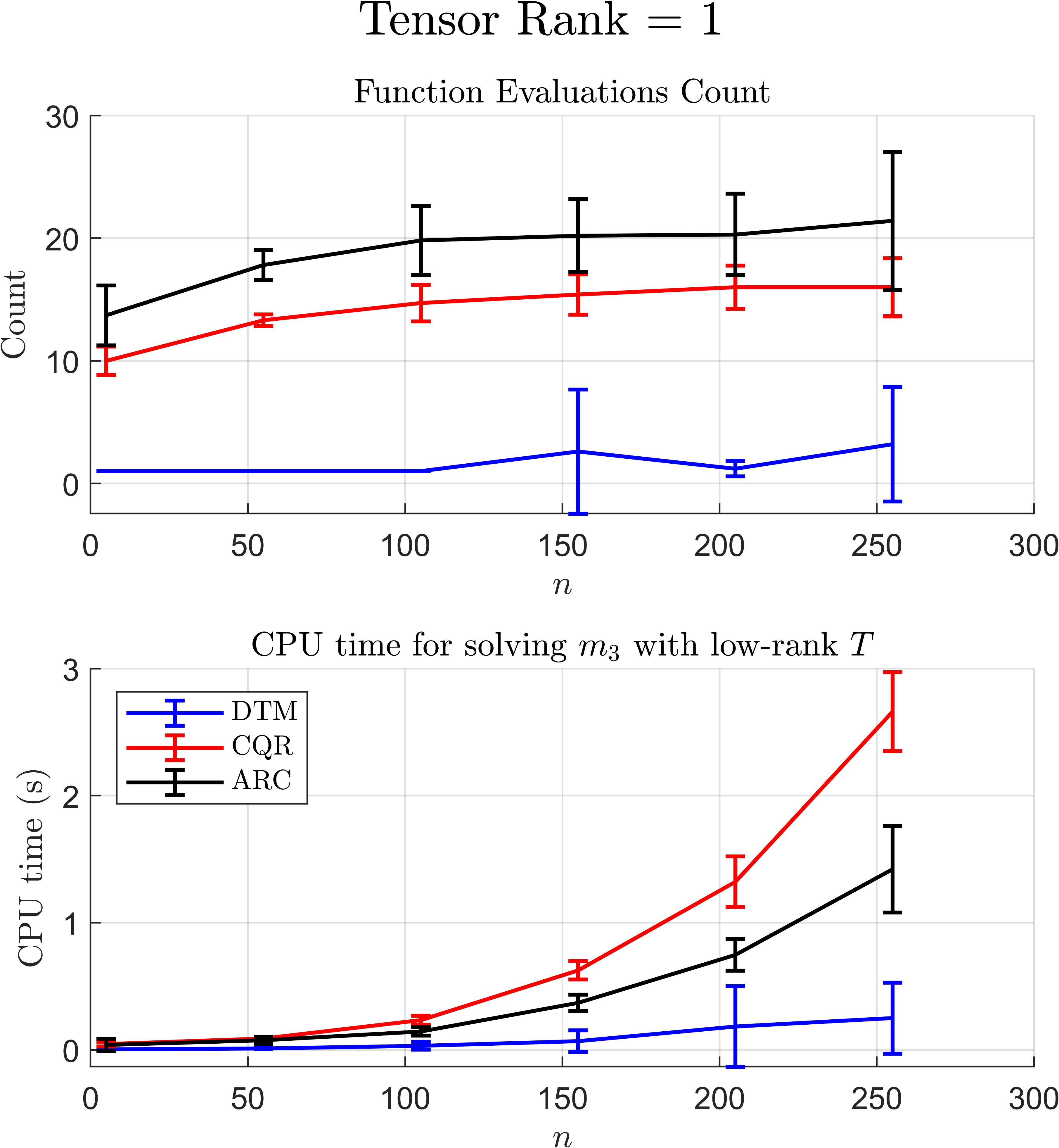}
    \includegraphics[width =5cm]{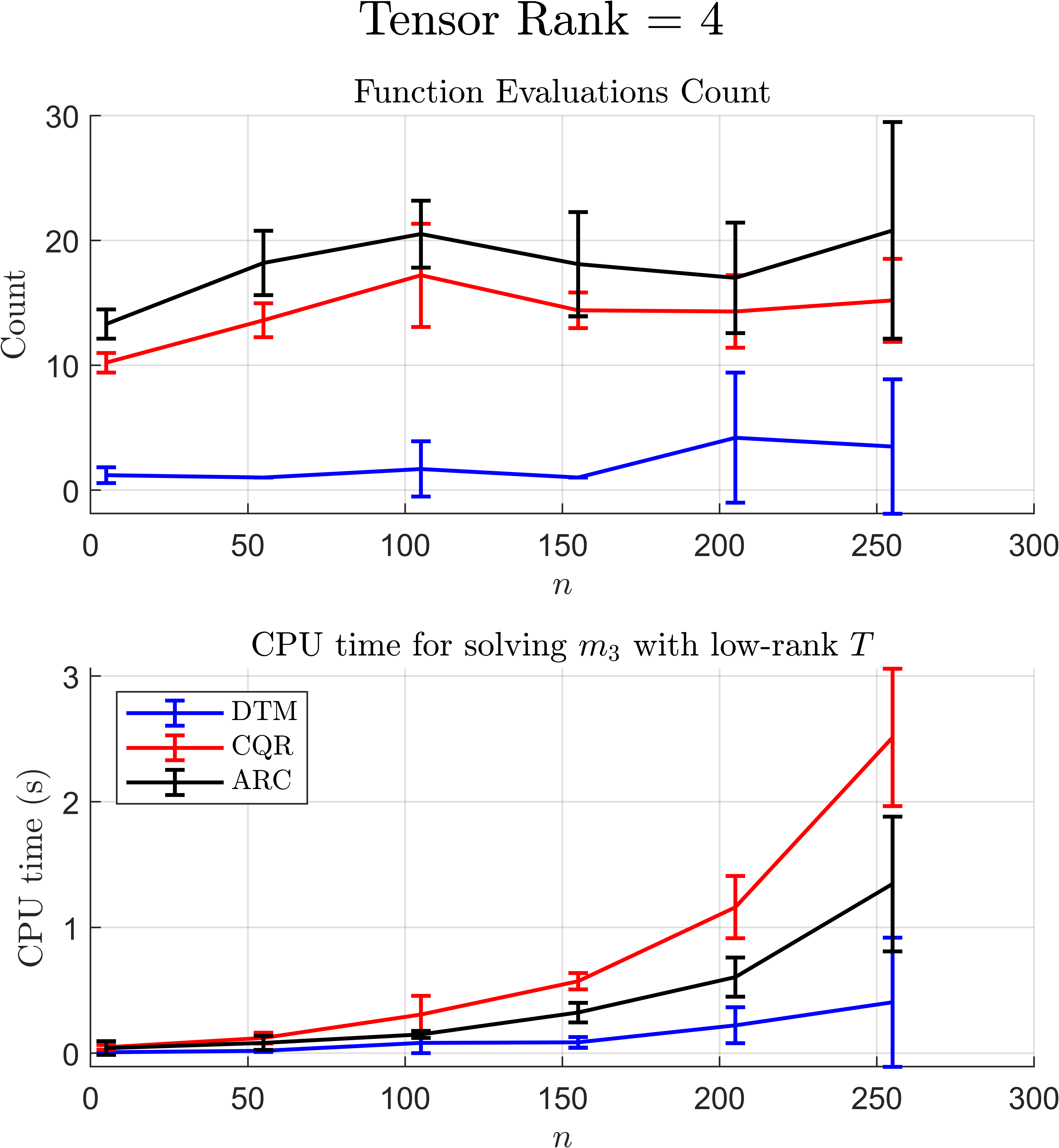}
\label{fig low-rank}
\end{figure}

\paragraph{Low tensor rank test set:} 
In \Cref{fig low-rank}, we tested {the DTM algorithmic framework (\Cref{TDTM}) with low-rank updates (\Cref{algo: pdtm} as subroutine)} on the low-rank tensor test set for tensor ranks $P = 1$ and $P = 4$. We compare DTM with both ARC algorithm and the CQR algorithm. 
Similar performance patterns were observed, with DTM requiring significantly fewer function evaluations and computation time than ARC/CQR in all cases. We give a brief justification of why this is the case. 
{When the tensor has $P = 1$ or $P = 4$, we use the change of basis described in Theorem~3.1:
\[
T^{\mathbbm{1}}[s]^3 
= \left(a^{(1)}_1 s_1 + \dotsb + a^{(1)}_n s_n\right)^3
\quad \underset{\text{change of basis}}{\leadsto} \quad 
\sum_{j=1}^P \tilde{s}_j^3.
\]
Under this transformation, the tensor essentially becomes a sparse diagonal tensor with a few nonzero entries, $T_{1:P,1:P,1:P} = 1$ where $T_{i,j,k}=0$  for all other indices. Consequently, the DTM subproblem can be solved in essentially one iteration, matching the
behaviour observed when the tensor is already diagonal (\Cref{fig: diagnonal}). Yet, the computation time is slightly longer than in the diagonal test set, because our computation time includes the cost of performing the change of variables. We notice that for \(n \le 100\), the performance of DTM and ARC is largely comparable, but for \(n > 100\), a clear distinction emerges. This is because, when \(n\) is small, the Cholesky factorization dominates the computation time. While for larger \(n\), the cost of function and derivative evaluations dominates. Therefore, the advantage of DTM---performing a single change of basis and solving the transformed problem in one or only a few iterations---becomes pronounced.}
Lastly, we note that in some DTM iterations, when solved to high accuracy, \Cref{Newton algo dtm} may encounter numerical issues due to the norm $\|\cdot\|_W$ being ill-conditioned, prompting us to revert to the ARC algorithm as a safeguard. This results in DTM converging in slightly more iterations, typically $3$-$4$ iterations. Further improvements in \Cref{Newton algo dtm} to enhance scalability are deferred to future work.

\paragraph{AR$3$ subproblem test set:} 
In \Cref{table:results}, we tested {DTM algorithmic framework (\Cref{TDTM}) with diagonal updates (\Cref{algo: pdtm} as subroutine)} on the full-rank AR$3$ subproblem with no dominant tensor direction and compared it with the ARC and CQR methods.   DTM performed competitively with ARC in terms of both iterations and function evaluations, although it required slightly more iterations and evaluations than the CQR method. Generally, if the quadratically regularized polynomial is low-rank or has a diagonal tensor term, DTM is preferred. On the other hand, for polynomials with a full-rank tensor and no dominant tensor direction, ARC or CQR is preferred.
Note that in this paper, we only provide a preliminary numerical illustration of applying the DTM method for minimizing the AR$3$ subproblem, particularly those with low-rank or diagonal tensor terms. The practical implementation of the AR$3$ methods to minimize the objective function is work in itself, and we refer the reader to \cite{cartis2024Efficient} for a full implementation of the AR$3$ algorithm. In \cite{cartis2024Efficient}, the AR$3$ solver is tested on the Moré, Garbow, and Hillstrom Test Problem Set, and acceleration techniques such as pre-rejection and pre-acceptance are applied to accelerate the convergence of the AR$3$ solver.

\begin{table}[ht]
\caption{\small \textbf{AR$3$ Subproblem test set}: \texttt{a = b = c = 80, TOL = $10^{-3}$},  stopping criterion as $\|\nabla^2 m_3(s_\epsilon) \| \le \texttt{TOL}$. The results are averaged over 10 randomly generated test functions \Cref{TDTM} with the Diagonal Rule in \Cref{algo: pdtm}.}
\centering
\small
\begin{tabular}{|c|c c c c c|c c c c c|}
\hline
\multirow{2}{*}{Method} & \multicolumn{5}{c|}{Successful Iterations} & \multicolumn{5}{c|}{Total Iterations} \\
& $n=15$ & $n=25$ & $n=50$ & $n=75$ & $n=100$ & $n=15$ & $n=25$ & $n=50$ & $n=75$ & $n=100$ \\
\hline
\textbf{ARC} & 10.5 &   12.9   & 15.4    &15.1   & 18.3   &  8.6  &  10.9  &  13.4    &13.1   & 15.9 \\
\textbf{CQR} & 7.3   &   11  &  10.8  &  11.6 &   13.5  &   6.6  &   8.9   &  9.8 &   10.4   &   12  \\
\textbf{DTM}  & 8.1 &    9.5  &  11.6  &  12.7   & 15.7&     8.1   &  9.5   & 11.6   & 12.7    &15.7 \\
\hline
\end{tabular}
\label{table:results}
\normalsize
\end{table}

\paragraph{Magnitude of $\|T\|$ or $\sigma$ and Global Optimality Characterization:} In \Cref{fig: increase t or sigma} in \Cref{appendix: size}, we examine how the magnitude of $\|T\|$ or $\sigma$ influences global optimality characterization. \Cref{fig increase sigma} illustrates that as $\sigma$ or $\|T\|$ increases, the likelihood that the minimum generated by the DTM algorithm satisfies the global sufficient condition in \Cref{thm sufficiency tight} also increases, indicating that the global minimum of $m_3$ is achieved. This finding aligns with the theoretical results in \Cref{sec: gap global Necessary and sufficient}, which demonstrate that the global necessary and sufficient conditions for the AR$3$ subproblem coincide when $\sigma$ is sufficiently large or $\|T\|$ is sufficiently small.

\paragraph{Applying \Cref{TDTM} directly to $m_3$:}
As discussed in \Cref{Sec Iterative Algorithm for Minimizing AR3 Model}, the tensor-free structure of the DTM offers advantages in terms of computational time and storage. 
Note that \Cref{TDTM} can be set directly with $\T_i = T_i$, where $T_i$ denotes the full third-order derivative rather than a diagonal tensor. 
In \Cref{fig computational time} in \Cref{appendix: time}, we compare the computational time of \Cref{TDTM} with $\T_i = T_i$ (the full third-order derivative) against the DTM framework where $\T_i$ is a sparse diagonal tensor. As the dimension increases, the variant with $\T_i = T_i$ (the full third-order derivative) becomes significantly slower. 
This slowdown is primarily due to the costly tensor--vector operations required by the full third-order derivative when neither tensor-free nor sparsity structures are employed.
% These numerical results highlight that, for a tensor method to remain computationally feasible, adopting a tensor-free framework is often essential.

% Potential numerical experiment to do: 
% Karl's two high-order examples.
% tensor free or not tensor free difference. 

{

\subsection{Numerical Examples with Moré–Garbow–Hillstrom test set}
\label{sec  MGH Test Set examples}

We test our method on a subset of the Moré–Garbow–Hillstrom (MGH) test set~\cite{more1981testing} and provide preliminary numerical observations in this subsection. All problem initializations, as well as the first-, second-, and third-order derivatives, were generated following 
the specification in \cite{birgin2018fortran}\footnote{The original implementations were in Fortran; 
our tests use a Matlab version obtained via a Fortran--Matlab converter \cite{cartis2024Efficient}.}.  We compare the performance of the following algorithms:  
\begin{itemize} \setlength{\itemindent}{0pt}
    \item \texttt{ARC}: Adaptive Regularization with Cubics \cite{cartis2011adaptive}.  
    \item \texttt{AR3 + ARC}: Third-order method with \texttt{ARC} as the inner solver \cite{cartis2020concise}.  
    \item \texttt{AR3 + DTM}: AR3 with DTM subproblem solver \Cref{TDTM} with diagonal rule in \Cref{algo: pdtm}.    
\end{itemize}
We first examine whether the third order derivative $T_0: = \nabla^3 f(x_0)$ in this test set exhibit diagonal or low-rank structure; 
the results are presented in \Cref{fig:mgh_distributions}. 
In our analysis of the MGH problems, we observe a broad range of tensor structures. 
Approximately $20\%$ of the problems have off-diagonal ratios below $7\%$, where
$\text{offdiag\_ratio}=\frac{\|T_{\mathrm{off}}\|_{F}}{\|T\|_{F}}$ where $T_{\mathrm{off}}$ denotes off-diagonal entries.  Several other problems, however, are dense. About $20\%$ of the problems have an off-diagonal ratio below $7\%$. The CP rank measures how many outer-product components are required to reconstruct a tensor.
The approximate CP rank is computed as the smallest $r$ such that $\|T_r - T\|_F / \|T\|_F < 10^{-2}$, where $T_r$ denotes a rank-$r$ approximation, and is then normalized by $r$ over number of tensor entries.  For the initial third-order tensors, we find that roughly $10\%$ of the problems are rank~1 and around $30\%$ are rank~2. 
Note that the size of the MGH test set is relatively small, with dimensions ranging from $n=2$ to $n=11$. Whether real-world large-scale problems exhibit similar diagonal or low-rank structure remains an open question for future investigation.

\begin{figure}[t]
    \centering
    \includegraphics[width=8cm]{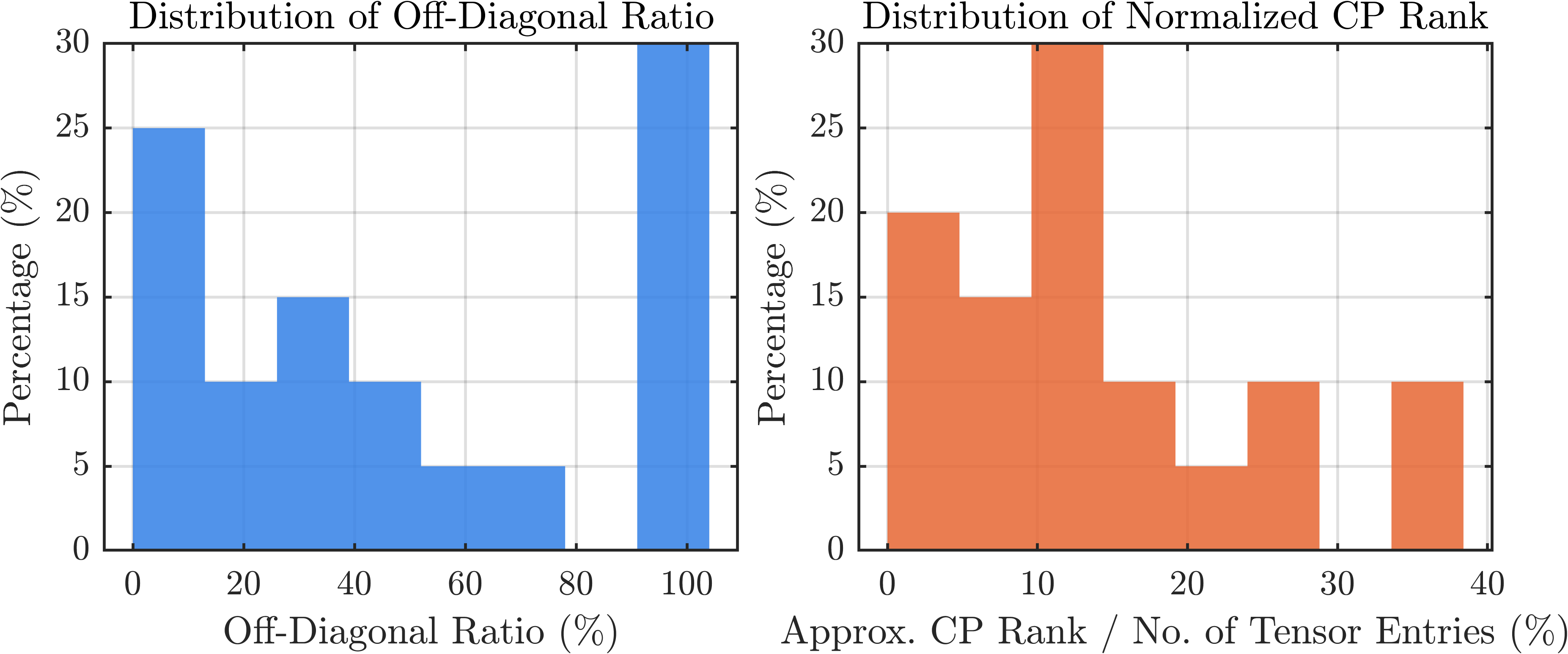}
    \caption{\small
        \textbf{Left:} Distribution of the off-diagonal ratio of the initial third-order tensors $T: = \nabla^3 f(x_0)$, expressed as a percentage of the total Frobenius norm $\|T\|_F$. 
        \textbf{Right:} Distribution of the approximate CP rank normalized by the number of tensor entries. 
    }
    \label{fig:mgh_distributions}
\end{figure}

The detailed numerical results for the problems in the MGH test set and comparisons is given \Cref{tab:full_comparison}.  
The outer-loop termination tolerance was set to $\|\nabla f\| < 10^{-3}$, and the maximum iteration count at 3000.
If an inner iteration was used, its tolerance was fixed at $10^{-4}$ and the maximum inner iteration count at 3000.  In performance tables, \texttt{eval.} refers to derivative evaluations (successful iterations only) and \texttt{iter.} denotes the total number of function evaluations (successful and unsuccessful). 
Performance profiles were used to compare the algorithms across the test set, following the framework of 
\cite{dolan2002benchmarking, cartis2024Efficient, gould2016note}.  
We adopt the notation of \cite{birgin2017use}:  
$\Gamma_{\hat{m}}(\tau_m)$ is the proportion of problems for which method $\hat{m}$ attained an 
$\epsilon_g$-approximate solution within a computational cost no more than $\tau_m$ times the cost of the most efficient method, where the cost is measured in iteration or evaluation counts\footnote{Some comments regarding CPU time are given at the end of the section.}.  Formally, if $t_{\hat{m}\hat{j}}(\epsilon_g)$ denotes the cost incurred by method $\hat{m}$ on problem $\hat{j}$ 
to reach an $\epsilon_g$-approximate solution, we define
$
t^{\hat{j}}_{\min}(\epsilon_g) := 
\min_{\hat{m}} t_{\hat{m}\hat{j}}(\epsilon_g).
$
where ${\hat{m} \in \{\texttt{ARC}, \texttt{AR3+ARC},  \texttt{AR3+CQR}, \texttt{AR3+DTM}\}}$.
If a method fails to satisfy the accuracy requirement within the iteration limit, we set $t_{\hat{m}\hat{j}}(\epsilon_g) = +\infty$.  

Across the test set, all methods consistently converged to the same minimizer whenever convergence occurred, and each computed solution satisfied the second-order necessary conditions. 
As shown in \Cref{fig performance plot long dtm} and \Cref{table:results}, third-order methods generally require no more outer function or derivative evaluations than the second-order method, and in some cases require even fewer. Among all methods, \texttt{ARC} remains the fastest in terms of overall computation time.
Between the third-order methods, \texttt{AR3 with DTM} requires the least number of function and derivative evaluations on average\footnote{The average excludes Problem 10.}, \texttt{AR3+DTM} requires 421 outer iterations and 410 derivative evaluations—fewer than \texttt{ARC} (507 and 505) and comparable to \texttt{AR3+ARC} (426 and 414). 
For each subproblem, the DTM solver is efficient: it typically converges in 1–3 inner iterations/evaluations, whereas \texttt{ARC} generally requires 7–9 inner iterations/evaluations per subproblem.
This advantage arises because each subproblem in DTM uses a diagonal tensor model, which more accurately captures the dominant diagonal structure of the AR3 subproblem tensor.
This effect is especially visible in problems that are approximately diagonal—for example, Problems 2, 14, and 17—where \texttt{AR3+DTM} demonstrates  stronger performance. 
Problem 10 and  Problem 20 are not solved by any of the three methods within the maximum limit of 3000 iterations; this behaviour of Problem 10 is also reported in \cite{birgin2017use, cartis2024Efficient}.
In Problem~20, when the 3000-iteration limit is reached, all solvers satisfy $\|\nabla f(x_k)\| \approx 10^{-2}$ and are clearly converging but have not yet reached high accuracy. Increasing the maximum iteration limit to 6000 allows third-order methods to converge fully, but \texttt{ARC} still fails to satisfy the gradient tolerance.

\begin{figure}[ht!]
  \begin{center}
    \includegraphics[width=12cm]{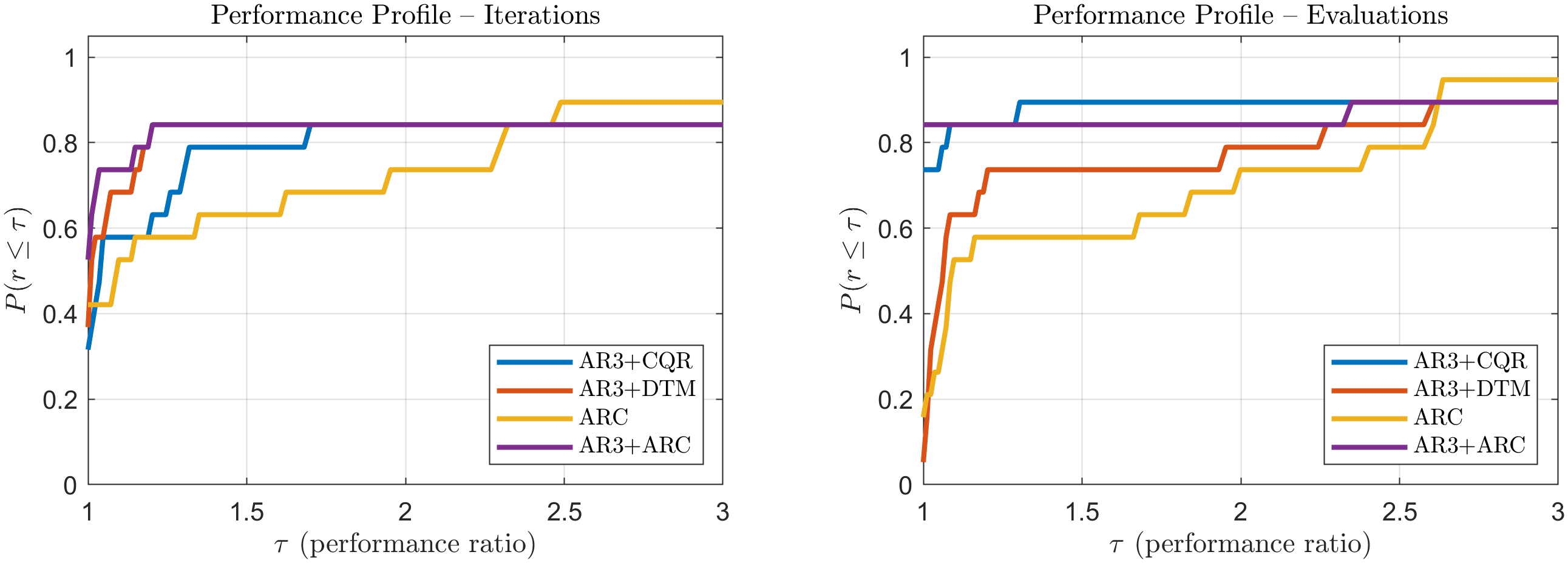}
    \caption{\small Performance profile plots for the three methods (ARC, AR3+ARC, AR3+CQR, and AR3+DTM). }
    \label{fig performance plot long dtm}
  \end{center}
\end{figure}

\subsection{A Low Rank Example: Packing Problem}
\label{sec: sphere}

We consider the problem of packing $N$ identical small spheres with radius $\tilde{r} > 0$ within a large sphere of radius $\tilde{R}  \ge \tilde{r}$. By packing, we mean
that the small spheres must be placed within the large sphere without overlapping (\cite{birgin2008minimizing, birgin2005optimizing, birgin2010new}). Following the formulation presented by Birgin et al., the problem can be modeled as follows
\begin{eqnarray}
\label{pack problem}
\min_{x \in \mathbb{R}^{n}} f_{pp}(x) = \sum_{i=1}^N \sum_{j=i+1}^N 
    \underbrace{\max\big\{ 0,\ (2\tilde{r})^2 - \| c_i - c_j \|^2 \big\}^4}_{\text{overlapping between the small spheres}}
+ \sum_{i=1}^N 
    \underbrace{\max\big\{ 0,\ \|c_i\|_2^2 - (\tilde{R} - \tilde{r})^2 \big\}^4}_{\text{fitting within the large sphere}} .
\end{eqnarray}
where $n=3N$, $x = (c_1^T, \ldots, c_N^T)^T$ represents $N$ the small balls are the centers $c_1,\ldots,c_N \in \mathbb{R}^3$. 
A solution of the packing problem corresponds to a solution of problem~\eqref{pack problem} at which $f_{pp}(x)$ is approximately minimized to zero. 
Following \cite[§13.1]{birgin2014practical}, the packing formulation considered here is central to the approach implemented in Packmol~\cite{martinez2009packmol}, widely used for initializing molecular dynamics systems. The radii $\tilde{r}$ and $\tilde{R}$ are related through $\tilde{R} \approx \tilde{r} \sqrt[3]{N}/0.3$, meaning that the small spheres occupy approximately $30\%$ of the volume of the large sphere. 

We initialize the centres of the small spheres near the origin using \(x_0 = 0.1\,\mathcal{N}(0,I_n)\) and examine how the CP rank and sparsity of the third-order tensor \(\nabla^{3} f(x)\) evolve throughout the optimization process.  As illustrated in \Cref{fig sphere rank}, the tensors arising from such random initializations are not low rank.  However, as the spheres gradually redistribute within the enclosing region and the iterates are driven toward a non-overlapping arrangement that fills the volume with radii $\tilde{R}$, the tensors computed at the final iterates begin to exhibit lower CP rank. 
An observation is that, for increasing dimension \(n\), 
the estimated CP rank grows approximately linearly with \(n\), rather than with the full \(n^{3}\) degrees of freedom in the tensor.  We also quantify sparsity by the proportion of nonzero tensor entries, and the results indicate that the tensor becomes sparser as the problem size increases. The sparsity structure becomes more pronounced as the iterates approach the minimizer.  
%Although the detailed evolution of rank and sparsity across iterations can be intricate, 
These preliminary findings suggest that third-order tensors arising in realistic optimization models may possess exploitable structure, such as approximate low rank or sparsity, even when such structure is absent at initialization. We compare the performance of \texttt{ARC} and \texttt{AR3+DTM} on this class of packing problems, with the results summarized in \Cref{fig sphere iterate}.  While \texttt{ARC} remains the fastest method in terms of computation time, \texttt{AR3+DTM} requires substantially fewer outer iterations and derivative evaluations to 
reach a minimizer of \(f_{pp}\). These experiments are preliminary numerical illustrations, highlighting the potential benefits of exploiting tensor structure within higher-order models. 
Note that we employed the standard \texttt{CP\_ALS} algorithm~\cite{bader2008efficient} to obtain 
an estimate of the best CP rank of a tensor via an alternating least-squares procedure. 
More advanced low-rank update strategies tailored to third-order methods remain an active area of research. Recent developments, such as the quasi-Newton–type updates proposed in \cite{welzel2023generalizing} for approximating third derivatives from Hessian evaluations, may offer significant in computational efficiency. \Cref{fig:sphere_packing_performance_and_tensor_structure} includes experiments up to \(n = 100\), corresponding to tensors with on the order of \(10^6\) entries. We note, however, that these experiments remain relatively small-scale and illustrative. A scalable, low-rank implementation of DTM is beyond the scope of the present work and constitutes a direction for our future research.

\begin{figure}[t]
    \centering
     \begin{subfigure}{0.48\textwidth}
        \centering
        \includegraphics[width=\linewidth]{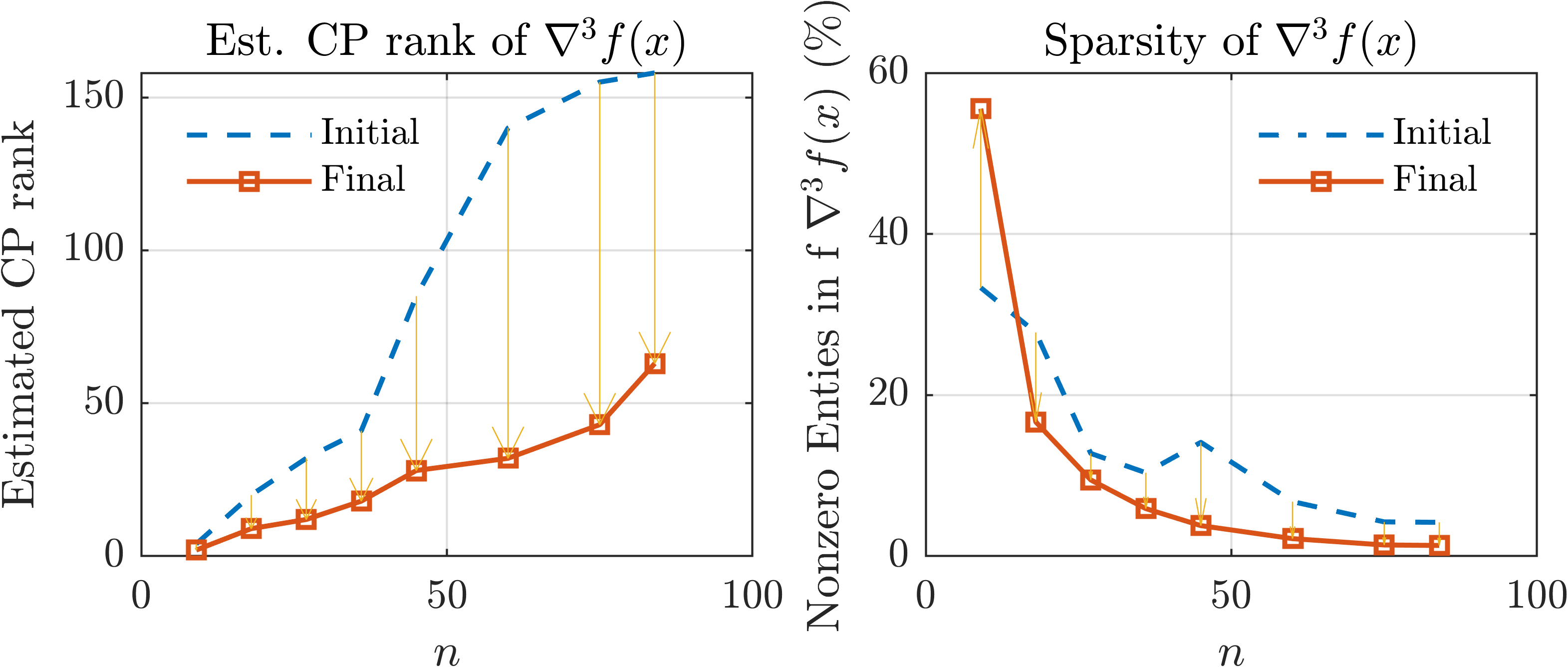}
        \caption{\small Approximate CP rank and sparsity of $\nabla^3 f(x)$ at initialization and after convergence.}
        \label{fig sphere rank}
    \end{subfigure}
    \hfill
   \begin{subfigure}{0.48\textwidth}
        \centering
        \includegraphics[width=\linewidth]{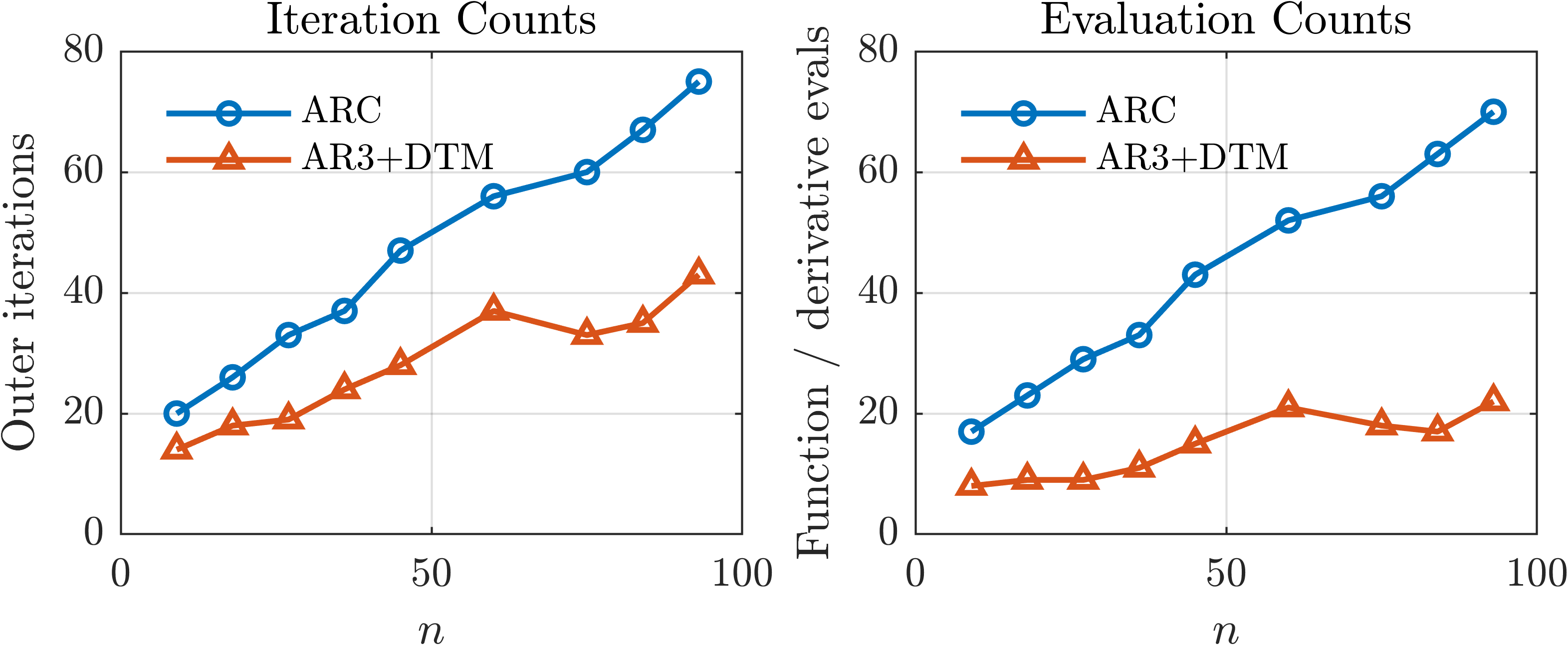}
        \caption{\small Iteration and evaluation counts for ARC and AR3+DTM.}
        \label{fig sphere iterate}
    \end{subfigure}
    \vspace{0.3em}
    \caption{
Performance of \texttt{ARC} and \texttt{AR3+DTM} on the sphere packing problem (right), 
and the structure of the third-order tensor $\nabla^{3} f(x)$ before and after optimization (left). 
In all experiments, the small-sphere radius is fixed at $r = 1$, and the $R = r\,(N/0.3)^{1/3}$. The problem dimension is $n = 3N$, initialized at $x_{0} = 0.1\,\mathcal{N}(0,I_{n})$. Each method is terminated once $\|f_{pp}(x)\| < 10^{-4}$. The approximate CP rank is computed as the smallest $r$ such that $\|T_r - T\|_F / \|T\|_F < 5\%$. }

    \label{fig:sphere_packing_performance_and_tensor_structure}
\end{figure}
}

\section{Conclusions and Future Work}
\label{sec: conclusion}

In this paper, we addressed the research gap in characterizing the global optimality conditions for the AR$3$ subproblem.  {The global optimality characterization is derived for the general form of \(m_3\), without assuming any special structure on the tensor or the Hessian. These conditions can be simplified in a straightforward manner and applied to quadratic polynomials with structured tensor terms or alternative regularisation norms. In particular, for special cases such as zero, low-rank, or diagonal third-order tensor terms, the optimality conditions admit further simplifications.} 
A highlight of this paper is identifying the conditions under which the necessary and sufficient conditions coincide.  
Specifically, under suitable assumptions for the Taylor polynomial, we derived a lower bound for the regularization parameter \( \sigma_k \) such that the necessary and sufficient criteria coincide, fully characterizing global optimality. Furthermore, we established a connection between this bound on the regularization parameter and polynomial sum-of-squares (SoS) convexification for locally convex subproblems, as well as the adaptive regularization framework for locally nonconvex subproblems. 
{We presented an algorithm for minimizing the quartic regularized polynomial using these optimality conditions, relying on Cholesky factorizations and Newton’s method for solving the associated secular equation. A full convergence analysis of this procedure is left for future work. In this algorithm, the first-order condition together with the global sufficient optimality condition yields a positive definite linear system. Although the method converges to a critical point of the quartic model, the global optimality conditions remain fully checkable a posteriori, allowing us to certify global minimality when they are satisfied.}

{By further exploiting only the diagonal entries of the third-order tensor, we introduce the Diagonal Tensor Method (DTM) framework, which iteratively solves AR3 subproblems.} 
DTM algorithm was specifically designed for the efficient minimization of nonconvex quartically-regularized cubic polynomials, such as the AR$3$ subproblem. The DTM framework minimizes the AR$3$ subproblem by iteratively solving a sequence of local quadratic models that incorporate both a diagonal cubic term and a quartic regularization term. Our theoretical analysis established the global convergence of the DTM algorithmic variant, proving it guarantees a first-order critical point for the AR$3$ subproblem within at most $\mathcal{O}(\epsilon^{-3/2})$ function and derivative evaluations, where $\epsilon$ represents the prescribed first-order optimality tolerance.

For subproblems with diagonal or low-rank third-order tensor information, the DTM model was able to exactly recover the subproblem, and therefore converges in one iteration after the change of basis. Preliminary numerical results further indicated promising performance by DTM compared to state-of-the-art methods like the ARC method. Notably, for subproblems with diagonal or low-rank third-order tensor information, DTM often required fewer iterations, evaluations, or computational time, underscoring its potential for practical applications.

Though the DTM model can exactly recover the AR3 subproblem in the low-rank third-order tensor case, we note that the cost of computing a CP decomposition is dominated by the standard alternating least-squares (ALS) procedure, whose cost is cubic in $n$ per ALS iteration for a dense third-order tensor in $\R^{n\times n\times n}$~\cite{kolda2009tensor}
% . In contrast, ARC can exploit Hessian sparsity to achieve subcubic complexity~\cite{cartis2011adaptive1}
, which partially explains the computational gap observed in our experiments.  
However, when the CP rank satisfies $r \ll n$, the DTM subproblem simplifies significantly, having a reduced the per-iteration cost of $O(rn^2)$. 
% More broadly, DTM and AR3 leverage third-order information and are particularly effective when the Lipschitz constant of the third derivative is small relative to that of the gradient~\cite{nesterov2006cubic}. In this regime, the cubic model provides a more accurate local approximation than quadratic models, leading to larger accepted steps and fewer outer iterations.
Therefore, closing the computational gap with ARC  depends critically on the presence of exploitable tensor structure. In particular, DTM can achieve competitive per-iteration cost when $T$ is sparse or admits a low-rank decomposition. Characterizing such problem classes and developing scalable implementations remain important directions for future work.
 Building on this perspective, a natural next step is to extend these methods to more scalable and practical tensor settings. This may involve integrating low-rank tensor approximations, homogeneous tensor formulations~\cite{chen2025tensor}, randomized techniques, and scalable iterative schemes, such as Krylov methods, eigenvalue-based formulations, or subspace optimization. Such enhancements would further improve the scalability of our solvers and broaden the applicability of high-order optimization methods to real-world problems. By proposing efficient algorithms for minimizing the AR3 model, this work takes a step toward making high-order tensor methods for general objectives more practically viable.

\smallbreak 
\noindent
\textbf{Acknowledgments}: This work was supported by the Hong Kong Innovation and Technology Commission (InnoHK Project CIMDA). Coralia Cartis's research was supported by the EPSRC grant EP/Y028872/1, Mathematical Foundations of Intelligence: An “Erlangen Programme” for AI.

\normalsize{}
\appendix

\section{Definition of Norms}
\label{appendix def of norm}

\begin{definition} (Following \cite{cartis2020concise}) Let $ T \in \R^{n^3}$ be a third-order tensor. The tensor norm of $ T$ is
\begin{eqnarray}
\Lambda_3 := \max_{\|v_1\|=\|v_2\|=\|v_3\|=1} \big|T[v_1][v_2][v_3]\big| 
\label{standard norm}
\end{eqnarray}
where $v_1, v_2 , v_3 \in \R^n$. $ T[v_1][v_2][v_3]$ stands for applying the third-order tensor $T$ to the vectors $v_1, v_2, v_3$. 
\label{tensor norm}
\end{definition}
\noindent 
According to \cite[Thm 2.4]{cartis2015branching}, 
$
\Lambda_3 \le \big[ \sum_{1 \le i,j,k \le, n} |T_{[i,j,k]}|^2. \big]^{1/2}
$
where  $T_{i,j,k}$ denotes the $i,j,k$th entry of $T$. 

\begin{corollary}  

Assume  $W$ is a symmetric, positive-definite matrix. Then
 \begin{eqnarray}
{\Lambda_W :=} \max_{v,u \in \R^n} \frac{|T[u][v]^2|}{\|u\|_W\|v\|_W^2 } \le \Lambda_3 [ \lambda_{\min}(W)]^{-3/2}
\end{eqnarray}
%The proof can be found in \Cref{sec: norm of T}.
\label{lemma Tensor norm}
\end{corollary}

\begin{proof} Using \eqref{standard norm}, $ \lambda_{\min}(W) \|v\|^2 \le v^TW v = \|v\|_W^2$ for all $v \in \R^n$, we deduce that
\begin{eqnarray}
 \max_{u, v \in \R^n} \frac{T[u][v]^2}{\|u\|_W\|v\|^2_W} = \max_{u, v \in \R^n} \frac{T[u][v]^2}{\|u\|\|v\|^2} \frac{\|u\|\|v\|}{\|u\|_W\|v\|^2_W} \le  \Lambda_3 [ \lambda_{\min}(W)]^{-3/2}.
\label{norm temp}
\end{eqnarray}
\end{proof}

\begin{remark}\label{remark explicit LambdaW}
An explicit upper bound is given by
\begin{equation}\label{eq:LambdaW_explicit}
    \Lambda_W 
    \leq \bigg[\sum_{1\leq i,j,k\leq n}
               |T_{i,j,k}|^2\bigg]^{1/2}
         [\lambda_{\min}(W)]^{-3/2},
\end{equation}
obtained by chaining the bound 
$\Lambda_3 \leq [\sum_{i,j,k}|T_{i,j,k}|^2]^{1/2}$ 
from \cite{cartis2020concise} 
with the inequality $\Lambda_W \leq \Lambda_3
[\lambda_{\min}(W)]^{-3/2}$ from 
Corollary~\ref{lemma Tensor norm}. 
This bound requires only 
the tensor entries $T_{i,j,k}$ and the smallest 
eigenvalue of $W$, both of which are available 
in the problem data.
\end{remark}

\begin{lemma}
\label{example special T}
In some special cases of $T$ and $W$, $\Lambda_W$ has a simplified expression. 
\begin{itemize}
    %\item If $T$ is a rank $r$ tensor with $0\le r \le n$, through an appropriate change ofz basis, we can convert $T$ into a diagonal tensor with $r$ non-zero entries. 
    \item If $T$ is a diagonal tensor with diagonal entries $ \{t_1, t_2, \dotsc, t_n\}$. Then, 
    $$
    \Lambda_W =   \max \{|t_1|, |t_2|, \dotsc, |t_n|\} [ \lambda_{\min}(W)]^{-3/2} .
    $$
\item If  $T$ is a diagonal tensor  with diagonal entries $\{t_1, \dotsc, t_n\}$ and $W = \diag\{|t_j|^{2/3}\}_{1 \le j\le n}$, then  $\Lambda_W =1$. 
\end{itemize}
%More details can be found in \Cref{sec: low-rank} and \Cref{sec: norm of T}.
\label{diagoal tensor norm}
\end{lemma} 

\begin{proof} 
If $T \in \R^{n^3}$ is a tensor with diagonal entries $\{t_1, \dotsc, t_n\}$, we have
\begin{eqnarray}
\max_{u, v \in \R^n} \big|T[u][v]^2\big|=  \sum_{j=1}^n |t_j| |v_j|^2 |u_j| \le t_*  \bigg(\sum_{j=1}^n |v_j|^2\bigg) \|u\| \le t_*\|v\|^2 \|u\| 
    \label{t norm}
\end{eqnarray}
where $t_* =   \max \{|t_1|, |t_2|, \dotsc, |t_n|\}$ and the first inequality uses $|u_i| \le \|u\|$. Substituting \eqref{t norm} into \eqref{norm temp} gives the first result in \Cref{diagoal tensor norm}. For the second result, if $W = \diag\{|t_j|^{2/3}\}_{1 \le j\le n}$, we have
$$
\max_{u, v \in \R^n} \big|T[u][v]^2\big| \le \sum_{j=1}^n (|t_j|^{2/3} |v_j|^2  )(|t_j|^{1/3} |u_j| )\le  \sum_{j=1}^n (|t_j|^{2/3} |v_j|^2  ) \|u\|_W  =\|v\|_W^2\|u\|_W
$$
where the first inequality uses $ |t_j|^{1/3} |u_j|  \le  \big[\sum_{j=1}^n |t_j|^{2/3} |u_j|^2\big]^{1/2}  = \|u\|_W.$ 
\end{proof}

\section{\texorpdfstring{Proof for \Cref{technical lemma dtm}}{Alternative Proof}}
\label{appendix alternative lemma}

\begin{proof}
For all vectors $v \in \R^n$, we have
\begin{eqnarray*}
   && m_3(s+v)  -f_0 = g^T (s+v)  +\frac{1}{2} H [s+v]^2   +  \frac{1}{6} T  [s+v]^3+ \frac{\sigma}{4} \big\|s+v\big\|_W ^4,
    \\  && =  \underbrace{g^T (s+v) + \frac{1}{2}\B_d(s) [s+v]^2 {- \frac{1}{12} T  [s]^3 - \frac{\sigma}{4} \|s\|_W^4 +\frac{1}{4} T [s][v]^2}}_{\mathcal{F}_1} +
     \\&&  \underbrace{ -\frac{1}{2} \bigg[\frac{1}{2}T [s] + \sigma W \|s\|_W^2 \bigg] [s+v]^2  +  { \frac{1}{6}  T  [s+v]^3+ \frac{\sigma}{4} \|s+v\|_W ^4  } +  {\frac{1}{12} T  [s]^3 + \frac{\sigma}{4} \|s\|_W^4 -\frac{1}{4} T [s][v]^2}}_{\mathcal{F}_2}. 
\end{eqnarray*}
Note that
\begin{eqnarray*}
\mathcal{F}_1   &=& g^T s + \frac{1}{2} \bigg[\B_d(s)  {- \frac{1}{6} T  [s]^3 - \frac{\sigma}{2} \|s\|^2}\bigg][s]^2 + \bigg[ g  + \B_d(s)s \bigg]^T  v+ \frac{1}{2} \B_d(s)[v]^2 + {\frac{1}{4}T [s][v]^2}, \notag
\\  &  =& (m_3(s) -f_0) + \bigg[g  + \B_d(s)s\bigg]^T  v+ \frac{1}{2} \G_d(s)[v]^2. 
\end{eqnarray*}
Rearranging $\mathcal{F}_2$ gives 
$
\mathcal{F}_2 =  \frac{1}{6} T  [v]^3 +  \frac{\sigma}{4} \big( \|s+v\|_W^2 - \|s\|_W^2\big) ^2.
$
In conclusion, we have $m_3(s+v)  =  m_3(s) + \mathcal{F}_1 + \mathcal{F}_2$ which gives \eqref{universal difference}. 
\end{proof}

An alternative proof for \Cref{technical lemma dtm} is also provided. 

\begin{proof}
By $4$th order Taylor expansion of $m_3(s+v)$, we have
$$
m_3(s+v)  - m_3(s) = \bigg[m_3(s) + \nabla m_3(s)^T v + \frac{1}{2}\nabla^2 m_3(s)[v]^2 + \frac{1}{6}\nabla^3 m_3(s)[v]^3 + \frac{1}{24}\nabla^4 m_3(s)[v]^4  \bigg] - m_3(s). 
$$
Let $\B_d(s)$ and  $\G_d(s)$ be defined as in \eqref{B and G}. Using $\nabla m_3(s) = g  + \B_d(s)s$,  $\nabla^2 m_3(s) = \G_d(s) + 2 \sigma(Ws_*)(Ws_*)^T $, $\frac{1}{6}\nabla^3 m_3(s)[v]^3 = \frac{1}{6}T[v]^3 +  \|v\|_W^2 s^TWv$ and $\frac{1}{24}\nabla^4 m_3(s)[v]^4  = \frac{\sigma}{4}\|v\|_W^4$, we deduce that
\begin{eqnarray*}
m_3(s+v)  - m_3(s) &=&    \bigg[g  + \B_d(s)s\bigg]^T  v + \frac{1}{2} \bigg[\G_d(s) + 2 \sigma sWs^T\bigg][v]^2 + \frac{1}{6}T[v]^3 +  \|v\|_W^2s^TWv  +  \frac{\sigma}{4}\|v\|_W^4.
\\ &=&  \bigg[g  + \B_d(s)s\bigg]^T  v + \frac{1}{2} \G_d(s) [v]^2 + \frac{1}{6}T[v]^3  +  \frac{\sigma}{4}\bigg[\|s\|_W^2 - \|s+v\|_W^2 \bigg]^2
\end{eqnarray*}
which is equivalent to \eqref{universal difference}. 
\end{proof}

\section{\texorpdfstring{Global Optimality Condition for \eqref{separable Polynomial} }{Global Optimality Condition of Separable Polynomial}}
\label{sec: CQR separable Global Opt}

In this section, we investigate the (global) optimality conditions of \eqref{separable Polynomial}. If $s_*$  is a local minimum, the first-order and second-order local optimality conditions of \eqref{separable Polynomial} yield
\begin{eqnarray}
 \nabla  m_{\textbf{SQR}}(s_*) = 0   \quad &\Rightarrow&  \quad     \hat{\B_d}(s_*) s_*: = \bigg(H + \frac{1}{2} \T [s_*] +   \Sigma [s_*]^2\bigg) s_* = -g_i,
 \label{local 1 seperable}
\\
 \nabla^2  m_{\textbf{SQR}}(s_*) \succeq 0  \quad &\Rightarrow&  \quad 
  H+  \T [s_*]  + 3  \Sigma [s_*]^2 \succeq 0 
   \label{local 2 seperable}
\end{eqnarray}
where $ \Sigma [s]^2 :=  \diag\big\{\sigma_j  s_j^2\big\}_{1\le j \le n} \in \R^{ n\times n}$.  
Before proving the global optimality conditions for  $m_{\textbf{SQR}}$, we need the following technical lemma.

\begin{lemma}
Let $\boldsymbol{\sigma} = [\sigma_1, \dotsc, \sigma_n]^T \in \R^n$, $t = [t_1 , \dotsc, t_n]^T \in \R^n$ and $
    \hat{\B_d}(s)  := H + \frac{1}{2}\T [s] +  \Sigma [s]^2. $
For any vector $s, v\in \R^n$, we have \eqref{universal difference seperable}.
\label{technical lemma seperable}
\end{lemma}

\begin{proof}
We assume wlog $f_0 $ = 0. For all vectors $v \in \R^n$, we have
\begin{eqnarray*}
   && m_{\textbf{SQR}}(s+v)  = g^T (s+v)  +\frac{1}{2} H [s+v]^2   +  \frac{1}{6} \T  [s+v]^3+ \frac{1}{4} \Sigma [s+v]^4,
    \\  && = \underbrace{g^T (s+v) + \frac{1}{2}\hat{\B_d}(s) [s+v]^2 - \frac{1}{12} \T  [s]^3 -  \frac{1}{4} \Sigma [s]^4 -\frac{1}{12} \T [s][v]^2}_{\hat{\mathcal{F}}_1} +
     \\&&  \underbrace{-\frac{1}{2} \bigg[\frac{1}{2}\T [s] + \Sigma [s]^2 \bigg] [s+v]^2  +  { \frac{1}{6}  \T  [s+v]^3+  \frac{1}{4} \Sigma [s+v]^4  } +  {\frac{1}{12} \T  [s]^3 +  \frac{1}{4}  \Sigma [s]^4 +\frac{1}{12} \T [s][v]^2}}_{\hat{\mathcal{F}}_2} . 
\end{eqnarray*}
\normalsize{}
Note that
\begin{eqnarray}
&&\hat{\mathcal{F}}_1   = g^T s + \frac{1}{2} \bigg[\hat{\B_d}(s)  {- \frac{1}{6} \T  [s]^3 - \frac{1}{2} \Sigma[s]^2}\bigg][s]^2 + \bigg[ g  + \hat{\B_d}(s)s \bigg]^T  v+ \frac{1}{2} \hat{\B_d}(s)[v]^2 - {\frac{1}{12}\T [s][v]^2}  , \notag
\\  && \quad =  m_{\textbf{SQR}}(s)+ \bigg[g  + \hat{\B_d}(s)s\bigg]^T  v+ \frac{1}{2} \hat{\B_d}(s)[v]^2 - {\frac{1}{12}\T [s][v]^2}. 
\label{F1 seperable}
\end{eqnarray}
Rearranging $\hat{\mathcal{F}}_2$ gives 
\begin{eqnarray}
\hat{\mathcal{F}}_2 &=& \frac{1}{3} \T [s][v]^2+ \frac{1}{6} \T  [v]^3 +   \frac{1}{4} \bigg[ -2 \Sigma[s]^2[s+v]^2 + \Sigma[s]^4 + \Sigma[v+s]^4\bigg].
\label{F2 seperable}
\end{eqnarray}
We add $\diag\big\{\frac{t_j^2}{36 \sigma_j}\big\}[v]^2$ to \eqref{F2 seperable} and give the SoS expression, 
\begin{eqnarray*}
\\ \hat{\mathcal{F}}_2 + \diag\big\{\frac{t_j^2}{36 \sigma_j}\big\}[v]^2 &=& \sum_{j=1}^n \bigg[\frac{1}{3} {s}_j v^2 +\frac{1}{6} v_j^3 +   \frac{\sigma_j}{4} \big[- 2( {s}_j+v)^2 {s}_j^2 +  {s}_j^4 + (v+ {s}_j)^4 \big]
 + \frac{t_j^2 v_j^2}{36 \sigma_j}\bigg]
 \\ &=& \frac{1}{4}  \sum_{j=1}^n  \bigg[\sigma_j v_j^2(v_j +2{s}_j+\frac{t_j}{3 \sigma_j})^2 \bigg] = \frac{1}{4}  \bigg\| \boldsymbol{\sigma}^{\frac{1}{2}} \cdot v \cdot (v + 2s +\frac{t \cdot \boldsymbol{\sigma} ^{-1}}{3}) \bigg\|^2.
\end{eqnarray*}
Combining \eqref{F1 seperable} and \eqref{F2 seperable}, we have,
$$m_{\textbf{SQR}}(s+v) = m_{\textbf{SQR}}(s)+ \big[g  + \hat{\B_d}(s)s\big]^T  v+ \frac{1}{2} \hat{\B_d}(s)[v]^2  - {\frac{1}{12}\T [s][v]^2}  -  \diag\big\{\frac{t_j^2}{36 \sigma_j}\big\}[v]^2 + \frac{1}{4} 
 \bigg\|\boldsymbol{\sigma}^{\frac{1}{2}} \cdot v \cdot (v + 2s +\frac{t \cdot \boldsymbol{\sigma} ^{-1}}{3}) \bigg\|^2. 
$$
Rearrange gives \eqref{universal difference seperable}. 
\end{proof}

\begin{theorem} (Global Sufficient Condition)
Assume \eqref{sufficient seperable} is satisfied. Then, \( s_* \) is a global minimum of \( m_{\textbf{SQR}}(s) \) over \( \mathbb{R}^n \).
\end{theorem}

\begin{proof}
Substituting \eqref{sufficient seperable} into \eqref{universal difference seperable} gives, for all \( v \in \mathbb{R}^n \),
$
m_{\textbf{SQR}}(s_*+v)  - m_{\textbf{SQR}}(s_*) \ge 0. 
$
\end{proof}

Note that the sufficient global optimality condition \eqref{sufficient seperable} implies the second-order local optimality condition \eqref{local 2 seperable}. To see this,
 for all \( v \in \mathbb{R}^n \), denote $s_*= [{s_*}_1, \dotsc, {s_*}_n]^T$, then
\begin{eqnarray*}
 \bigg[ H+  \T [s_*]  + 3 \Sigma [s_*]^2 \bigg][v]^2 \underset{\eqref{sufficient seperable}}{\ge}   \diag\bigg[ 2\sigma_j {s_*}_j^2 + \frac{2}{3}t_j{s_*}_j  + \frac{t_j^2}{18 \sigma_j}\bigg][v]^2 = 2 \diag\bigg[  \sigma_j\big({s_*}_j  + \frac{t_j}{6 \sigma_j}\big)^2\bigg][v]^2\ge0
\end{eqnarray*}

\begin{theorem} (Global Necessary Condition)
Assume $s_*$ is a global minimum of $m_{\textbf{SQR}}(s)$ over $\R^n$. Then, \eqref{necessary seperable} is satisfied.
\end{theorem}

\begin{proof}
Since $s_*$ is the global minimum we have $ m_{\textbf{SQR}}(s_*+v)  - m_{\textbf{SQR}}(s_*) \ge 0$ for all $v \in \R^n$.  $s_*$ is also the local minimum, therefore $\hat{\B_d}(s_*)s_* =-g$. If $2{s_*}_j +\frac{t_j}{3\sigma_j} \neq 0$, let $ \tilde{v} = \tilde{k} e_j := -( 2{s_*}_j +\frac{t_j}{3\sigma_j})e_j$.  Then, $\tilde{v}  \cdot (\tilde{v}  + 2s_* +\frac{t \cdot \boldsymbol{\sigma} ^{-1}}{3}) =0$. Substituting $\hat{\B_d}(s_*)s_* =-g$ and $ \tilde{v} = -\tilde{k} e_j$ into \eqref{universal difference seperable}  in \Cref{technical lemma seperable}, we have
\begin{eqnarray}
 0\le  m_{\textbf{SQR}}(s_*+\tilde{v} )  - m_{\textbf{SQR}}(s_*) =\frac{1}{2} \hat{\G_d}(s_*)[e_j]^2  \tilde{k}^2
 \label{temp step 3}, \qquad \forall j = 1, \dotsc, n
\end{eqnarray}
where\footnote{Note that  $\hat{\G_d}(s_*)[e_j]^2 \ge 0$ for all $1\le j \le n$ does not necessarily mean that $\hat{\G_d}(s_*) \succeq 0$, a counter example is $[0, 1; 1,0]$.} $ \hat{\G_d}$ is defined in \eqref{hat B and GS}. 
If $2{s_*}_j +\frac{t_j}{3\sigma_j} = 0$, set $ \tilde{v}  = e_j$, the same results follows. 

Now, take any unit vector $v = [v_1, \dotsc, v_n]^T = \sum_{j=1}^nv_je_j$, then
\begin{eqnarray}
\hat{\G_d}(s_*)[v]^2 = && \sum_{j=1}^n v_j^2\underbrace{ \hat{\G_d}(s_*) [e_j]^2}_{\ge 0 \text{ by \eqref{temp step 3}}} +2 \underbrace{\sum_{\iota, j=1, \iota \neq j}^n  (v_\iota v_j) H[e_\iota][e_j]}_{:=\textbf{H}}  + \notag \\+ && \sum_{\iota, j=1, \iota \neq j}^n 2 (v_\iota v_j) \underbrace{\bigg[ \frac{1}{3}\T [s_*] +\sigma_j   \Sigma [s_*]^2-\diag \big\{\frac{t_j^2}{18 \sigma_j} \big\} \bigg][e_\iota][e_j]}_{=0}.
\label{temp 7}
\end{eqnarray}
Using $H[e_\iota][e_j] = |H_{\iota,j}|$, we deduce that
\begin{eqnarray}
 \textbf{H}  \ge - \max_{\|v\| =1} |v_\iota v_j| \sum_{\iota, j=1, \iota \neq j}^n |H_{\iota, j}| \ge -\|H_0\|_1 
 \label{temp 8}
\end{eqnarray}
where  $\|H_0\|_1 =  \sum_{\iota, j=1, \iota \neq j}^n |H_{\iota, j}| $. Substituting \eqref{temp 8} in \eqref{temp 7} gives the result.  
%We construct a vector $\hat{v} \in \R^{n(n-1)}$ with each entry as $v_\iota v_j$ for $\iota, j=1, \dotsc, n$ and $\iota \neq j$. Also, construct vector  $\hat{h} \in \R^{n(n-1)}$ with each entry $= H_{\iota, j}$ ordered in the same way as $\hat{v}$.  $$ = -\hat{v}^T \hat{h} \le   \|\hat{v}^T \hat{h}\|_1  \le \|\hat{v} \|_1\|\hat{h}\|_{\infty}  = \max_{\|v\|=1, \iota  \neq j} v_\iota v_j  \max_{\iota  \neq j}  H_{\iota, j} \le  \frac{nh_*}{2}.$$
\end{proof}
If the second order term $H$ is a diagonal matrix, then $ \|H_0\|_1 =  \sum_{\iota, j=1, \iota \neq j}^n |H_{\iota, j}|  =0$. Consequently, the global necessary condition \eqref{necessary seperable} and the global sufficient condition \eqref{sufficient seperable} coincide. Furthermore, if  $H$ is a diagonal Hessian, minimizing \eqref{separable Polynomial} simplifies to the minimization of $n$ univariate CQR polynomials, as illustrated in \cite[Sec 3.2]{zhu2023cubic}.

\section{Optimality of Subproblems in the AR\texorpdfstring{$3$}{3} algorithm}
\label{appendix Difference in Objective Function and Taylor Expansion)}

\begin{assumption} (Lipschitz Property)
Assume that the $p$th derivative ($p\ge 1$) of $f$ at $x$ is globally Lipschitz; that is, there exists a constant $L > 0$ such that, for all
$x, y \in \R^n$, $\|\nabla^p_x f(x) - \nabla^p_x f(y)\| \le (p-1)! L \|x-y\|_W.$ 
\label{assumption Liptz}
\end{assumption}

\begin{remark}
There are slight differences in the definition of the Lipschitz constant across the literature. For instance, the Lipschitz constant in \cite{Nesterov2021implementable} differs from that in \cite{cartis2020concise} by a constant factor of the order of $p$. The Lipschitz properties in \cite{cartis2020concise, Nesterov2021implementable} are defined with respect to the Euclidean norm, specifically, $\|\nabla^p_x f(x) - \nabla^p_x f(y)\| \le (p-1)! \tilde{L} \|x-y\|$. In contrast, \Cref{assumption Liptz} uses the $W$-norm. The two definitions are related by the inequalities $ \lambda_{\min}(W) \tilde{L} \|x-y\| \le L \|x-y\|_W \le \lambda_{\max}(W) \tilde{L} \|x-y\|$.
\end{remark}

\begin{lemma}  \cite[Lemma 2.1 and Appendix A.1 ]{cartis2020concise}
Let $f\in \mathcal{C}^{p, 1}(\R^n)$ for $p\ge 2$ which means that $f$ is $p$-times continuously differentiable, bounded below and the $p$th derivative of $f$ is globally Lipschitz continuous, $p \ge 1$. $T_p(x, s)$ is the Taylor approximation of $f (x + s)$ about $x$. Then, under \Cref{assumption Liptz}, for all $x, s \in \R^n$, 
\begin{eqnarray}
      f(x+s) \le T_p(x, s) + \frac{L}{p}\|s\|_W^{p+1} 
       , \qquad  \|\nabla_x f(x+s) - \nabla_s T_p(x, s) \| \le L \|s\|_W^p. 
\end{eqnarray}
\end{lemma}

\begin{lemma}(Adapted from \cite[Lemma 3.1]{cartis2020concise})
 Let $f \in \mathcal{C}^{p, 1}(\R^n)$ for $p=3$, meaning that $f$ is $p$-times continuously differentiable, bounded below, and the $p$th derivative of $f$ is globally Lipschitz continuous (\Cref{assumption Liptz}).  If $\big\|\nabla_x f(x+s_*)\big\|>\epsilon_g$ for a positive tolerance $\epsilon_g$, then,  a first order critical point $s_*$ satisfies
 $\|s_*\|_W \ge \big[\frac{\epsilon_g}{L+\sigma}\big]^{1/3}$.
      \label{lemma lower bound for s_*}
\end{lemma}
\begin{proof}
Since $s_*$ is  a first-order critical point, therefore  $\nabla m_3(x_k, s_*) = 0$. We deduce that
    \begin{eqnarray*}
\epsilon &<& \big\|\nabla_x f(x+s_*)\big\| = \big\|\nabla_x f(x+s_*) - \nabla m_3(x_k, s_*)  \big\| 
\\ &\le & \big\|\nabla_x f(x+s_*) - \nabla T_3(x, s_*) \big\| + \big\|\nabla m_3(x_k, s_*)  \ - \nabla T_3(x, s_*) \big\| 
\\ &\le &  L \|s_*\|_W^3 + \big\|\nabla m_3(x_k, s_*)  \ - \nabla T_3(x, s_*) \big\| = (L + \sigma) \|s_*\|_W^3. 
\end{eqnarray*}
where $\nabla$ denote derivative with respect to $s$. Note that the last inequality uses Lemma 2.1 and Appendix A.1  \cite{cartis2020concise}  with slight modification to $W$-norm. 
\end{proof}

\begin{theorem}
\label{thm: how large sigma}
Assume that \Cref{assumption Liptz} holds and that $\|g\| \geq \epsilon_g$ for some positive tolerance $\epsilon_g$. 
Let $\lambda := \min\{0, -\lambda_{\min}[H]\} \lambda_{\min}(W)^{-1}$, and if $\sigma$ satisfies 
\begin{eqnarray}
 \sigma \ge \max\bigg\{L,  4 (3 \lambda)^{3} \epsilon_g^{-2}, \sqrt{2} (7 \Lambda_W)^{3/2} \epsilon_g^{-1/2}\bigg\},
 \label{sigma bound for nonconvex}
\end{eqnarray}
then \eqref{necessary tight DTM} and \eqref{sufficient tight DTM} are equivalent. Moreover, if $s_*$ and $\sigma$ satisfy $ \B_d(s_*) s_* := -g$, \eqref{necessary tight DTM} and \eqref{sigma bound for nonconvex}, then $s_*$ is a global minimum. 
\end{theorem}

\begin{proof}
If $ \sigma$ satisfies \eqref{sigma bound for nonconvex}, then
\begin{eqnarray*}
&&\sigma (L + \sigma)^{-2/3}\ge \sigma (2 \sigma)^{-2/3} = 2^{-2/3} \sigma^{1/3} \ge 3 \lambda   \epsilon_g^{-2/3}, 
\\ 
&& \sigma (L + \sigma)^{-1/3} \ge \sigma (2 \sigma)^{-1/3} = 
 2^{-1/3} \sigma^{2/3} \ge 7 \Lambda_W   \epsilon_g^{-1/3}. 
\end{eqnarray*}
We deduce that
\begin{eqnarray*}
\sigma \ge   3 \max \bigg\{ \lambda  (L + \sigma)^{2/3} \epsilon_g^{-2/3},   \frac{7}{3}\Lambda_W  (L + \sigma)^{1/3} \epsilon_g^{-1/3} \bigg\} \ge   3 \max \bigg\{-\lambda \|s_*\|_W^{-2},   \frac{7}{3}\Lambda_W  \|s_*\|_W^{-1}  \bigg\}
\end{eqnarray*}
where the last inequality uses $\|s_*\|_W^{-1} \leq (L + \sigma)^{1/3} \epsilon_g^{-1/3}$ from \Cref{lemma lower bound for s_*}. 
Lastly, using \Cref{thm sigma large enough}, we obtain the desired result.
\end{proof}

\begin{remark}
If $\lambda_{\min}[H] < 0$, then $\lambda := -\lambda_{\min}[H] \lambda_{\min}(W)^{-1} > 0$. For sufficiently small $\epsilon_g > 0$, the dominant term in \eqref{sigma bound for nonconvex} is $\mathcal{O}(\epsilon_g^2)$, leading to
$\sigma \geq 4 (3 \lambda)^{3} \epsilon_g^{-2}.$ 
\label{remark noconvex iff}
\end{remark}

%Combining the result of partial convexification of a locally nonconvex $m_3$ with \Cref{thm: how large sigma}, we deduce that, if $\sigma \geq 4(-3 \lambda)^3 \epsilon_g^{-2} = \mathcal{O}(\epsilon_g^{-2})$,  then $m_3$ is convex except in a small region around the origin $\|s\| < \mathcal{O}(\epsilon_g)$. Assuming that \Cref{assumption Liptz} holds and that $\|g\| \geq \epsilon_g$, then the global minimum of $m_3$ is characterized by $\B_d(s_*) s_* = -g$ and \eqref{necessary tight DTM} (\Cref{thm: how large sigma} and \Cref{remark noconvex iff}). Additionally, by \Cref{lemma lower bound for s_*}, we have $\|s_*\| \geq \big[\frac{\epsilon_g}{L+\sigma}\big]^{1/3} = \mathcal{O}(\epsilon_g^{1/3})$, which ensures that the global minimum lies within the convex region of $m_3$.

\section{Convergence and Complexity and DTM}
\label{appendix details of proof}

\begin{theorem}
\textbf{(An upper bound on step size, adapted from \cite[Thm 3.1]{zhu2023cubic})}  Suppose that  \Cref{assumption1 DTM}  holds and \Cref{algo: DTM variant 1}  is employed. Then, the norm of the iterates (both successful and unsuccessful) are uniformly bounded above independently of $i$, such that 
$\|s^{(i)}\|_W<r_c$ for all $i \ge 0$, where 
$r_c$ is a constant that depends only on the coefficients of $m_3$. 
\label{thm upper bound on step size DTM}
\end{theorem} 

\begin{corollary}\textbf{(Upper bounding  the tensor term, adapted from \cite[Lemma 3.3]{cartis2023second})} Let $L_H: = \Lambda_W+6 \sigma r_c$, where $\Lambda_W$ defined in \Cref{tensor norm}. Then, for every iteration, $i$, we have
$
  \|T_i [s]^2\|< L_H \|s\|_W^2
$ and 
$
  \|T_i [s]^3\|< L_H \|s\|_W^3
$ for all $s \in \R^n$. 
Note that $L_H$ is an iteration-independent constant that depends only on the coefficients in $m_3$. 
\label{corollary upper bound for Ti DTM}
\end{corollary}

\begin{lemma} 
\label{lemma: lower bound of norms DTM}
\textbf{(A lower bound on the step, adapted from \cite[Lemma 3.1]{zhu2023cubic})} Suppose that  \Cref{assumption1 DTM}  holds and \Cref{algo: DTM variant 1}  is used. Assume that $\|\nabla m_3(s^{(i)} + s_d^{(i)})\| > \epsilon$ and $\|g_i\| > \epsilon$, then if \footnote{If $d_i = 0$, $\|s_d^{(i)}\|_W > (B  +L_H )^{-1/2}\epsilon^{1/2}$.} $d_i \ge 0$, 
\begin{eqnarray}
\|s_d^{(i)}\|_W >   \min\bigg\{(B  +L_H )^{-1/2}\epsilon^{1/2}, \frac{1}{2}{d_i}^{-1/3}  \epsilon^{1/3} \bigg\}.
\label{Lower bound on step size DTM}
\end{eqnarray}
\end{lemma}

\begin{proof}
The proof is the same as the proof for \cite[Lemma 3.1]{zhu2023cubic}, with a slight modification: we replace $ \big|\beta_i \|s_d\|^2 \big|\le B \|s_d\|^2 $ by $ \|\T_i [s_d]^2\| \le B \|s_d\|_W^2$ and change the regularization norm from $\|\cdot\|$ to $\|\cdot\|_W$. 
\end{proof}

\begin{comment}

\begin{proof}
For notational simplicity, in this proof, we use $s_d$ to represent $s_d^{(i)}$. 
Since $s_d$ is the global minimum of $M_d$, we have $\nabla M_d(s^{(i)}, s_d) = 0$. Thus, 
$$
\epsilon<\|\nabla m_3(s^{(i)} + s_d)\| = \| \nabla  M(s^{(i)}, s_d)\| = \|\nabla M_d(s^{(i)}, s_d) - \nabla  M(s^{(i)}, s_d)\|. 
$$
Using the expression of $\nabla M_d(s^{(i)}, s_d) = 0$ and  $\nabla M(s^{(i)}, s) = g_i + H_i s +\frac{1}{2}T_i[s]^2 + \sigma\|s\|^2 s$, we have
\begin{eqnarray}
\epsilon< \bigg\| \frac{1}{2}  \T_i [s_d]^2-  \frac{1}{2}  T_i [s_d]^2 +   4d_i   \|s_d\|^2 s_d \bigg\| \le \frac{B}{2}\|s_d\|^2 +    \frac{L_H}{2} \big\|s_d\big\|^2+ 4d_i\|s_d\|^3
\label{diff grad}
 \end{eqnarray}
 where the last inequality uses norm properties, $ \T_i  > -B$ from \Cref{assumption1 DTM} and $\big\|T_i\big\| \le L_H$ from \Cref{corollary upper bound for Ti DTM}. 
If $d_i \ge 0$, either one of the following must be true
\begin{eqnarray*}
\frac{ \epsilon}{2} < \bigg(\frac{B}{2}  +  \frac{L_H}{2}  \bigg) \|s_d\|^2 \Rightarrow   \|s_d\| > (B+ L_H)^{-1/2}\epsilon^{1/2},
\quad \text{or} \quad  \frac{ \epsilon}{2} <  4 d_i \|s_d\|^3   \Rightarrow\|s_d\| >\frac{1}{2}{d_i}^{-1/3}  \epsilon^{1/3}
\end{eqnarray*}
which gives \eqref{Lower bound on step size DTM}. 
\end{proof}
    
\end{comment}

\begin{theorem}
\label{thm cubic Upper Bound DTM}
\textbf{(A local upper bound for $m_3$, adapted from \cite[Thm 3.2]{zhu2023cubic})} Suppose that  \Cref{assumption1 DTM}  holds and \Cref{algo: DTM variant 1}  is used.
Then, if
    \begin{eqnarray}
        d_i > \frac{B + L_H}{6}\|s_d^{(i)}\|_W^{-1}>0,
        \label{dc DTM}
    \end{eqnarray}
     we have 
    $
    m_3(s^{(i)}+s_d^{(i)}) = M(s^{(i)},s_d^{(i)}) \le M_d(s^{(i)},s_d^{(i)}). 
    $ 
with the equality only happens at $s_d^{(i)} = 0$. 
\end{theorem}

\begin{proof}
For notational simplicity, in this proof, we use $s_d$ to represent $s_d^{(i)}$. By definition of $M_d$ and $M$, we have
\begin{eqnarray}
 M_d(s^{(i)}, s_d) -  M(s^{(i)}, s_d) =  \frac{1}{6} \bigg[ [\T_i-T_i] [s_d]^3 \bigg]+ d_i\|s_d\|_W^4. 
\label{difference in Mc and M DTM}
\end{eqnarray}
If $s_d = 0$, clearly $M(s^{(i)},0) = M_d(s^{(i)},0). $  Assume that $\|s_d\|_W > 0$, we choose $d_i > \frac{1}{6}(B + L_H)\|s_d\|_W^{-1}>0$. Using 
$\T[s_d]^3  > -B\|s_d\|_W^3$ from \Cref{assumption1 DTM} and \Cref{remark for assumption}, and $T[s_d]^3 \ge -L_H \|s_d\|_W^3$ from \Cref{corollary upper bound for Ti DTM}, we deduce from \eqref{difference in Mc and M DTM} that 
$
    M_d(s^{(i)}, s_d) -  M(s^{(i)}, s_d)> \frac{1}{6}\big(-B  - L_H \big)\|s_d\|_W^3 + d_i\|s_d\|_W^4 
    \underset{\eqref{dc DTM}}{>} 0.
$
\end{proof}

\begin{remark} If $  d_i$  satisfies \eqref{dc DTM}, then $ \sigma_d^{(i)}$ satisfies $
    \sigma_d^{(i)} \ge \frac{2}{3} \big(-\beta_i + \alpha\big)\|s_d^{(i)}\|_W^{-1} $ in \Cref{algo: DTM variant 1}. To see this, since $B > -\beta_i$ , we have $d_i  > \frac{1}{6}(B + L_H)\|s_d^{(i)}\|_W^{-1} >  \frac{1}{6} \big(-\beta_i + \alpha\big)\|s_d^{(i)}\|_W^{-1}  $
    and, consequently,  $\sigma_d^{(i)}  = \sigma + 4d_i\ge \frac{2}{3} \big(-\beta_i + \alpha\big)\|s_d^{(i)}\|_W^{-1}.$ 
\end{remark}

\begin{lemma}
\label{lemma upper bound on di DTM}
\textbf{(An upper bound for $d_i$, adapted from \cite[Lemma 3.4]{zhu2023cubic})}
Suppose that  \Cref{assumption1 DTM} holds and \Cref{algo: DTM variant 1}  is used. 
Then, for all $i \ge 0$
\begin{eqnarray}
 d_i \le d_{\max} :=  \gamma (B + L_H)^{3/2}\epsilon^{-1/2},
\label{upper bound on di} 
\end{eqnarray}
where $\gamma$ is a fixed constant from \Cref{algo: DTM variant 1} . 
\end{lemma}

\subsection{Proof of \texorpdfstring{\Cref{thm M3 value decrease DTM}}{Theorem M3 value decrease DTM}}
\label{appendix proof of complexity thm}

For notational simplicity, in this proof, we use $s_d$ to represent $s_d^{(i)}$ and $\sigma_d$ to represent $\sigma_d^{(i)}$. Since $i$ is a successful iteration, \eqref{ratio test} gives
$$
\mathcal{D}:=\eta^{-1}\left[m_3(s^{(i)}) - m_3(s^{(i+1)})\right] \ge M_d(s^{(i)}, 0) - M_d(s^{(i)}, s_d) 
= - g_i^T s_d - \frac{1}{2} H_i [s_d]^2 - \frac{1}{6} \T_i [s_d]^3 - \frac{\sigma_d}{4}\|s_d\|_W^4.
$$
$s_d$ is the global minimum of $M_d$  and satisfies $-g_i  = \B_d_i(s_d)s_d$. We have
\begin{eqnarray} 
M_d(s^{(i)}, 0) - M_d(s^{(i)}, s_d) =\frac{1}{2} H_i [s_d]^2 +  \frac{1}{3} \T_i [s_d]^3+ \frac{3\sigma_d}{4}\|s_d\|_W^4.
\label{decrease on Mc DTM}
\end{eqnarray}
Using \Cref{technical lemma dtm} with $-g_i  = \B_d_i(s_d)s_d$, we have
$$
 M_d(s^{(i)}, s_d+v) - M_d(s^{(i)}, s_d) =  \frac{1}{2} \bigg[ H_i +  \T_i[s] +  \sigma_d W\|s\|_W^2\bigg][v]^2 + \frac{1}{6} T  [v]^3 +  \frac{\sigma}{4} \bigg[ \|s+v\|_W^2 - \|s\|_W^2\bigg] ^2.
$$
\begin{comment}
Let $v = k s_d$, we have
\begin{eqnarray*}   
0 &\le& \frac{1}{2}  k^2 H[s_d]^2 +  \frac{1}{2}  k^2 \T[s_d]^3 +  
 \frac{1}{2}  k^2 \sigma_d \|s_d\|_W ^4 + \frac{1}{6}k^3 \T[s_d]^3 +  \frac{\sigma_d}{4} \bigg[ (1+k)^2-1 \bigg]^2 \|s_d\|_W ^4.
\\
0 &\le& \frac{1}{2}  H[s_d]^2  + \frac{1}{6}(k+3) \T[s_d]^3 +  \frac{\sigma_d}{4} \bigg[(k+2)^2 + 2\bigg]  \|s_d\|_W ^4.
\end{eqnarray*}
we have
\begin{eqnarray}
0 \le M_d(s^{(i)}, -s_d) - M_d(s^{(i)}, s_d) = \frac{1}{2} H_i [s_d]^2 +  \frac{1}{3} \T_i [s_d]^3 +  \frac{\sigma_d}{2} \|s_d\|_W^4.
 \label{decrease on Mc 2}
\end{eqnarray}
\end{comment}
Let $v =-2s_d$,  
\begin{eqnarray}
 \frac{1}{2} H_i [s_d]^2 +  \frac{1}{6} \T_i [s_d]^3 +  \frac{\sigma_d}{2} \|s_d\|_W \ge 0.    
  \label{decrease on Mc 2}
\end{eqnarray}
\begin{itemize}
    \item Assume that $\T_i \big[\frac{s_d}{\|s_d\|_W} \big]^3 \ge \alpha >0$, substituting \eqref{decrease on Mc 2} into \eqref{decrease on Mc DTM}, we have
$
\mathcal{D} \ge  \frac{1}{6} \T_i [s_d]^3+ \frac{\sigma_d}{4}\|s_d\|_W^4 \ge \frac{\alpha }{3}  \|s_d\|_W^3.
$
\item Assume that $\T_i \big[\frac{s_d}{\|s_d\|_W} \big]^3 \le \alpha$ and $\sigma_d \ge \frac{2}{3} \bigg(- \T_i \big[\frac{s_d}{\|s_d\|_W} \big]^3 + \alpha\bigg)\|s_d\|_W^{-1} > 0$, 
$$
\mathcal{D}\ge\frac{1}{6} \T_i [s_d]^3+ \frac{1}{6} \bigg(- \T_i \bigg[\frac{s_d}{\|s_d\|_W} \bigg]^3 + \alpha\bigg)\|s_d\|_W^3 \ge \frac{\alpha}{6}  \|s_d\|_W^3.
$$
\item Assume that $\T_i \big[\frac{s_d}{\|s_d\|_W} \big]^3  \in  [-B, -4\alpha] $ and $0 \le \sigma_d \le \frac{1}{6} \bigg(- \T_i \big[\frac{s_d}{\|s_d\|_W} \big]^3 + \alpha\bigg)\|s_d\|_W^{-1}$, according to the second-order local optimality condition, 
$
0 \le \nabla^2 M_d(s^{(i)}, s_d)[s_d]^2 = H_i[s_d]^2 + \T_i[s_d]^3 + 3 \sigma_d \|s_d\|_W^4. 
$
Substituting the local optimality condition into \eqref{decrease on Mc DTM}, we have
\begin{eqnarray*}
\mathcal{D} \ge  -\frac{1}{6} \T_i [s_d]^3 - \frac{3\sigma_d}{4}\|s_d\|_W^4 \ge  -\frac{1}{6} \T_i [s_d]^3 - \frac{1}{8} \bigg(- \T_i \bigg[\frac{s_d}{\|s_d\|_W} \bigg]^3 + \alpha\bigg)\|s_d\|_W^3 \ge - \frac{1}{24} \T_i [s_d]^3 -  \frac{\alpha}{8} \|s_d\|_W^3 \ge \frac{\alpha}{24} \|s_d\|_W^3
\end{eqnarray*}
where the last inequality uses  $ \T_i \big[\frac{s_d}{\|s_d\|_W} \big]^3 \le -4\alpha.$ 
\end{itemize}
In all cases, we have 
$
    m_3(s^{(i)})-  m_3(s^{(i+1)})  \ge   \frac{\alpha \eta}{24}\|s_d\|_W^3.
$
These three cases correspond to the requirements of $\beta_i$ and $\sigma_d$ in \Cref{algo: DTM variant 1}. 
Using the lower bound on $\|s_d\|_W$ in \Cref{lemma: lower bound of norms DTM} and the bound on $d_i$ in \Cref{lemma upper bound on di DTM}, 
\begin{eqnarray*}
    m_3(s^{(i)})-  m_3(s^{(i+1)}) \ge   \frac{\alpha \eta}{24}  \max\bigg\{(B  +L_H )^{-3/2}\epsilon^{3/2}, \frac{1}{8}{d_{\max}}^{-1}  \epsilon \bigg\} \ge \frac{\alpha \eta}{24}  \min\bigg\{1, \frac{1}{8\gamma}   \bigg\} (B + L_H)^{-3/2} \epsilon^{3/2}. 
\end{eqnarray*}
Using $\gamma>1$, gives the desired result in \eqref{general model decrease DTM}.

\begin{remark}
\label{remark cauchy analysis DTM}
Note that the $\sigma_d^{(i)}$ condition in  \Cref{algo: DTM variant 1} is meant to ensure that the decrease in values of $m_3$ is at least $\mathcal{O}(\epsilon^{3/2})$ as indicated in \eqref{general model decrease DTM}. If we only have the condition $\rho_i \ge \eta$, \Cref{TDTM} still results in a decrease in the value of $m_3$ at each successful iteration. We can verify the decrease in value of $m_3$ using a standard Cauchy analysis argument as illustrated in \cite[Remark 3.4]{zhu2023cubic}.
\end{remark}

\subsection{Proof of \texorpdfstring{\Cref{thm: dtm complexity}}{thm dtm complexity}}
\label{appendix dtm complexity}

\begin{proof}
The complexity proof follows the same approach as in \cite[Lemma 3.3, Lemma 3.4]{zhu2023cubic}. We begin by bounding the number of successful iterations, denoted as 
$
|\mathcal{S}_i|:= \kappa_s \frac{m_3(s^{(0)})-m_{\text{low}}}{\epsilon^{3/2}}
$ 
and then proceed to bound the total number of iterations, given by
$
i \le |\mathcal{S}_i|  + \frac{\log(d_{\max})}{\log \gamma}.
$ 
Substituting  $|\mathcal{S}_i|$ and applying the bound on 
$d_{\max}$ from \eqref{upper bound on di} yields the desired result with $\kappa_c:=  1  + \frac{3 \log(B + L_H)}{2\log\gamma}$. A detailed proof can be found in \cite[Thm 3.4]{zhu2023cubic}.
\end{proof}

\section{Supplement for Numerical Results}

\subsection{Size of Tensor and Global Optimality Characterization}
\label{appendix: size}

\begin{figure}[ht]
    \centering
\caption{\small Global optimality condition and size of $\|T\|$ or $\sigma$ using the Diagonal Tensor Test Set:  \\\textbf{First plot}: \texttt{a = 10, b = 20}, $\sigma = 100$, \texttt{c} $\in [5, 80]$. 
 \textbf{Second plot}: \texttt{a = 10, b = 20, c= 50}, $\sigma$ $\in [50, 300]$.}
    \label{fig: increase t or sigma}
    \includegraphics[width =7cm]{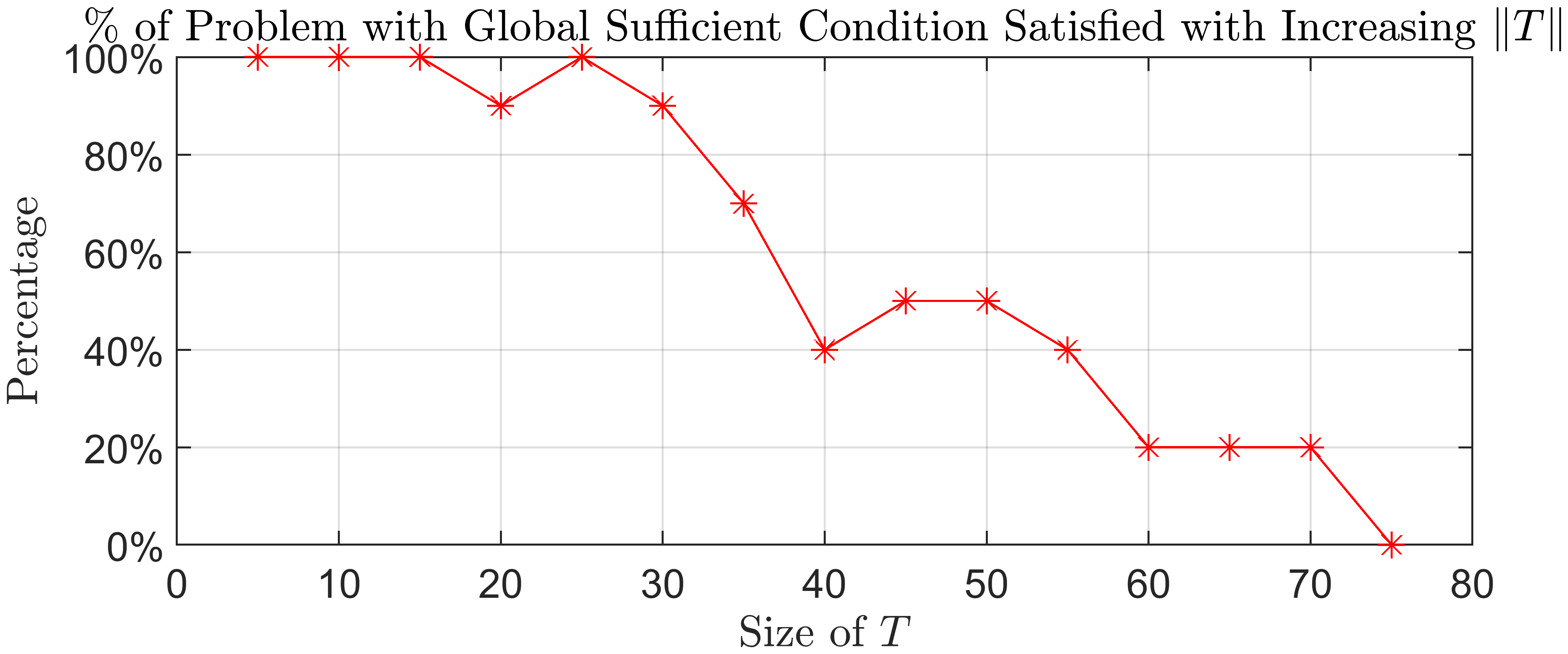}
    \includegraphics[width =7cm]{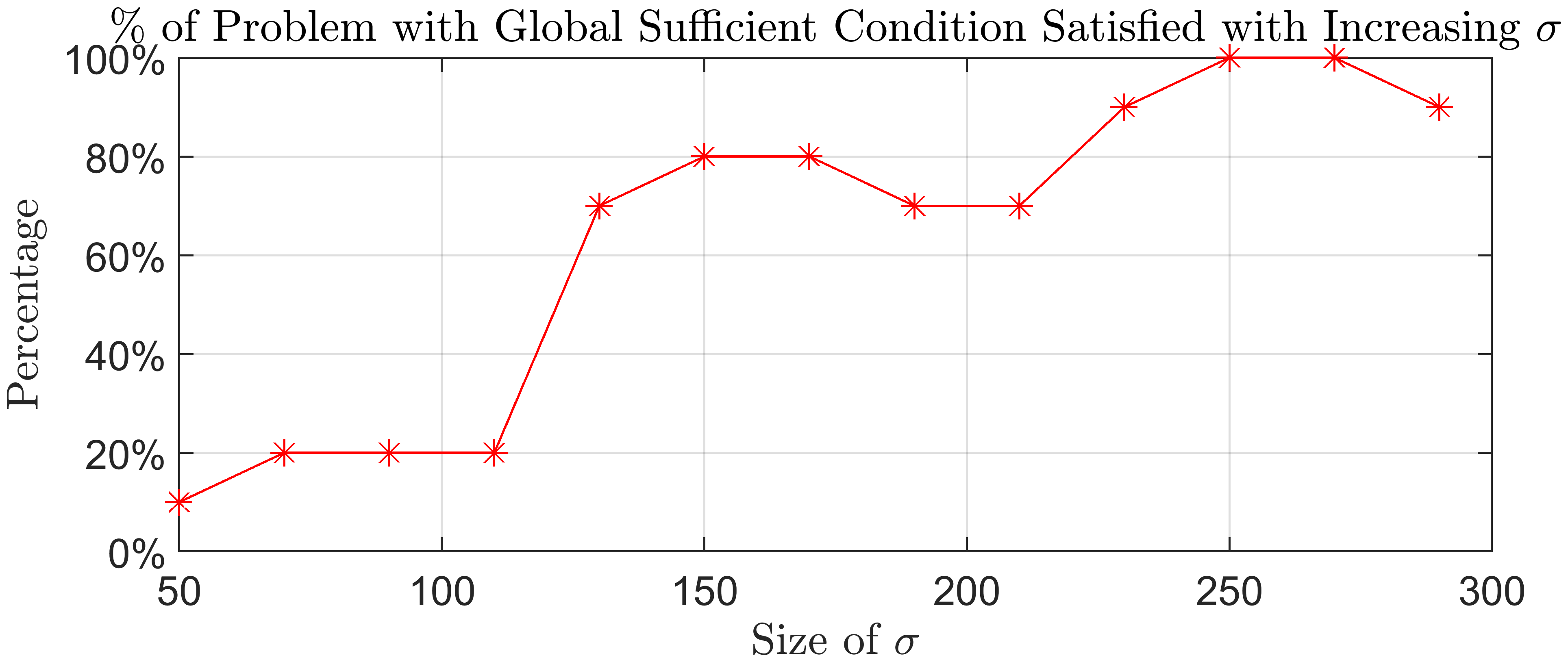}
\label{fig increase sigma}
\end{figure}

\subsection{\texorpdfstring{Computational Time Comparison for Applying \Cref{TDTM} Directly to AR3 Subproblem}{Computational Time Comparison for Applying Ago 2 Directly to AR3 Subproblem}}
\label{appendix: time}

\begin{figure}[ht]
    \centering
\caption{\small
Computational time comparison: {\Cref{TDTM} with $\T_i = T_i$, where $T_i$ is the full third-order derivative, 
versus with $\T_i$ chosen as the diagonal tensor using the Diagonal Rule in \Cref{algo: pdtm}. 
The test functions and experimental setup are identical to those in \Cref{table:results}.}}
\includegraphics[width =7cm]{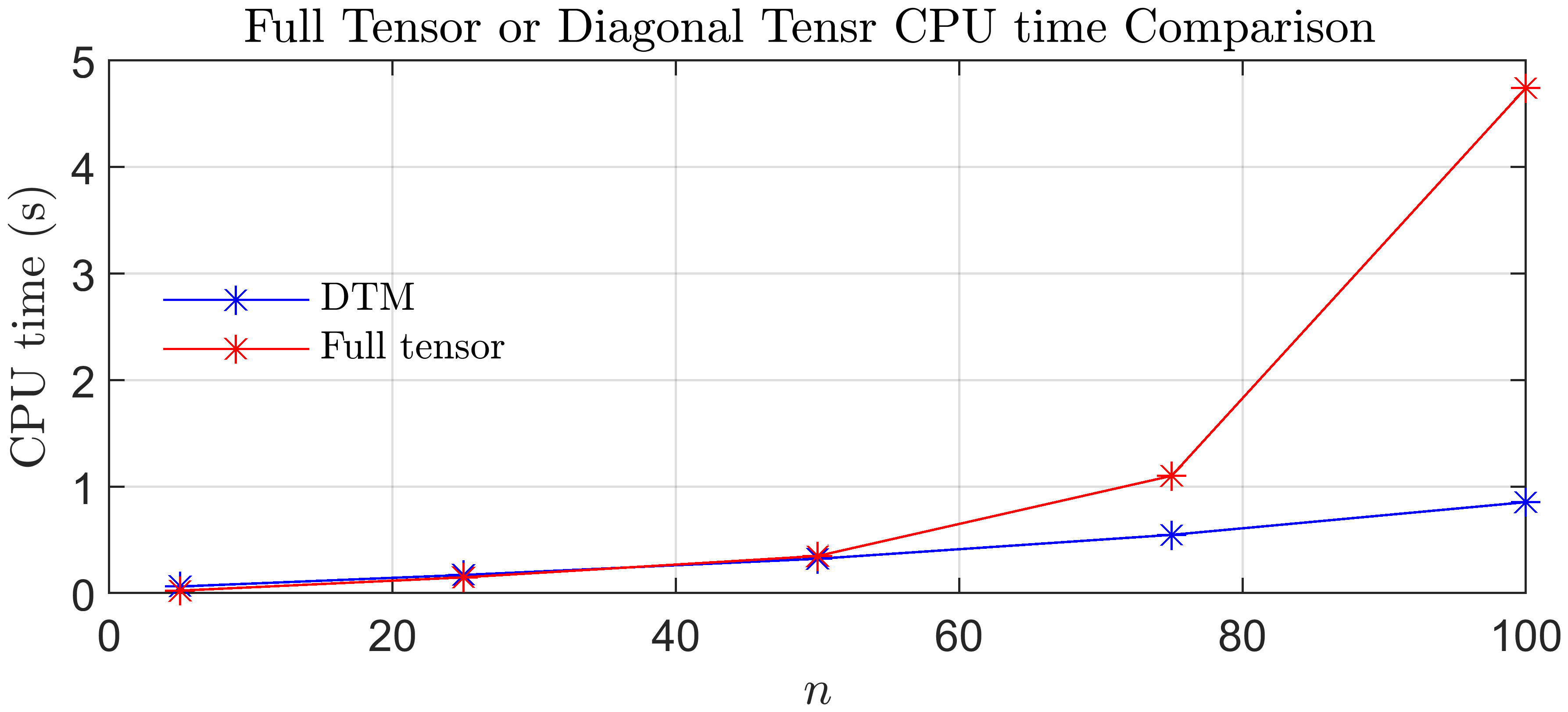}
\label{fig computational time}
\end{figure}

\subsection{Experiment with MGH Test Set}
\label{Appendix Experiment with MGH Test Set dtm}

In \Cref{tab:full_comparison},  $-$ indicates that the algorithm failed to make progress after 3000 iterations. 
${\dagger}$ indicates that the algorithm timed out:  $\|\nabla f\|$ was reduced below $10^{-2}$ but did not reach the prescribed outer-loop tolerance of $10^{-3}$. 

\begin{table}[ht]
\centering
\scriptsize
\setlength{\tabcolsep}{2pt}
\renewcommand{\arraystretch}{1.2}

\caption{\small Comparison across ARC, AR3+ARC, and AR3+DTM over the 20 CUTEst problems, including total iteration counts, function evaluations, and per-subproblem efficiency. {The average excludes Problem 10.}}
\begin{tabular}{c c c | c c | c c c c | c c c c}
\hline
\textbf{Prob.} & \textbf{CP Rank.} & \textbf{Off-Diag.} &
\multicolumn{2}{c|}{\textbf{ARC}} &
\multicolumn{4}{c|}{\textbf{AR3+ARC}} &
\multicolumn{4}{c}{\textbf{AR3+DTM}} \\
 & & &
\textbf{Iter} & \textbf{Eval} &
\textbf{Iter} & \textbf{Eval} & \textbf{Iter/Sub} & \textbf{Eval/Sub} &
\textbf{Iter} & \textbf{Eval} & \textbf{Iter/Sub} & \textbf{Eval/Sub} \\
\hline
 1 &  1 &  8.0\%   & 2322 & 2322 & 2354 & 2311 &  4.4  &  4.1  & 2356 & 2312 &  2.4 &  1.3 \\
 2 &  1 &  0.3\%   &  113 &  113 &   58 &   57 &  2.9  &  2.9  &   61 &   58 &  1.6 &  1.7 \\
 3 &  2 & 14.3\%   &   26 &   13 &   31 &   12 &  7.4  & 18.2  &   31 &   12 &  2.1 &  5.3 \\
 4 &  3 &100.0\%   & 1421 & 1421 &   32 &   31 & 10.9  &  6.0  &   32 &   31 & 11.2 &  2.1 \\
 5 &  3 & 32.3\%   &   66 &   66 &   58 &   57 &  2.6  &  2.7  &   58 &   57 &  1.4 &  1.4 \\
 6 &  2 &  7.9\%   &  104 &  104 &  119 &   99 & 15.0  & 16.3  &  119 &   99 & 16.0 &  3.2 \\
 7 &  3 & 51.3\%   &  782 &  781 &  717 &  716 &  2.2  &  2.2  &  717 &  716 &  1.2 &  1.2 \\
 8 &  3 & 49.8\%   &  745 &  745 &  726 &  725 &  2.0  &  2.0  &  726 &  725 &  1.0 &  1.0 \\
 9 &  5 & 52.9\%   &   31 &   31 &   23 &   12 &  5.3  & 10.3  &   27 &   13 &  1.3 &  2.6 \\
10$^\dagger$ &  5 &100.0\%   &  --  &  --  &  --  &  --  &  --   &  --   &  --  &  --  &  --   &  --   \\
11 &  3 &100.0\%   &   16 &   13 &   60 &   49 & 20.6  & 25.2  &   57 &   47 &  1.3 &  1.5 \\
12 &  6 & 30.3\%   &   42 &   42 &   26 &   25 &  2.8  &  2.9  &   26 &   25 &  1.4 &  1.5 \\
13 &  2 & 30.5\%   &  223 &  223 &  208 &  206 &  2.2  &  2.3  &  209 &  207 &  1.2 &  1.2 \\
14 &  2 &  3.2\%   &  116 &  116 &  119 &   63 & 13.5  & 21.9  &  123 &   64 &  4.5 &  7.8 \\
15 &  8 & 98.0\%   &   50 &   50 &   22 &   21 &  3.0  &  3.1  &   22 &   21 &  1.1 &  1.1 \\
16 &  3 & 16.7\%   &   91 &   91 &   92 &   85 &  3.5  &  3.5  &   92 &   85 &  4.8 &  2.2 \\
17 &  2 &  3.2\%   &  400 &  387 &  173 &  148 & 13.4  & 14.9  &  173 &  148 &  2.7 &  3.1 \\
18 & 20 &100.0\%   &   42 &   42 &   17 &   16 &  3.2  &  3.4  &   17 &   16 &  1.9 &  2.1 \\
19 &  4 & 94.0\%   &   38 &   38 &  251 &  238 & 28.8  & 30.4  &  159 &  146 &  1.1 &  1.2 \\
20$^\dagger$ & 20 & 67.9\%   & 3000 & 3001 & 3001 & 3001 &  2.0  &  2.0  & 3001 & 3001 &  1.0 &  1.0 \\
\hline
\textbf{Average} & & &
\textbf{507} & \textbf{505} &
\textbf{426} & \textbf{414} &
\textbf{7.7} & \textbf{9.2} &
\textbf{421} & \textbf{410} &
\textbf{3.1} & \textbf{2.2} \\
\hline
\end{tabular}
\label{tab:full_comparison}
\end{table}

\scriptsize{
\bibliographystyle{plain}
\bibliography{sample.bib}
}

\end{document}